\documentclass[11pt]{article}
\newcommand		{\comment}[1]		{}

\usepackage{amsmath,amsthm} 			
\usepackage{mathtools}

\usepackage[left=2.54cm,right=2.54cm,top=2.54cm,bottom=2.54cm]{geometry}
\usepackage[nouppercase]{scrlayer-scrpage}
\usepackage{xspace}
\usepackage{enumitem}

\usepackage{avant} 				
\usepackage[sc,osf]{mathpazo} 			
	\AtBeginDocument{
		\DeclareSymbolFont{AMSb}{U}{msb}{m}{n}
		\DeclareSymbolFontAlphabet{\mathbb}{AMSb}
	}
\usepackage{extpfeil} 
\usepackage{mathabx}  
\usepackage{mathrsfs}	
\usepackage{fancyvrb}	

\usepackage{rotating,pdflscape}
\usepackage{afterpage}
\usepackage{graphicx}
\usepackage{float}
\usepackage{subcaption}
\usepackage{caption}
\newcommand{\mockalph}[1]{\!}

\usepackage[utf8]{inputenc}
\usepackage[T1]{fontenc}

\usepackage{tocloft}
\makeatletter
\renewcommand{\l@figure}{\@dottedtocline{1}{1em}{3.5em}}
\renewcommand{\l@table}{\@dottedtocline{2}{1em}{3.5em}}
\makeatother
\newcommand*{\noaddvspace}{\renewcommand*{\addvspace}[1]{}}
\addtocontents{lof}{\protect\noaddvspace}

\usepackage[hyphens]{url}		%
\usepackage{hyperref}	
\usepackage{cleveref} 	
\newcommand		{\myred}		{BrickRed}
\hypersetup{
			colorlinks		= true, 	
			urlcolor		= \myred, 
			linkcolor		= \myred, 	
			citecolor		= PineGreen 	
}
\newcommand		{\hyref}[1]		{\hyperref[#1]{\ref*{#1}}}

\newif\ifdebug
\ifdebug\usepackage{lineno}\linenumbers\else\fi

\usepackage{etoolbox}
\PassOptionsToPackage{dvipsnames,svgnames,x11names,table}{xcolor}
\newcommand		{\defd}[1]	{\textcolor{RoyalBlue}{\textbf{\textit{#1}}}}
\newcommand		{\defm}[1]	{\textcolor{RoyalBlue}{#1}}

\usepackage[all]{xypic}\xyoption{rotate}
	\newdir{ >}{{}*!/-7pt/\dir{>}}
	\newdir{i}{{}*!/-7pt/\dir{(}	}
\usepackage{tikz}
\usetikzlibrary{matrix,fit,backgrounds,calc,arrows} 
\tikzstyle{image}=[rectangle,fill=Red!20,inner sep=-2pt]
\tikzstyle{nonzero}=[rectangle,fill=Navy!20,inner sep=0pt]
\tikzstyle{nonzerosm}=[rectangle,fill=Navy!20,inner sep=-2pt]

\usepackage{thmtools}
\makeatletter
\let\c@figure\c@table
\let\c@equation\c@table
\makeatother

\numberwithin{table}{section}
\numberwithin{figure}{section}
\newtheorem{abcthm}[table]{Theorem}

\newtheorem{theorem}[table]{Theorem}

\newtheorem{proposition}[table]{Proposition}

\newtheorem{corollary}[table]{Corollary}

\newtheorem{lemma}[table]{Lemma}

\newtheorem{claim}[table]{Claim}

\theoremstyle{definition}

\newtheorem{definition}[table]{Definition}

\newtheorem{notation}[table]{Notation}

\newtheorem{observation}[table]{Observation}

\newtheorem{conjecture}[table]{Conjecture}

\theoremstyle{remark}
\newtheorem{fact}[table]{Fact}

\newtheorem{example}[table]{Example}
\newtheorem{counterexample}[table]{Counterexample}
\newtheorem{exercise}[table]{Exercise}

\newtheorem{problem}[table]{Problem}
\newtheorem{histrmks}[table]{Historical remarks}
\newtheorem{remark}[table]{Remark}
\newtheorem{remarks}[table]{Remarks}

\theoremstyle{plain}
\newtheorem*{thm*}{Theorem}
\newtheorem*{theorem*}{Theorem}
\newtheorem*{prop*}{Proposition}
\newtheorem*{proposition*}{Proposition}
\newtheorem*{lemma*}{Lemma}
\newtheorem*{corollary*}{Corollary}
\newtheorem*{cor*}{Corollary}

\theoremstyle{definition}
\newtheorem*{definition*}{Definition}
\newtheorem*{defn*}{Definition}
\newtheorem*{QQ*}{Question}
\newtheorem*{obs*}{Observation}
\newtheorem*{notation*}{Notation}
\newtheorem*{discussion*}{Discussion}

\theoremstyle{remark}
\newtheorem*{rmk*}{Remark}
\newtheorem*{remark*}{Remark}
\newtheorem*{examples*}{Examples}
\newtheorem*{example*}{Example}
\newtheorem*{EG*}{Example}
\newtheorem*{EGs*}{Examples}
\newtheorem*{fact*}{Fact}
\newtheorem*{prob*}{Problem}

\newcommand{\bthm}{\begin{theorem}}
\newcommand{\ethm}{\end{theorem}}
\newcommand{\bprop}{\begin{proposition}}
\newcommand{\eprop}{\end{proposition}}
\newcommand{\bcor}{\begin{corollary}}
\newcommand{\ecor}{\end{corollary}}
\newcommand{\bconj}{\begin{conjecture}}
\newcommand{\econj}{\end{conjecture}}
\newcommand{\blem}{\begin{lemma}}
\newcommand{\elem}{\end{lemma}}
\newcommand{\bclm}{\begin{claim}}
\newcommand{\eclm}{\end{claim}}
\newcommand{\bpf}{\begin{proof}}
\newcommand{\epf}{\end{proof}}
\newcommand{\bdetails}{\begin{details}}
\newcommand{\edetails}{\end{details}}
\newcommand{\bdefi}{\begin{definition}}
\newcommand{\edefi}{\end{definition}}
\newcommand{\bdefn}{\begin{definition}}
\newcommand{\edefn}{\end{definition}}
\newcommand{\bex}{\begin{example}}
\newcommand{\eex}{\end{example}}
\newcommand{\bprob}{\begin{problem}}
\newcommand{\eprob}{\end{problem}}
\newcommand{\bob}{\begin{observation}}
\newcommand{\eob}{\end{observation}}
\newcommand{\bexer}{\begin{exercise}}
\newcommand{\eexer}{\end{exercise}}
\newcommand{\bexers}{\begin{exercises}}
\newcommand{\eexers}{\end{exercises}}
\newcommand{\brmk}{\begin{remark}}
\newcommand{\ermk}{\end{remark}}
\newcommand{\bhist}{\begin{histrmks}}
\newcommand{\ehist}{\end{histrmks}}
\newcommand{\brmks}{\begin{remarks}}
\newcommand{\ermks}{\end{remarks}}
\newcommand{\bverb}{\begin{verbatim}}
\newcommand{\everb}{\end{verbatim}}
\newcommand{\bntn}{\begin{notation}}
\newcommand{\entn}{\end{notation}}
\newcommand{\bfct}{\begin{fact}}
\newcommand{\efct}{\end{fact}}
\newcommand{\bfcts}{\begin{facts}}
\newcommand{\efcts}{\end{facts}}
\newcommand{\benum}{\begin{enumerate}}
\newcommand{\eenum}{\end{enumerate}}
\newcommand{\bitem}{\begin{itemize}}
\newcommand{\eitem}{\end{itemize}}

\tracingpatches
\makeatletter
\patchcmd{\@setref}{\bfseries ??}{\bfseries\color{red} FIX ME!}{}{}
\patchcmd{\@setcite}{\bfseries ?}{\bfseries\color{red} FIX ME!}{}{}
\patchcmd{\@setcref}         {??}{\color{red} FIX ME!}{}{}
\patchcmd{\@setcref}         {??}{\color{red} FIX ME!}{}{}
\patchcmd{\@setcrefrange}    {??}{\color{red} FIX ME!}{}{}
\patchcmd{\@setcrefrange}    {??}{\color{red} FIX ME!}{}{}
\patchcmd{\@setcrefrange}    {??}{\color{red} FIX ME!}{}{}
\patchcmd{\@setcrefrange}    {??}{\color{red} FIX ME!}{}{}
\patchcmd{\@setcrefrange}    {??}{\color{red} FIX ME!}{}{}
\patchcmd{\@setcrefrange}    {??}{\color{red} FIX ME!}{}{}
\patchcmd{\@setnamecref}     {??}{\color{red} FIX ME!}{}{}
\patchcmd{\@setnamecref}     {??}{\color{red} FIX ME!}{}{}
\patchcmd{\@setcpageref}     {??}{\color{red} FIX ME!}{}{}
\patchcmd{\@setcpageref}     {??}{\color{red} FIX ME!}{}{}
\patchcmd{\@setcpagerefrange}{??}{\color{red} FIX ME!}{}{}
\patchcmd{\@setcpagerefrange}{??}{\color{red} FIX ME!}{}{}
\patchcmd{\@setcpagerefrange}{??}{\color{red} FIX ME!}{}{}
\patchcmd{\@setcpagerefrange}{??}{\color{red} FIX ME!}{}{}
\patchcmd{\@setcpagerefrange}{??}{\color{red} FIX ME!}{}{}
\patchcmd{\@cref}            {??}{\color{red} FIX ME!}{}{}
\def\blx@citation@entry#1#2{%
  \blx@bibreq{#1}%
  \ifinlist{#1}{\blx@cites}
    {}
    {\listgadd{\blx@cites}{#1}%
     \blx@auxwrite\@mainaux{}{\string\abx@aux@cite{#1}}}%
  \ifinlistcs{#1}{blx@segm@\the\c@refsection @\the\c@refsegment}
    {}
    {\listcsgadd{blx@segm@\the\c@refsection @\the\c@refsegment}{#1}}%
  \blx@ifdata{#1}%
    {}%
    {\ifcsdef{blx@miss@\the\c@refsection}%
       {\ifinlistcs{#1}{blx@miss@\the\c@refsection}%
          {{\bfseries\color{red} cite:} }%
          {\blx@logreq@active{#2{#1}}}}%
       {\blx@logreq@active{#2{#1}}}}}
\def\blx@citeadd#1{%
  \ifcsdef{blx@keyalias@\the\c@refsection @#1}
    {\edef\blx@realkey{\csuse{blx@keyalias@\the\c@refsection @#1}}}
    {\def\blx@realkey{#1}}%
  \expandafter\blx@citation\expandafter{\blx@realkey}\blx@msg@cundefon
  \expandafter\blx@ifdata\expandafter{\blx@realkey}
    {\advance\blx@tempcnta\@ne
     \listeadd\blx@tempa{\blx@realkey}}
    {\ifnum\blx@tempcntb>\z@\multicitedelim\fi
     \expandafter\abx@missing\expandafter{\blx@realkey}%
     \advance\blx@tempcntb\@ne}}
\makeatother

\newcommand{\presectionskip}{-1.5\baselineskip}
\newcommand{\postsectionskip}{0.3\baselineskip}
\makeatletter
\usepackage{titlesec}
\renewcommand{\section}{\@startsection
  {chapter}{0}{0mm}
  {\presectionskip}
  {\postsectionskip}
  {\sffamily\huge}}
\renewcommand{\section}{\@startsection
  {section}{1}{0mm}
  {\presectionskip}
  {\postsectionskip}
  {\sffamily\LARGE}}
\renewcommand{\subsection}{\@startsection
  {subsection}{2}{0mm}
  {\presectionskip}
  {\postsectionskip}
  {\sffamily\Large}}
\renewcommand{\subsubsection}{\@startsection
  {subsubsection}{3}{0mm}
  {\presectionskip}
  {\postsectionskip}
  {\sffamily\normalsize}}
\renewcommand{\@seccntformat}[1]{\csname the#1\endcsname.\quad}
\newcommand\HUGE{\@setfontsize\Huge{30}{47}} 
\makeatother
\usepackage{titling}
  \titleformat{\chapter}[display]
  {\sffamily\Large}
  {Chapter {\HUGE\normalfont\thechapter}}    
  {1em}
  {\huge}
  \titleformat{\section}
  {\sffamily\Large}
  {\normalfont{\@setfontsize\Large{18}{47}\thesection.}}    
  {1em}
  {}
  \titleformat{\subsection}
  {\sffamily\large}
  {\normalfont{\Large\thesubsection}.}    
  {1em}
  {}
  \titleformat{\subsubsection}
  {\sffamily\normalsize}
  {\normalfont{\large\thesubsubsection}.}    
  {1em}
  {}

\renewcommand		{\SS}				{\textsection}

\newcommand		{\ang}[1]			{\langle #1 \rangle}

\newcommand		{\quation}[1]			{\begin{equation} #1 \end{equation}}

\newcommand		{\eqn}[1]			{\begin{align*} #1 \end{align*}}

\newcommand		{\eqnl}[1]			{\begin{equation}\begin{aligned}#1 \end{aligned}\end{equation}}

\def			\SPSB#1#2			{\rlap{\textsuperscript{#1}}\textsubscript{#2}}
\makeatletter
\def			\smallunderbrace#1		{\mathop{\vtop{\m@th\ialign{##\crcr
							   $\hfil\displaystyle{#1}\hfil$\crcr
							   \noalign{\kern3\p@\nointerlineskip}%
							   \tiny\upbracefill\crcr\noalign{\kern3\p@}}}}\limits}
\makeatother

\makeatletter
\newcommand{\subalign}[1]{%
  \vcenter{%
    \Let@ \restore@math@cr \default@tag
    \baselineskip\fontdimen10 \scriptfont\tw@
    \advance\baselineskip\fontdimen12 \scriptfont\tw@
    \lineskip\thr@@\fontdimen8 \scriptfont\thr@@
    \lineskiplimit\lineskip
    \ialign{\hfil$\m@th\scriptstyle##$&$\m@th\scriptstyle{}##$\crcr
      #1\crcr
    }%
  }
}
\makeatother

\makeatletter
\newcommand		{\oset}[3][0ex]			{%
								\raisebox{.175ex}{$%
								  \mathrel{\mathop{#3}\limits^{
								    \vbox to#1{\kern-2\ex@
								    \hbox{$\scriptstyle#2$}\vss}}}
								    $}%
							    }
\makeatother

\usepackage{faktor,xfrac}

\newcommand		{\nd}		{\noindent}

\newcommand		{\bs}		{\bigskip}

\newcommand		{\mn}		{\mspace{-2mu}}
\newcommand		{\mnn}		{\mspace{-1mu}}

\newcommand		{\ol}			{\overline}
\newcommand		{\os}			{\overset}
\newcommand		{\us}			{\underset}

\newcommand		{\ul}			{\underline}

\newcommand		{\wt}			{\widetilde}

\newcommand		{\mr}			{\mathrm}
\newcommand		{\bb}			{\mathbb}

\newcommand		{\ms}			{\mathscr}
\newcommand		{\sans}			{\mathsf}

\newcommand		{\g}		{\gamma}

\renewcommand		{\epsilon}	{\varepsilon}
\newcommand		{\e}		{\epsilon}
\newcommand		{\z}		{\zeta}
\newcommand		{\h}		{\eta}

\renewcommand		{\l}		{\lambda}

\newcommand		{\s}		{\sigma}

\newcommand		{\w}		{\omega}

\newcommand		{\W}		{\Omega}

\newcommand		{\D}		{\Delta}

\DeclareSymbolFont{cmletters}{OT1}{cmr}{m}{n}
\DeclareMathSymbol{\Ups}{\mathalpha}{cmletters}{"7}
\renewcommand		{\Upsilon}	{\Ups}

\newcommand		{\ceq}		{\coloneqq}

\renewcommand		{\cup}		{\mspace{-1mu}\smile\mspace{-1mu}}

\DeclareRobustCommand	{\lq}		{\text{\reflectbox{$/$}}}		

\makeatletter
\DeclarePairedDelimiterX
			{\pmodx}[1]	{(}{)}{{\operator@font mod}\mkern6mu#1}
					\renewcommand{\pmod}{%
					  \allowbreak
					  \if@display\mkern18mu\else\mkern8mu\fi
						  \pmodx
					}
\makeatother

\renewcommand		{\:}		{\colon}
\renewcommand		{\-}		{^{-1}}

\renewcommand		{\o}		{\circ}

\newcommand		{\adj}		{\dashv}

\ExplSyntaxOn
\NewDocumentEnvironment{adjunctions}{O{}}
{
	\cs_set_eq:cN {@arraycr} \farin_arraycr:
	\keys_set:nn { farin/adjunction } { #1 }
	\begin{array}
		{
			@{ \hspace { \dim_eval:n { \l_farin_left_shift_dim + \l_farin_padding_dim } } }
			r
			@{ {\farin_strut:} \l_farin_symbol_tl {} }
			l
			@{ \hspace { \dim_eval:n { \l_farin_right_shift_dim + \l_farin_padding_dim } } }
		}
	}
	{
	\end{array}
}
\keys_define:nn { farin/adjunction }
{
	leftshift       .dim_set:N = \l_farin_left_shift_dim,
	leftshift       .initial:n = 0pt,
	rightshift      .dim_set:N = \l_farin_right_shift_dim,
	rightshift      .initial:n = 0pt,
	padding         .dim_set:N = \l_farin_padding_dim,
	padding         .initial:n = 6pt,
	symbol          .tl_set:N  = \l_farin_symbol_tl,
	symbol          .initial:n = \longrightarrow,
	verticalspacing .dim_set:N  = \l_farin_vertspac_dim,
	verticalspacing .initial:n = {3pt},
}
\cs_new_protected:Npn \farin_strut:
{
	\vrule height \dim_eval:n { \ht\strutbox + 1.2\l_farin_vertspac_dim }
	depth  \dim_eval:n { \dp\strutbox + \l_farin_vertspac_dim }
	width 0pt
}
\makeatletter
\exp_args:NNo \cs_new:Npn \farin_arraycr:
{
	\@arraycr\hline
}
\makeatother
\ExplSyntaxOff

\renewcommand		{\.}		{\cdot}

\newcommand		{\x}		{\times}
\newcommand		{\xu}[3]	{\smash{{#2}\us{#1}\times{#3}}}

\newcommand		{\oplushigher}	{\mathbin{\raisebox{.85pt}{$\displaystyle\oplus$}}}

\DeclareMathOperator*	{\otimesvariable}{%
			\mathchoice {\raisebox{.85pt}{$\displaystyle\otimes$}}
						{\raisebox{.85pt}{$\otimes$}}
						{\raisebox{0.7pt}{$\scriptstyle\otimes$}}
						{\raisebox{0.2pt}{$\scriptscriptstyle\otimes$}}
						}
\newcommand		{\tensor}	{\otimesvariable}

\newcommand		{\direct}	{\oplushigher}
\newcommand		{\ox}		{\tensor}
\newcommand		{\ot}		{\direct}
\newcommand		{\+}		{\direct}

\newcommand		{\Direct}	{\bigoplus}

\newcommand		{\ext}		{\exterior}

\newcommand		{\bideg}	{\operatorname{bideg}}

\DeclareMathOperator	{\diag}		{diag}

\DeclareMathOperator	{\rk}		{rk }
\DeclareMathOperator	{\im}		{im }

\DeclareMathOperator	{\coker}	{coker }

\DeclareMathOperator	{\Tor}		{Tor}

\newcommand		{\cupone}	{\mathbin{{\cup}_1}}

\DeclareMathOperator	{\Map}		{Map}

\newcommand		{\U}		{\mr{U}}
\newcommand		{\SU}		{\mr{SU}}

\makeatletter
\newbox\xrat@below
\newbox\xrat@above
\newcommand		{\xrightarrowtail}[2][]	{%
						  \setbox\xrat@below=\hbox{\ensuremath{\scriptstyle #1}}%
						  \setbox\xrat@above=\hbox{\ensuremath{\scriptstyle #2}}%
				  \pgfmathsetlengthmacro{\xrat@len}{max(\wd\xrat@below,\wd\xrat@above)+.6em}%
  						\mathrel{\tikz [>->,baseline=-.55ex]
              					   \draw (0,0) -- node[below=-2pt] {\box\xrat@below}
                            					    node[above=-2pt] {\box\xrat@above}
                    						   (\xrat@len,0) ;}
						}
\newbox\xrat@below
\newbox\xrat@above
\renewcommand		{\xtwoheadrightarrow}[2][]{%
						  \setbox\xrat@below=\hbox{\ensuremath{\scriptstyle #1}}%
						  \setbox\xrat@above=\hbox{\ensuremath{\scriptstyle #2}}%
				  \pgfmathsetlengthmacro{\xrat@len}{max(\wd\xrat@below,\wd\xrat@above)+.6em}%
						 \mathrel{\tikz [->>,baseline=-.55ex]
					                 \draw (0,0) -- node[below=-2pt] {\box\xrat@below}
					                                node[above=-2pt] {\box\xrat@above}
						                       (\xrat@len,0) ;}
		       				}
\makeatother

\newcommand		{\xepi}		{\xtwoheadrightarrow}

\newcommand		{\longepi}	{\xepi[]{\ \ \ \ }}

\newcommand		{\longto} 	{\longrightarrow}
\newcommand		{\lt}		{\longto}
\newcommand		{\xtoo}		{\xrightarrow} 
\newcommand		{\from}		{\leftarrow}
\newcommand		{\longfrom}	{\longleftarrow}
\newcommand		{\lf}		{\longfrom}

\newcommand		{\lmt}		{\longmapsto}

\newcommand		{\simto}	{\xrightarrow{\sim}}

\newcommand		{\longsimto}	{\os\sim\longto}
\newcommand		{\isoto}	{\longsimto}

\newcommand		{\longbij}	{\longleftrightarrow}

\newcommand		{\vertsim}	{\rotatebox{90}{$\sim$}}

\newcommand		{\hmt}		{\simeq}
\newcommand		{\iso}		{\cong}
\newcommand		{\homeo}	{\approx}

\newcommand		{\F}		{\bb F}

\newcommand		{\Z}		{\bb Z}

\newcommand		{\R}		{\bb R}
\newcommand		{\C}		{\bb C}

\newcommand		{\RP}		{\bb R \mr P}

\newcommand		{\RPi}		{\RP^\infty}

\DeclareMathOperator	{\id}		{id}

\DeclareMathOperator	{\Sq}		{Sq}

\newcommand		{\CGA}		{\textsc{cga}\xspace}

\newcommand		{\DGA}		{\textsc{dga}\xspace}
\newcommand		{\DG}		{\textsc{dg}\xspace}
\newcommand		{\DGAs}		{\textsc{dga}s\xspace}
\newcommand		{\DGC}		{\textsc{dgc}\xspace}
\newcommand		{\DGCs}		{\textsc{dgc}s\xspace}
\newcommand		{\CDGA}		{\textsc{cdga}\xspace}
\newcommand		{\CDGAs}	{\textsc{cdga}s\xspace}

\newcommand		{\SHC}		{\textsc{shc}\xspace}
\newcommand		{\Ai}		{$A_\infty$}

\newcommand		{\kk}		{k}

\renewcommand 		{\H}		{H^*}

\newcommand 		{\HK}		{{\H_K}}

\newcommand		{\EMSS}		{Eilenberg--Moore spectral sequence\xspace}

\newcommand		{\quism}	{quasi-isomorphism\xspace}
  \debugfalse
\usepackage{tracklang}
\usepackage{multirow,accents}
\usepackage{array}\newcolumntype{R}{>{$}l<{$}}
\usepackage{scalerel}[2016/12/29]
\usepackage{chngcntr}
\counterwithin{table}{section}
\counterwithin{equation}{section}

\theoremstyle{definition}
\newtheorem*{layout*}{Outline}
\newtheorem{sublayout}{}[table]

\usepackage{silence}
\colorlet{jlabel}{Crimson!75!Black}		
\colorlet{jred}{Crimson!90!Black}	
\colorlet{nonisored}{Crimson!90!Black}	
\colorlet{jhmt}{Dandelion!80!Indigo}	
\colorlet{jcmp}{Black!35!White}
\colorlet{jblue}{RoyalBlue!80!Black}	
\colorlet{jviolet}{Violet!60!Black}
\colorlet{compromise}{rgb:Magenta!85!Indigo,3;white,3;Black,2}	
\colorlet{jQ2}{rgb:Indigo!90!Black,2;red,2}	
\colorlet{jA}{RoyalBlue!40!Indigo}		
\colorlet{jX3}{BlueViolet!100!Black}
\colorlet{jX2}{BlueViolet!100!Black}
\colorlet{jX2p}{BrickRed!60!Dandelion!50!Magenta!60!White}	
\colorlet{jgreen}{SeaGreen!90!Black}
\colorlet{jcyan}{Indigo!30!Cyan!95!White}

\xyoption{line}

\let		\epsilon	\varepsilon
\newcommand		{\B}		{\mathbf{B}}
\renewcommand	{\W}		{\BW}
\newcommand		{\Y}		{\Upsilon}
\renewcommand	{\ot}		{^{\ox \mnn 2}}
\newcommand		{\ott}		{^{\ox \mnn 3}}
\renewcommand	{\2}		{\ot}
\newcommand		{\hd}		{h^\dagger}
\newcommand		{\desusp}		{s^{-1}}
\def				\iter#1#2{#1^{[#2]}}
\newcommand		{\muA}		{\mu_{\mn A}}
\renewcommand	{\C}		{C^*}
\newcommand		{\I}		{I^*}
\newcommand		{\Dpra}		{\pi_0 \x 0}
\newcommand		{\Dprb}		{\pi_1 \x 0 = 0 \x \pi_0}
\newcommand		{\Dprc}		{0 \x \pi_1}
\newcommand		{\Miso}		{\Xi}

\newcommand		{\Algs}		{\sans{DGA}}
\newcommand		{\Coalgs}		{\sans{DGC}}
\newcommand		{\GM}		{\sans{Mod}}

\newcommand		{\Ysub}		{}
\renewcommand	{\EMSS}		{\textsc{emss}\xspace}
\newcommand		{\SHCA}		{\SHC-algebra\xspace}
\newcommand		{\SHCAs}		{{\SHCA}s\xspace}
\newcommand		{\WHC}		{\textsc{whc}\xspace}
\newcommand		{\WHCA}		{\WHC-algebra\xspace}

\renewcommand{\ext}{\mr{ext}}

\DeclareMathOperator*{\T}{%
	\mathchoice{\raisebox{.85pt}{$\displaystyle\underline{{\otimes}}$}}
	{\mathbin{{\raisebox{.85pt}{$\underline{{\otimes}}$}}}}
	{\mathbin{{\raisebox{.7pt}{$\scriptstyle\underline{{\otimes}}$}}}}
	{\mathbin{{\raisebox{.2pt}{$\scriptscriptstyle\underline{{\otimes}}$}}}}
}

\DeclareFontFamily{U}{wncy}{}
\DeclareFontShape{U}{wncy}{m}{n}{<->wncyr10}{}
\DeclareSymbolFont{mcy}{U}{wncy}{m}{n}
\DeclareMathSymbol{\Sha}{\mathord}{mcy}{"58}

\def				\BW			{\boldsymbol{\Omega}}

\newtheorem{manualeginner}{Example}
\newenvironment{manualeg}[1]{%
  \IfBlankTF{#1}
    {\renewcommand{\themanualeginner}{\unskip}}
    {\renewcommand\themanualeginner{#1}}%
  \manualeginner
}{\endmanualeginner}

\newtheorem{manualpropinner}{Proposition}
\newenvironment{manualprop}[1]{%
  \IfBlankTF{#1}
    {\renewcommand{\themanualpropinner}{\unskip}}
    {\renewcommand\themanualpropinner{#1}}%
  \manualpropinner
}{\endmanualpropinner}

\newtheorem{manualrmkinner}{Remark}
\newenvironment{manualrmk}[1]{%
  \IfBlankTF{#1}
    {\renewcommand{\themanualrmkinner}{\unskip}}
    {\renewcommand\themanualrmkinner{#1}}%
  \manualrmkinner
}{\endmanualrmkinner}

\newcommand{\Tr}[1]{	\us{#1}	{\mathrm{Tor}}	}

\newcommand{\BA}{\B A}
\newcommand{\BX}{\B X}
\newcommand{\BY}{\B Y}

\newcommand		{\new}[1]		{\nd\textcolor{JungleGreen}{#1}}
\newcommand		{\wrong}[1]		{\nd\textcolor{Red}{[#1]}}

\DeclareMathVersion{normal2}

\begin{document}
\title{A ring structure on Tor}
\author{Jeffrey D. Carlson}



\maketitle

\abstract{{\footnotesize 
We prove that within a natural class of $E_3$-algebras,
	the graded group $\Direct \Tor^i_A(X,Y)$
	induced by a pair of $E_3$-algebra maps $X \from A \to Y$
	carries a 
	graded algebra structure
	generalizing the classical structure
	when $A$, $X$, $Y$ 
	are genuine commutative differential graded algebras. 

	\smallskip
	
	We \new{attempt to} prove, as a topological corollary,
	that Munkholm's Eilenberg--Moore collapse
	result for pullbacks of spaces with polynomial cohomology
	can be enhanced to a ring isomorphism.
	\new{This is not achieved, 
	and in fact the claim as stated in the previous drafts is false.
	If additionally, $2$ is assumed to be a unit of the base ring,
	then that claim is true (not that the results in this paper establish it)
	and is known due to previous work of the author
	and Franz, and also, as it turns out,
	to Huebschmann's unpublished 1983 habilitation work.}
}
}


\vspace*{14pt}

\section*{Caveat lector}
An earlier draft of this paper passed refereeing,
but it was afterwards pointed out to the author by Matthias Franz
that the main topological result, \Cref{thm:main-cor-topological},
as stated in the previous draft, 
is false
(without an additional hypothesis which updates its status to ``known'').
The main algebraic result, \Cref{thm:main-algebraic},
is true subject to additional conditions not observed in the previous draft, 
but these are hard to verify, and 
are definitely untrue 
(see \Cref{eg:Baum})
in the intended use case, that of the main topological result---%
unless $2$ is a unit in the base ring,
the hypothesis which the entire motivation of this paper was originally to avoid.
When $2$ \emph{is} a unit, 
the main topological result is indeed true,
as the author belatedly found, in writing a literature review, 
was already established long ago in unpublished work 
of Huebschmann~\cite[Cor.~6.2, p.~78]{huebschmann1983perturbation},
and had more recently and independently been shown 
in work of the author and Franz~\cite{carlsonfranzlong};
however the results of the present paper do \emph{not} provide yet a third proof.

To communicate what is true within this document
without spending more time than necessary revising it
until and unless it can be genuinely repaired,
we adopt a highlighting scheme in which \new{additions} are made in 
\texttt{xcolor}'s ``{JungleGreen}''
and \wrong{falsehoods from the last draft} are reddened and bracketed.
We will call particular attention to the error, 
delineate the portions of the paper it contaminates,
and indicate what corrections need to be made to salvage the rest.
We suspect it is possible to prove 
  the additional compatibility conditions on homotopies
    needed to make \Cref{thm:main-algebraic} hold
  do obtain in the situation of \Cref{thm:main-topological}
if we additionally assume $2$ is a unit,
but actually verifying as much would seem to require 
explicit computations the author is not prepared or probably equipped to make.

For anyone keeping score, 
the natural guess on scanning this paper would have been that 
its central error is in claiming the commutativity 
of one of the margin-straining diagrams toward the end.
These are actually fine;
the issue is with a much simpler diagram toward the beginning.

\bigskip

\nd \emph{Acknowledgments}.
The author's interest in this question 
originated in joint work with Matthias Franz.
The author sketched a proof 
generalizing Franz's homogeneous space result~\cite{franz2019homogeneous} 
to biquotients and asked what was left to do.
Franz said the missing ingredient was a natural product structure on 
a two-sided bar construction of homotopy Gerstenhaber algebras
and within a month gave a mod-two reduction of what he was convinced
were the correct formulas. 
The author then found the right signs and proved the formulas worked. 
A draft has by now appeared on the {arXiv}~\cite{carlsonfranzlong} 
and a briefer version will appear in \emph{Algebraic \& Geometric Topology}.

The author would like to thank Omar Antol\'in Camarena for listening 
to the ongoing drama of this proof 
and for detailed explanations of several sections of Lurie's \emph{Higher Algebra}
back when it seemed that approach might be fruitful,
Manuel Rivera for discussion of Husemoller--Stasheff--Moore,
Jim Stasheff for real magnanimity in discussion of his paper,
Peter May for discussing the history of his papers and requesting a sample computation,
Bernhard Keller for directing the author to Jesse Burke,
Jesse Burke for preventing him from spending any more time 
trying to find a model category structure on coalgebras over a general ring,
Bj\"orn Eurenius and David White for showing him such a structure 
does exist,
Markus Szymik for sending him a copy of Stasheff--Halperin~\cite{halperinstasheff1970},
Larry Smith and Pedro Tamaroff for detailed comments on an earlier version of this paper,
Richard Thomas, John Nicholson, 
and an anonymous, very thorough referee
for advice on presentation, and
Joanna Quigley for rigorous copyediting of an earlier draft.
\new{The author would like additionally, now, to
thank Matthias Franz for pointing out Baum's counterexample 
to the original statement of \Cref{thm:main-topological},
Anja Randecker for going to the library at Heidelberg 
to scan enough
of Johannes Huebschmann's 1983 \emph{Habilitationsschrift}
to determine the sought-after result was indeed there,
and Huebschmann for correspondence 
and for promising to revise and publish his results.}

\tableofcontents

\section*{Introduction}

Among the most fundamental objects of
homological algebra are the 
derived functors $\Tor^i$ of the tensor product.
These are individually merely modules,
but if $M \from A \to N$ are maps of 
commutative differential graded algebras (\defm{\CDGA}s),
then the graded group $\Direct_i \Tor^i_A(M,N)$ carries
the structure of a bigraded ring denoted simply $\defm{\Tor_A(M,N)}$,
essentially because under the hypothesis of commutativity,
the multiplications
$A \ox A \lt A$
and so on can themselves be seen as maps
of differential graded algebras (\defm{\DGA}s).
When the input rings are instead cochain algebras 
$\C(X) \from \C(B) \to \C(E)$,
it is again classical that there exists a ring structure on Tor,
but this is 
because of the Eilenberg--Zilber theorem,
on the face of it an entirely unrelated reason.
It is thus natural to wonder under what 
general conditions on the input \DGAs
a ring structure on Tor should exist.
The question, however, seems never to have been seriously considered.

This situation is the more surprising because 
there does exist a candidate for such a product. 
In 1974, Hans J. Munkholm realized the products 
on $\Tor_{\H B}(\H X, \H E)$
and $\Tor_{\C B}(\C X, \C E)$
could both be described in terms of a 
structure generalizing \CDGAs,
called \emph{strongly homotopy commutative (\defm{\SHC}-) algebras},
of which 
both the cochain algebra $\C(-)$ 
and cohomology ring $\H(-)$
are examples.%
\footnote{\ 
		This language is now unfamiliar,
		but should not be intimidating:
		we will be able to state all the relevant facts
		about \SHCAs in less than two pages in \Cref{sec:SHC}.
		Although they are not strictly comparable, 
		the main examples show
		they should be thought of as somewhere
		between $E_2$- and $E_3$-algebras;
		see \Cref{thm:Franz-SHC} and \Cref{rmk:operads}.
		}
He had used these \SHCA structures to obtain a 
collapse result
for certain Eilenberg--Moore spectral sequences (\defm{\EMSS}s),
and went on in the last section
of his paper
to note, almost as an afterthought, 
that with some additional homotopy data,
this alternative construction produces a product on Tor 
more generally%
~\cite[\SS9]{munkholm1974emss}.
The construction is involved and rather speculative, and
Munkholm was not sanguine about his product's prospects~\cite[p.~49]{munkholm1974emss}:

%
%
%

\begin{quote}{\footnotesize
	The composition of (1), (2) and (3) now gives some sort of a product. 
	We have no specific applications of this in mind, 
	so we have not tried to investigate the properties of this product. 
	Presumably they are relatively bad, 
	because of the dependence of (1) on our choice of the homotopies.
	}
\end{quote}

\nd 
\wrong{But he was wrong to doubt.}\footnote{\ 
	Maddeningly, there is a 1976 reference~\cite{munkholm1976shmII} 
	to an unpublished paper with the word ``multiplicative''
	in its title~\cite{munkholm1976shmIII}.
	It is 
	hard not to wonder what Munkholm
	would have said or why this paper never appeared.
}
\new{He was absolutely right to speculate that the this product would only be 
functorial relative to the choice of homotopies defining it,
and that this would be a significant impediment to using it.}
It is the main task of the present paper to vindicate this product.\footnote{\ 
	More accurately speaking, in this manuscript we construct and 
	valorize our own product;
	to show it agrees with Munkholm's is in some ways
	more technical and can wait until a sequel.
	\new{There does not seem to be anything wrong with this other manuscript,
	but there does not seem to be any great urgency 
	to its publication, either.}
}

\begin{theorem}\label{thm:Tor-CGA}
	Let {\SHCA}s $A$, $X$, $Y$, and \SHCA maps 
	$\B X \from \B A \to \B Y$ be given.
	Then the product 
	defined in \Cref{def:product}
	is independent of the choice of homotopies used to define it
	and renders $\Tor_A(X,Y)$
	a commutative graded algebra.
	Moreover, the product is functorial in triples of \SHCA maps
	making the two necessary squares commute up to homotopy.
\end{theorem}


	There are other plausible approaches 
	to defining a product on Tor
	running through highly structured ring spectra 
	(see \Cref{rmk:Cinfty} for sketches of how this might work),
	but the the  
	advantage of our approach 
	is that it obtains functoriality 
	\new{(of a sort---which is not as strong as one would want)}
	of the product under minimal hypotheses on maps of input \DGAs.
	This is critical for us because these hypotheses
	are already known to hold in the motivating situation,
	Munkholm's original Eilenberg--Moore collapse result,
	whereas the stronger hypotheses required for the functoriality
	of other probable products are not.

Casting our minds back, the advent of the \EMSS%
~\cite{EMSS1965,smith1967emss} 
had made 
the cohomology of a wide range of fiber products 
accessible to computation, 
loop spaces and two-stage Postnikov systems being popular examples. 
Cartan~\cite{cartan1950transgression}
had famously shown that for coefficient ring $\kk = \R$
and $K < G$ compact, connected Lie groups, one has
\quation{\label{eq:G/H}
	\H(G/K) \iso \Tor_{\H BG}(\kk, \H BK)\mathrlap,
}
and with a view especially toward generalizing this
result to harder coefficient rings,
authors including Baum, May, Gugenheim, Munkholm, Halperin, Stasheff, Husemoller, Moore, and Wolf---with eventual success---set themselves the task
of proving collapse results,
if not for the \EMSS of a general pullback,
then at least for the one-sided variant applying
to a fibration of the form
$F \to E \to B$ with $\H(E)$ and $\H(B)$ polynomial
and converging to $\H(F)$.

Munkholm's 1974 result was 
the most far-reaching of these, 
not only showing that the sequence collapses, 
but resolving the additive extension problem.
There was no general multiplicative result until 2019,\footnote{\
	excepting the special case
	$G/H \to BH \to BG$ 
	when $\rk G = \rk H$ and $\kk$ is a field,
	which Borel and Baum proved in different ways in their theses~%
	\cite[Prop.~3.2]{borelthesis}%
	\cite[Cor.~7.5]{baum1968homogeneous}.
	}
when Franz~\cite{franz2019homogeneous} 
established multiplicativity in the ``one-sided'' case,
meaning  \eqref{eq:G/H} holds
as an isomorphism of graded rings when $\H(BG;\kk)$ and $\H(BK;\kk)$
are polynomial and $2$ is a unit of $\kk$.
While this represented the greatest progress on the problem 
in the forty-five years since Munkholm,
in the two-sided case,
Munkholm's additive isomorphism 
was not still known to be multiplicative.

It is our second central observation 
that it is
(and one need not invert~$2$).
Since many existing \EMSS collapse results factor through Munkholm's,
it follows 
that 
we have ring isomorphisms 
much more commonly than it
had seemed reasonable to hope.
This is Munkholm's result:

\begin{theorem}[Munkholm]
\label{thm:main-additive}
	Let $X \to B \from E$ be a diagram of topological spaces 
	with $E \lt B$ a Serre fibration 
	such that $\pi_1(B)$ acts trivially on $\H(E;\kk)$
	and suppose that $\H(X;\kk)$, $\H(B;\kk)$, and $\H(E;\kk)$ 
	are polynomial rings on at most countably many generators.
	If the characteristic of the principal ideal domain $\kk$ is $2$, 
	assume as well that the $\cupone$-square vanishes 
	on some selection of polynomial generators for $\H(X;\kk)$
	and $\H(E;\kk)$.
	Then there is a graded {$\kk$-module} isomorphism
	\begin{equation}\label{eq:iso}
		\Tor_{\H(B;\kk)}\big(\mnn\H(X;\kk),\H(E;\kk)\mnn\big) 
			\isoto 
		\H(\xu B X E;\kk)
		\mathrlap.
	\end{equation}
\end{theorem}

\begin{corollary}
	In the situation 
of \Cref{thm:main-additive},
the	Eilenberg--Moore spectral sequence
of $X \to B \from E$
	collapses 
	with no additive extension problem.
\end{corollary}	


\wrong{Our}
\new{An} enhancement 
\wrong{assumes 
only what Munkholm does, 
except in characteristic $2$,
and}
concludes the isomorphism is multiplicative:

\begin{theorem}[\new{Huebschmann~\cite[Cor.~6.2]{huebschmann1983perturbation}}]
\label{thm:main-topological}
	Assume the hypotheses of \Cref{thm:main-additive} and moreover,
	\wrong{if $\mr{char}\, \kk = 2$,
	that the $\cupone$-square vanishes 
	on some selection of polynomial generators for $\H(B;\kk)$.}
	\new{that $2$ is a unit of $\kk$}.
	Then
	\eqref{eq:iso} is a $\kk$-algebra isomorphism.
\end{theorem}

\begin{corollary}
	In the situation 
of \Cref{thm:main-topological},
the	Eilenberg--Moore spectral sequence
of $X \to B \from E$
	collapses 
	with no additive or multiplicative extension problem.
\end{corollary}	

%
%

\subsection*{Applications of \Cref{thm:main-topological}}
The apparently restrictive hypothesis that the input spaces have polynomial 
cohomology holds for the motivating classical example of a homogeneous space 
$G/K$ as in \eqref{eq:G/H}. In fact, we can substantially generalize this.

As noted by Singhof~\cite{singhof1993},
if $H$ and $K$ are subgroups of a topological group $G$,
then the homotopy orbit space \[EK \ox_K G/H \ceq \frac{EK \x G/H}{(ek,gH)\sim(e,kgH)\mathrlap{, \quad e \in EG, \ k \in K, \ g \in G}}\]
is the homotopy pullback of $BK \to BG \from BH$,
realizable as the pullback of the $G/H$-bundle $EG/H \lt EG/G$ for a fixed
model for $EG$.
If the two-sided action of $K \x H$ on $G$ by $(k,h)\.g \ceq kgh\-$
is free, this can be identified with the so-called \defd{biquotient}
$K\lq G/H$, which if $G$ is a Lie group and~$K,H < G$ are closed
is naturally a smooth manifold.
(For $K = 1$, of course, this reduces to the statement that~$G/H$ 
is the fiber of $EG/H \lt EG/G$.)
Then if $K$ and $H$ are connected, so that after inverting finitely 
many primes $\H(BK)$, $\H(BG)$, and $\H(BH)$ are all polynomial rings,
\Cref{thm:main-topological} immediately gives the
Borel cohomology ring $\HK(G/H) \ceq \H(EK \ox_K G/H)$ with suitable coefficients.%
\footnote{\ It can happen that $\H(BK)$ and $\H(BH)$ are polynomial 
	even if $K$ or $H$ is disconnected,
	for example if the component group $K/K_0$ 
	acts as a reflection group on $\H(BK_0)$
	and $\kk$ is a field of characteristic relatively prime
	to~$|K/K_0|$.
	}

	\begin{theorem}\label{thm:main-cor-topological}
	Let $G$ be a compact, connected Lie group, $K$ and $H$ closed subgroups,
	and $\kk$ a principal ideal domain \new{in which $2$ is unit},
	over which the cohomology of the classifying spaces
	$BG$, $BH$, and $BK$ is polynomial.
	Then
	we have 
	an isomorphism of graded $\kk$-algebras 
	\[\HK(G/H) 
	\isoto
	\Tor_{\H BG}(\H BK, \H BH)
	\mathrlap.\]
	In particular,
	this determines $\H(K\lq G/ H)$ if $K \x H$ acts freely 
	on $G$ under the standard action
	and $\H(G/H)$ if $K=1$.
\end{theorem}

\wrong{
This improves even on the result 
of Franz~\cite{franz2019homogeneous} for $G/H$,
in that $2$ no longer need be a unit of $\kk$.}
\new{As noted above, work of Huebschmann and later of Carlson--Franz already 
    establishes this, and the present paper does not.
}
\wrong{%
A minimal example in which one recovers previously unavailable torsion 
is the following.}

\begin{example}
	Let $H$ be the $\U(1)$ subgroup of $\SU(4)$ with diagonal entries ${\diag}(z^{-3},z,z,z)$.
	One previously knew from Franz's result~\cite{franz2019homogeneous}
	that,
	indexing generators by degree,
	\[
	\H\Big(\SU(4)/H;\Z\big[{\textstyle\frac 1 2}\big]\Big) \,\iso\,
		\frac{\Z\big[\frac 1 2\big][s_2] \otimes \Lambda[a_5,b_7]}
	{(3s^2,s^3,s^2 a)}\mathrlap.
	\]
	\wrong{
	Now one sees that in fact
	\eqn{
	\H\big(\SU(4)/H;\Z\big) \,&\iso\, \frac{\Z[s_2] \otimes \Lambda[y_5,z_7]}
	{(6s^2,2s^3,s^4,2s^2 y,3s^2 z)}\mathrlap,\\
	\H\big(\SU(4)/H;\F_2\big) \,&\iso\, \frac{\F_2[s_2]}{(s^4)} 
	\otimes \Lambda[x_3,y_5].%
	}
	} 
	\new{The preceding could be true, 
	but known results do not actually help us prove it, if it is.}
\end{example}

\begin{manualeg}{0.9a}[{\cite[p.~38]{baum1968homogeneous}}]\label{eg:Baum}
\new{
	The \EMSS with input $\ast \to B\U(2) \from BZ\U(2)$,
	converging to the cohomology ring 
	of $\U(2)/Z\U(2) = P\U(2) \homeo \RP^3$,
	starts and ends with $E_2 = \Tor_{\kk[c_1,c_2]}\big(\kk,\kk[s]\big)$,
	where $|s| = |c_1| = 2$ and $ |c_2| = 4$,
	with $c_1 \lmt 2s$ and $c_2 \lmt s^2$.
	One can compute the Tor as the cohomology of the Koszul complex
		$\kk[s]\ \ox \ \Lambda[z_1,w_3]$,
	where $\bideg s = (0,2)$, 
		$\bideg z_1 = (-1,2)$, 
		$\bideg w_3 = (-1,4)$.
	For $\kk = \F_2$
	this works out to $\F_2[s]/(s^2) \ox \Lambda[z_1]$,
	whereas $\H(\RP^3;\F_2)$ is the nonisomorphic ring $\F_2[z_1]/(z_1^4)$.
	}

\new
{This already appears, embarrassingly,
	in a 1968 paper (and the original 1962 dissertation~\cite[p.~3.42]{baumthesis}
	it derives from)
	that the present author read in 2014,
	and shows that \Cref{thm:main-cor-topological}
	is not true as originally stated.}
\end{manualeg}

Another obvious application is to the cohomology of a free loop space $LX$,
immediately recovering the main result of a paper of Saneblidze~\cite{saneblidze2009bitwisted}:

\begin{theorem}[Saneblidze]\label{thm:free-loop}
Let $X$ be a space and $\kk$ a principal ideal domain such that $\H X$ is polynomial
on an at most countable set $Q$ of generators whose $\cupone$-squares vanish. 
Then
we have ring isomorphisms
\[
	\H LX 	\iso \Tor_{\H(X\x X)}(\H X, \H X)
			\iso \H X \ox \Lambda[\desusp Q] \iso \H X \ox \H \Omega X\mathrlap,
\]
where $\desusp Q$ is a set of generators $\desusp q$, for $q \in Q$,
with degrees $|\desusp q| = |q| - 1$. 
\end{theorem}
\begin{proof}It is well known 
that the free loop space $LX = \Map(S^1,X)$
can be identified with the homotopy pullback of the diagonal map
$\D\: X \lt X \x X$ along itself. 
\end{proof}

The third isomorphism in \Cref{thm:free-loop} implicitly
used the following result on a based loop space $\Omega B$,
which also follows from \Cref{thm:main-topological}:

\begin{theorem}[Probably Borel]\label{thm:loop}
Let $B$ be a space and $\kk$ a principal ideal domain such that $\H B$ is polynomial
on an at most countable set $Q$ of generators whose $\cupone$-squares vanish. 
Then
we have ring isomorphisms
\[
	\H \Omega B 	\iso \Tor_{\H (B)}(\kk,\kk)
			\iso \Lambda[\desusp Q] \mathrlap,
\]
where $\desusp Q$ is a the set of generators $\desusp q$, for $q \in Q$,
with degrees $|\desusp q| = |q| - 1$. 
\end{theorem}
\begin{proof}
In this case, we use the fact $\Omega B$
is the homotopy pullback of ${\ast} \to B \from {\ast}$.
\end{proof}

\begin{counterexample}\label{ex:RPi}
	We really need the added hypothesis on $\H(B)$ for \Cref{thm:main-topological}
	to go through.
	To see that \Cref{thm:loop} fails without this hypothesis,
note that $B^2(\Z/2) = K(\Z/2,2)$ 
does \emph{not} satisfy the hypotheses over $\kk = \F_2$,
its cohomology being the polynomial ring generated by 
the iterated $\cupone$-squares 
$x_{1+2^\ell} = \Sq^{2^\ell}\Sq^{2^{\ell-1}}\cdots\,\Sq^2 \Sq^1 \iota_2 
\in H^{1+2^\ell}K(\Z/2,2)$
of the fundamental class $\iota_2 \in H^2 K(\Z/2,2)$.
Thus, although Munkholm's theorem holds for the 
homotopy pullback $\RPi = K(\Z/2,1) = \Omega K(\Z/2,2)$
of~${*} \from K(\Z/2,2) \to {*}$,
affording us an isomorphism
\[
\F_2[\iota_1] = \H K(\Z/2,1) \iso 
\Tor_{\H K(\Z/2,2)}(\F_2,\F_2) \iso \Lambda_{\F_2}[\desusp \iota_2,
							\desusp x_3,
							\desusp x_5,
							\desusp x_9,\ldots]
\]
of graded vector spaces, this isomorphism is not multiplicative.
For more on loop spaces, see \Cref{rmk:exterior}.
See Saneblidze~\cite{saneblidze2017loop} 
for a detailed account of what can happen for $\H(\Omega X;\F_2)$
when $\cupone$-squares do not vanish.
\end{counterexample}

\Cref{thm:main-topological} also recovers the easiest cases of group cohomology:

\begin{proposition}[Classical]
	Let a finitely-generated abelian group $A$ 
	and principal ideal domain $\kk$ be given.
	Let $0 \to R \to F \to A \to 0$ be a presentation of $A$,
	which is to say a short exact sequence of groups
	with $R$ and $F$ free abelian.
	Then 
	we have a ring isomorphism
	\[
	\H(BA;\kk) \iso \Tor_{\H(B^2 F;\kk)}\big(\kk,\H(B^2 R;\kk)\big)\mathrlap.
	\]
\end{proposition}
\bpf
	There is an evident fiber sequence $BA \to B^2 R \to B^2 F$
	of Eilenberg--Mac Lane spaces, and~$\H B^2 R$ and $\H B^2 F$
	are polynomial rings on generators of degree $2$,
	hence \Cref{thm:main-topological} applies
	with $X = {*}$ and $(E \to B) = (B^2 R \to B^2 F)$.
\epf

\subsection*{Outline}
\stepcounter{table}
	The plan of the work is as follows. 
\begin{sublayout}
\Cref{sec:CDGA} defines algebras and coalgebras and the bar--cobar adjunction,
as well as the intermediary notion of a \emph{twisting cochain}.
\end{sublayout}
\begin{sublayout}
\Cref{sec:tensor} brings in the tensor product 
and some of its interactions with the adjunction.
\end{sublayout}
\begin{sublayout}
\Cref{sec:homotopy} discusses how homotopies of algebra maps 
can themselves be realized by algebra maps.
Particularly, \Cref{sec:path} introduces the {path object}
receiving such homotopies.
This material is classical, dating back to Munkholm's work or earlier,
until the
critical new \Cref{thm:hmt-P} and \Cref{thm:homotopies-homotopic}.
\Cref{sec:concatenation}, which is also new, 
then develops an array of categorical machinery allowing us to manipulate 
homotopies diagrammatically without leaving the category of algebras.
\end{sublayout}
\begin{sublayout}
\Cref{sec:Tor} discusses conditions under which maps on Tor 
of a span of \DGAs can be defined,
which are classical, and establishes their homotopy-invariance 
and functoriality, which are new. 
\end{sublayout}
\begin{sublayout}
\Cref{sec:SHC} recalls the notion of an \SHC-algebra.
\end{sublayout}
\begin{sublayout}
\Cref{sec:product} motivates and defines our reformulation of Munkholm's product.
\end{sublayout}
\begin{sublayout}
\Cref{sec:CGA} establishes a \CGA structure on Tor,
the first clause of \Cref{thm:Tor-CGA}\new{,
subject to additional conditions vastly restricting feasible use cases}.
\end{sublayout}
\begin{sublayout}
\Cref{sec:ring-map} proves this algebra structure is functorial in the input data
(the second clause of \Cref{thm:Tor-CGA})
\new{under restrictions drastically limiting applicability}
\wrong{and homotopy-invariant}.
\wrong{In particular, it proves \Cref{thm:main-topological}.}
\end{sublayout}

\nd These two sections involve some micromanagement of \DGA
homotopies and a number of diagrams,
but mostly rely on the formal properties of the homotopy categories of \DGAs and \DGCs
discussed in the preliminary sections,
without recourse to the cochain level.
In particular, it proves possible
to almost entirely black-box the \SHCA technology,
and prior familiarity with notions other than 
\DGAs and \DGCs is not assumed.

\section{Algebras, coalgebras, and twisting cochains}\label{sec:CDGA}
Fix for all time a commutative base ring $\kk$ with unity
with respect to which all tensor products
and hom-modules are taken.
We take as understood the notions of 
differential graded $\kk$-modules and quasi-isomorphisms,
of differential graded $\kk$-algebras (henceforth \textcolor{RoyalBlue}{\DGA}s)
and differential graded $\kk$-coalgebras (\textcolor{RoyalBlue}{\DGC}s)
and maps between them, 
tensor products, 
and the Koszul sign convention.
A commutative \DGA is a \textcolor{RoyalBlue}{\CDGA}.

All algebras we consider are graded and
associative and all coalgebras graded and coassociative.
\emph{\textbf{All algebras we consider are augmented}}, 
with an important exception that we modify in short order to be augmented as well,
\emph{\textbf{and coalgebras coaugmented}}.
All \emph{\textbf{differentials}} $\defm d$ 
\emph{\textbf{increase degree}} by $1$,
and we use the terms \emph{\DG $\kk$-module} and \emph{cochain complex}
interchangeably.
Our \DGAs and \DGCs are \emph{\textbf{nonnegatively-graded}}. 
The multiplication $A \ox A \lt A$, unit $\kk \lt A$, and augmentation $A \lt \kk$
of a \DGA $A$ 
are respectively denoted $\defm \mu$, $\defm \h$, and $\defm \e$,
decorated with a subscript $A$ when necessary,
and the augmentation ideal $\ker \e \iso \coker \h$ is denoted $\defm{\bar A}$.
The comultiplication $C \lt C \ox C$, counit $C \lt \kk$, and coaugmentation $\kk \lt C$
of a \DGC $C$ are respectively denoted $\defm \D$, $\defm \e$, and $\defm \h$,
and the coaugmentation coideal $\coker \h \iso \ker \e$ is denoted $\defm{\bar C}$.

We write $\defm{\GM}$ for the category of graded $\kk$-modules 
(with no differential) and maps of fixed but arbitrary degree,
$\defm{\Algs}$ for the category of augmented $\kk$-\DGAs 
and augmentation-preserving \DGA maps
and $\defm{\Coalgs}$ for the category of coaugmented, \emph{cocomplete} 
$\kk$-\DGCs (cocompleteness will be explained shortly)
and coaugmentation-preserving \DGC maps.
All \DGA and \DGC maps will be of degree $0$,
but maps $C \lt A$ from a coalgebra to an algebra are allowed to be 
homogeneous of varying degrees,
as will also be explained momentarily.
The base ring $\kk$ itself
is considered a \DG Hopf algebra concentrated in degree zero, 
the differential, multiplication, and comultiplication being what they must.


We briefly rehearse some well-known generalities,
taking the opportunity to establish notation and 
conventions which will be leaned on throughout.
General background 
resources include Munkholm~{\cite[\SS1]{munkholm1974emss}},
Husemoller--Moore--Stasheff~\cite[Pt.~II]{husemollermoorestasheff1974},
the thesis of Prout\'e~\cite{proute2011}, and
the book of Loday--Valette~\cite[Chs.~1--2]{lodayvalette}.
One must mind the direction of the differential: 
while {our differentials \emph{increase} degree}, others' do not.

Given two graded $\kk$-modules $C$ and $A$,
we denote by $\defm{\GM_n}(C,A)$ 
the $\kk$-module of $\kk$-linear maps~$f$ sending each $C_{\mn j}$ 
to~$A_{\mn j+n}$,
and set the degree $\defm{|f|}$ to $n$ for such a map.
The hom-set $\GM(C,A) = \Direct_{n \in \Z} \GM_n(C,A)$ 
then becomes itself a graded $\kk$-module.
If $C$ and $A$ are cochain complexes,
then $\GM(C,A)$ becomes
a cochain complex under the
differential~$d = \defm{d_{\GM(C,A)}}$
given by $\defm{d(f)} \ceq d\mn_A f - (-1)^{|f|}f\mnn d_C$~{\cite[\SS1.1]{munkholm1974emss}},
cochain maps being described by
the condition 
$d(f) = 0$.
%

If $C$ is a \DGC and $A$ a \DGA,
then $\GM(C,A)$ becomes a \DGA
when endowed with the \defd{cup product}
$\defm{f \cup g} \ceq \muA(f \tensor g)\D_C$%
~{\cite[\SS1.8]{munkholm1974emss}}.\footnote{\ 
	To explain the nomenclature, 
	write out the definition of the product in the cohomology theory $E^*$
	represented by a ring spectrum $E$;
	the name \emph{convolution} is also popular.
}
An element  $t \in \GM_1(C,A)$
satisfying the three conditions
\[
\e_{\mn A} t = 0 = t \h_C,\qquad d(t) = t \cup t
\]
is called a \defd{twisting cochain}%
~\cite[\SS{1.8}]{husemollermoorestasheff1974}%
\cite[Prop.~3.5(1)]{husemollermoorestasheff1974}%
\cite[\SS\SS1.5,\,4]{proute2011}.
Twisting cochains compose with \DGC and \DGA maps in the sense that
given \DGAs $A'$, $A$ and \DGCs $C$, $C$'
and maps
\[
	C' \os g\lt C \os t\lt A \os f\lt A'\mathrlap,
\]
$g$ a \DGC map, $t$ a twisting cochain, and $f$ a \DGA map,
the maps $f\mnn t$, $tg$, and hence $f\mnn tg$ are all again twisting cochains.
Given a \DGA $A$,
there is a \emph{final} twisting cochain
$\defm{t^A}\: \defm {\B A} \lt A$
defined by the property that any 
twisting cochain $t\: C \lt A$
factors uniquely through a \DGC
map $\defm{g_t}\: C \lt \B A$
such that $t = t^A \o g_t$. 
We denote this conversion in the input--output ``deduction rule''
format borrowed from proof theory:
\[
\begin{adjunctions}
	g_t\: C & \B A\\
	t\: C & A \mathrlap.
\end{adjunctions}
\]

\nd 
The \DGC $\B A$ is referred to as the \defd{bar construction},
and gives the object component of a functor 
${\textcolor{RoyalBlue}{\B}}\:$$\Algs \lt \Coalgs$%
~{\cite[\SS1.6]{munkholm1974emss}\cite[\SS2.5]{proute2011}}.
The tautological twisting cochain $t^{(-)}\: \B \lt \id$
is a natural transformation.


More explicitly, 
the bar construction is the tensor coalgebra on 
the desuspension ${\desusp \bar A}$ of $\bar A$,
equipped with the sum of the tensor differential
and the unique coderivation extending
the ``bar-deletion'' map
$\smash{(\desusp \bar A)\ot \simto {\bar A}\ot \xtoo{\mu} \bar A \simto \desusp \bar A}$.
The tautological twisting cochain $t^A$ is the composition
of the projection $\B A \lt \desusp \bar A$ 
and the resuspension $\desusp \bar A \isoto \bar A$.

As we have defined $\Coalgs$ to contain only \emph{cocomplete} \DGCs, 
the foregoing assertions include the statement that 
$\B A$ is {cocomplete}, in the following sense.
Given a coaugmented \DGC $C$, 
the comultiplication $\D\: C \lt C \ox C$
defines by reduction a map $\defm{\bar \D}
= ({\id} - \h\e)^{\otimes 2}\D\: \bar C \lt \bar C \ox \bar C$
on the coaugmentation coideal, 
and by coassociativity, the iterates $\defm{\iter {\bar\D{}} n}\: \bar C \lt \bar C{}^{\otimes n}$
starting with $\iter {\bar \D{}} 1 = \id_{\bar C}$, \ 
$\iter {\bar \D{}} 2 = \bar \D$, \ 
$\iter {\bar \D{}} 3 = (\bar \D \ox \id)\bar\D = (\id \ox \bar\D)\bar\D$ are well-defined.
We say $C$ is \defd{cocomplete}\footnote{\ 
	\emph{Conilpotent} is probably currently more popular.
	} 
if $\bar C$ is exhausted by  
the increasing filtration by kernels $\ker \iter {\bar \D{}} n$.

Cocompleteness is the condition needed to extend
a twisting cochain $C \lt A$
to a \DGC map $C \lt \B A$,
and holds of any connected \DGC
and of the bar construction.
Cocompleteness is also important for another reason,
whose relevance will become clearer when we discuss homotopy%
~\cite[\SS1.3]{munkholm1974emss}.
If $C$ is cocomplete
and $h \in \GM_0(C,A)$ satisfies
$h \eta_C = \eta_A$,
then $(\eta_A \e_C - h)\eta_C = 0$.
Since $\im \eta_C = \ker({\id} - \h_C\e_C)$
and $\iter{\bar \D_C{}}{\ell} = 
({\id} - \h_C\e_C)^{\otimes \ell}\iter{\D_C{}}{\ell}$,
%
the cup-power 
\[
(\eta_A \e_C - h)^{\cup \ell} 
= \iter \muA \ell(\eta_A \e_{C} - h)^{\tensor \ell} \iter {\D_C} \ell
\]
annihilates the kernel of 
$\iter{\bar \D_C{}}{\ell}$.
These kernels exhaust $C$ by cocompleteness,
so 
$\sum_{\ell = 0}^\infty (\eta_A \e_C - h)^{\cup \ell}$ 
is a finite sum on any element of $C$
and hence gives a sensible two-sided cup-inverse $\defm{h^{\cup\,-1}}$ to $h$.

Given a \DGC $C$, there is also a twisting cochain
$\defm{t_C}\: C \lt \defm {\W C}$
\emph{initial} in the sense that any twisting cochain 
$t\: C \lt A$
factors uniquely through a $\DGA$ map 
$\defm{f^t}\: \W C \lt A$ 
such that $t = f^t t_C$:
\[
\begin{adjunctions}
t\: C & A \\
f^t\: \W C & A\mathrlap.
\end{adjunctions}
\]
The \DGA $\W C$ is referred to as the \defd{cobar construction},
and gives the object component of a functor $\defm\W\:\Coalgs \lt \Algs$%
~{\cite[\SS1.7]{munkholm1974emss}}.
The tautological twisting cochain $t_{(-)}\: \id \lt \W$
is a natural transformation.
Thus the two functors $\W \adj \B$ form an adjoint pair~{\cite[\SS1.9--10]{munkholm1974emss}}:
\[
\begin{adjunctions}
	g_t\: \W C & A \\
	f^t\:  C & \B A\mathrlap.
\end{adjunctions}
\]
We will have frequent recourse to the unit and counit of the adjunction $\W \adj \B$,
\[
	\defm\h\: \id \lt \B\W
\qquad\qquad\mbox{and}\qquad\qquad
	\defm\e\: \W\B \lt \id
\]
respectively. 
These are both natural quasi-isomorphisms 
and homotopy equivalences on the level of \DG modules%
~\cite[Thm.~II.4.4--5]{husemollermoorestasheff1974}%
\cite[Cor.~2.15]{munkholm1974emss}%
\cite[Lem.~1.3.2.3]{LH}.

There are a few elementary, purely categorical properties of this adjunction we will use,
here included for easy reference.

\begin{lemma}\label{thm:twisting-adjunction}
	Given a \DGA $A$ and a \DGC $C$, 
	one has $\e \o t_{\B A} = t^A\: \BA \lt A$
	and $t^{\W C} \o \h = t_C\: C \lt \W C$.
\end{lemma}

\begin{lemma}\label{thm:unit-counit}
	Given a \DGA map $f\:\W C \lt A$,
	we have \[\e_{\mn A} \o \W\B f\o \W\h_C = f\mathrlap;\] 
	in particular, $\e_{\W C} \o \W\h_C = \id_{\W C}$.
	Dually, given a \DGC map $g\: C \lt \B A$,
	we have \[\B\e_{\mn A} \o \B\W g \o \h_C = g\mathrlap;\]
	in particular, $\B\e_{\mn A} \o \h_{\B A} = \id_{\B A}$.
\end{lemma}

\begin{corollary}\label{thm:induction}
	If $C$ is a \DGC and $A$ a \DGA, 
	then any \DGA map $f\:\W C \lt A$ naturally factors through 
	$\e\: \W\B A \lt A$.
	We write $\defm{f^\#} = \W\B f \o \W\h_C\: \W C \lt \W \B A$ for the first factor.
	The transformation $f \lmt f^\#$ is natural in that
	given another \DGA map $\phi\: A \lt B$, 
	one has $(\phi \o f)^\# = \W\B\phi \o f^{\#}$.
\end{corollary}

We will sometimes say $f^\#\: \W C \lt \W \B A$ is 
\defd{induced up} from $f\: \W C \lt A$.

\begin{corollary}\label{thm:W-bar-factor}
	If $C$ is a \DGC and $A$ a \DGA, 
	and $g\: C \lt \B A$ is a \DGC map,
	then $g$ factors through $\B \W C$
	and $\W g\: \W C \lt \W \B A$
	factors as 
	\[\W C 
		\xtoo{\W\h} 
	\W\B\W C 
		\xtoo{\W\B\W g} 
	\W\B\W\B A 
		\xtoo{\W\B\e} 
	\W\B A
	\mathrlap.\]
\end{corollary}

\section{The tensor product}\label{sec:tensor}
Much of the material in this section is not proven
in the source literature,
and in longer preliminary drafts of this document,
proofs of each result were spelled out.
For space reasons we have again suppressed these, 
but it is still convenient to at least
gather the statements in one place.

The functor $\B\:$ $\Algs \lt \Coalgs$ is lax monoidal with respect to the 
monoidal structure given on both categories by the appropriate tensor products,
and $\W\: \Coalgs \lt \Algs$ is lax comonoidal.

\begin{definition}[See Husemoller \emph{et al.}~{\cite[Def.~IV.5.3]{husemollermoorestasheff1974}}]%
\label{def:shuffle}
There exist natural transformations 
\[
		\defm\nabla \: \B A_1 \ox \B A_2 \lt \B(A_1 \ox A_2),	
\qquad\qquad
		\defm\gamma \: \W( C_1 \ox C_2) \lt \W C_1 \ox \W C_2
\]
\nd of functors ${\Algs\,} \x {\,\Algs} \lt {\Coalgs}$ 
and $\Coalgs \,\x\, \Coalgs \lt \Algs$,
respectively, the 
\defd{shuffle maps}, 
determined by the twisting cochains
%
\[
			t^{A_1 \ox A_2}\nabla  = t^{A_1} \ox \h_{A_2}\e_{\B A_2} +
			\h_{A_1}\e_{\B A_1} \ox t^{A_2},
\qquad\qquad
		\gamma t_{C_1 \ox C_2}  = t_{C_1} \ox \h_{\W C_2}\e_{C_2} +
				\h_{\W C_1}\e_{C_1} \ox t_{C_2}\mathrlap.
\]
These are homotopy equivalences of cochain complexes and hence quasi-isomorphisms.
\end{definition}

\brmk
Although we manage to sidestep cochain-level computations completely in this paper
with the exception of the easy check in \Cref{thm:hmt-P},
it may be psychologically helpful to have an idea of what 
some of these maps do.
The bar shuffle~$\nabla$ on $\B_p A \ox \B_p B$
and cobar shuffle~$\g$ on $\W_\ell (C \ox D)$ 
respectively take
\eqn{
[a_1|\cdots|a_p] \otimes [b_1|\cdots|b_q]
&\,	 \lmt\,
\sum_{\s} \s\big(
[a_1 \otimes 1|\cdots|a_p \otimes 1|1 \otimes b_1|\cdots|1 \otimes b_q] 
\big)\mathrlap,\\
\ang{c_1 \ox d_1 ; \cdots ; c_\ell \ox d_\ell}
&	\,\lmt\,
\big(\ang{c_1} \ox \e(d_1) + \e(c_1) \otimes \ang{d_1}\big)
\cdots
\big(\ang{c_\ell} \ox \e(d_\ell) + \e(c_\ell) \otimes \ang{d_\ell}\big)
\mathrlap,
}
where the shuffles $\s$ are the $(p+q)!/p!q!$ permutations of the bar-word 
interlacing the ``$a$'' and ``$b$'' letters while retaining
the relative order of the $a_i$ and that of the $b_{\! j}$,
multiplied by a Koszul sign\footnote{\ That is,
$-1$ to the power of $\sum \big(|a_i|-1\big)\big(|b_{\! j}|-1\big)$,
where the sum runs over $a_i$ and $b_{\! j}$ that the shuffle
	has moved past one another.
	}.
Hence a typical term of the former sum for $p = q = 2$
is $\pm [1 \ox b_1 | a_1 \ox 1 | a_2 \ox 1 | 1 \ox b_2]$.
Sample values of the latter are
\eqn{
\g\ang{c_1 \ox d_1} &= 
    0,
	\phantom{\,\ang{c_1;c_2} \ox(-1)^{(|c_2|+1)(|d_1|+1)}}
        \mathrlap{\qquad c_1 \in \bar C, \ d_1 \in \bar D,}
        \\
\g\ang{c_1 \ox 1; c_2 \ox 1} &= 
	\ang{c_1;c_2} \ox 1,\phantom{(-1)^{(|c_2|+1)(|d_1|+1)}}
        \mathrlap{\qquad c_i \in \bar C,}
        \\ 
\g\ang{1 \ox d_1; c_2 \ox 1} &=
    (-1)^{(|c_2|+1)(|d_1|+1)} \ang{c_2} \otimes \ang{d_1}, 
        \mathrlap{\qquad c_2 \in \bar C,\  d_1 \in \bar D.}
}
Of course any identity involving these maps 
needs to be checked for all $p,q,\ell$ simultaneously
and thus involves infinitely many equations,
so computations involving multiple such maps
rapidly scale beyond tractability.

Fortunately, however,
explicit formulae are actually beside the point for our purposes.
All we really need to know about the
natural transformations $\nabla$ and $\g$  
(and also  $\psi$ and $\T$, to be introduced momentarily)
is that they allow us to
move tensor products in and out of $\B$ and $\W$
without changing homotopy type.
That is, we may treat them 
as devices for packaging
and formally manipulating homotopical information.
We will not really need to look under the hood,
and given the complexity of the machinery,
most of the time it is actually safer not to. 
\ermk

%


There is another important natural transformation on $\Algs$ 
we will rely on.

\begin{theorem}[{%
\cite[Prop.~IV.5.5]{husemollermoorestasheff1974}%
\cite[{$k_{A_1,A_2}$, p.~21, via Prop.~2.14}]{munkholm1974emss}%
}]%
\label{def:psi}
	There exists a natural transformation 
		\[\defm \psi\: \W\B(A_1 \ox A_2) \lt \W\B A_1 \ox \W\B A_2\]
	of functors $\Algs \x \Algs \lt \Algs$.
	This transformation satisfies 
	\[
		(\e_{\mn A_1} \ox \e_{\mn A_2}) \o \psi 
			= 
		\e_{\mn A_1 \ox A_2}\: \W\B(A_1 \ox A_2) \lt A_1 \ox A_2		
	\]
	and reduces to the identity if $A_1$ or $A_2$ is $\kk$.
\end{theorem}

\brmk
Husemoller--Moore--Stasheff's construction of such a map 
relies on a splitting arising
from the notion of an \emph{injective class} in their categorical 
reformulation of differential homological algebra, and is not very explicit.
Munkholm's map arises from the fact 
$\e_{\!A_1 \ox A_2}\: \W\B(A_1 \ox A_2) \lt A_1 \ox A_2$
is the initial object in a category of \emph{trivialized extensions}
of $A_1 \ox A_2$, of which $\e_{\!A_1} \ox \e_{\!A_2}$
is another object;
each is a \DGA quasi-isomorphism with a \DG-module quasi-inverse and 
a contracting homotopy.
The initiality follows  
the homotopy transfer theorem for \DGAs
and the finality of the tautological twisting cochain 
$t^{A_1 \ox A_2}\: \B(A_1 \ox A_2) \lt A_1 \ox A_2$,
so that Munkholm's $\psi$ is not given terribly explicitly either.
However, the homotopy transfer theorem
is proven via the inductive
construction of a certain twisting cochain 
from a given first stage, 
and Munkholm gives an explicit contracting homotopy for $\e$,
so that the construction is at least explicit enough that in a sequel paper,
we are able to show by computation that $\psi \o \W\nabla = \g$,
allowing us to see the product on Tor defined in this paper
agrees with that defined by Munkholm.
\ermk



\begin{definition}
  Let $A$ and $B$ be \DGAs.
A \DGC map $\B A \lt \B B$ 
is called an $\defm{A_\infty}$\defd{-map} from $A$ to $B$.
\end{definition}

Evidently  for $f\: A \lt B$ a \DGA map,
$\B f$ is an \Ai-map, 
but most \Ai-maps are not of this type.
The natural transformation $\psi$ allows us to also 
take tensor products of $A_\infty$-maps.

\begin{definition}[{\cite[Prop.~3.3]{munkholm1974emss}}]\label{def:internal-tensor}
	Let $A_1,A_2,B_1,B_2$ be \DGAs and 
	$g_{\mn j}\: \B A_{\mn j} \lt \B B_{\mn j}$ be \DGC maps for $j \in \{1,2\}$.
	Then we define the \defd{internal tensor product} 
	$\defm{g_1 \T g_2}\: \B(A_1 \ox A_2) \lt \B(B_1 \ox B_2)$ 
	as the
	composition
	\[
	\B(A_1 \ox A_2) 
		\xtoo{\!\h\!} 
	\B\W\B(A_1 \ox A_2)
		\xtoo{\B\psi}
	\B(\W\B A_1 \ox \W\B A_2)
		\xtoo{\B(\e \,\W g_1 \ox \,\e \, \W g_2)}
	\B(B_1 \ox B_2)
		\mathrlap.
	\]
\end{definition}

This construction exhibits some functoriality:

\begin{lemma}[{\cite[Prop.~3.3(ii)]{munkholm1974emss}}]%
\label{thm:T-compose}
Given \DGA maps $f_{\mn j}: A_{\mn j} \lt B_{\mn j}$ for $j \in \{1,2\}$,
	we have 
	\[
	\B\mnn f_1 \T \B\mnn f_2 = \B(f_1 \ox f_2)\: \B(A_1 \ox A_2) \lt \B(B_1 \ox B_2)
	\mathrlap.
	\]
	If $A'_j$ and $B'_j$ are further \DGAs and
	$g_j\: \B A'_j \lt \B A_j$ and
	$\ell_j\: \B B_j \lt \B B'_j$ \DGC maps,
	then 
	\[
		(\ell_1 \T \ell_2) \o \B(f_1 \ox f_2)
			=
		(\ell_1 \o \B\mnn f_1) \T {}(\ell_2 \o \B\mnn f_2) 
		\ \ \ \,
		\mbox{and}
		\ \ \ \,
		\B(f_1 \ox f_2)\o (g_1 \T g_2)
			=
		(\B\mnn f_1 \o g_1) \T {}(\B\mnn f_2 \o g_2) 
		\mathrlap.
	\]
\end{lemma}

The internal tensor product is related as one
would hope with the classical:

\begin{lemma}[{\cite[Lem.~4.4]{franz2019homogeneous}}]%
	\label{thm:tensor-psi}

	Let $A_{\mn j}$ and $B_{\mn j}$ be \DGAs for $j \in \{1,2\}$
	and $g_{\mn j}\: \B A_{\mn j} \lt \B B_{\mn j}$ \DGC maps.
	Then $\nabla \o (g_1 \ox g_2) = (g_1 \T g_2) \o \nabla$.
\end{lemma}

%

\counterwithin{figure}{subsection}
\counterwithin{definition}{subsection}
\counterwithin{theorem}{subsection}
\counterwithin{lemma}{subsection}
\counterwithin{corollary}{subsection}
\counterwithin{proposition}{subsection}
\counterwithin{equation}{subsection}
\counterwithin{notation}{subsection}
\numberwithin{remark}{subsection}

\section{Formal manipulation of homotopies}\label{sec:homotopy}

In this section we define the relevant notions of homotopy
and discuss how to package homotopies into representing path
(and path-allied) objects.

\begin{definition}%
	[{\cite[\SS1.11]{munkholm1974emss}\cite[\SS4.1]{munkholm1978dga}}]%
	\label{def:homotopy}
\ \\
	Given two \DGA maps $f_0,f_1\: A' \lt A$,
	a \textit{\textcolor{RoyalBlue}{\DGA}} \defd{homotopy}
	$f_0 \hmt f_1$
	is a $\kk$-linear map $h\: A' \lt A$ of degree $-1$
	such that 
	\[
	\e_{\mn A} h = 0,\qquad\quad
	h\eta_{A'} = 0,\footnote{\ %
		In the definition from our main source~\cite{munkholm1974emss}, 
		the unit and counit conditions are omitted;
		in later work dealing more specifically
		with $\Algs$ as a category,
		Munkholm includes them~\cite[4.1]{munkholm1978dga}.
		These are actually critical for the adjunction 
		to preserve homotopy and hence
		later to our verification of the path object.
	}\qquad\quad
	d(h) = f_0 - f_1,\qquad\quad
	h\mu_{\mn A'} = \mu_{\mn A}(f_0 \ox h + h \ox f_1)
	\mathrlap.
	\]
	Given two \DGC maps $g_0,g_1\: C \lt C'$,
	a \textit{\textcolor{RoyalBlue}{\DGC}} \defd{homotopy}
	$g_0 \hmt g_1$
	is a $\kk$-linear map $j\: C \lt C'$ of degree $-1$
	such that 
	\[
	\e_{C'} j = 0,\qquad\quad
	j\eta_{C} = 0,\qquad\quad
	d(j) = g_1 - g_0,\qquad\quad
	\D_{C'} j = (g_0 \ox j + j \ox g_1)\D_{C}
	\mathrlap.
	\]
	Given two twisting cochains $t_0,t_1\: C \lt A$,
	a \defd{twisting cochain homotopy}
	$t_0 \hmt t_1$
	is a $\kk$-linear map $x\: C \lt A$ of degree $0$
	such that 
	\[
	\e_{\mn A} x = \e_C,\qquad\quad
	x\h_{A} = \h_C,\qquad\quad
	\phantom{d(h) = f_0 - f_1,}\qquad
	d(x) = t_0 \cup x - x \cup t_1
	\mathrlap.
	\]
\end{definition}

\bs

These three notions evidently 
each compose well with maps in the appropriate categories.

\begin{lemma}\label{thm:homotopy-composition}
The postcomposition of a \DGA map to a \DGA or twisting cochain homotopy,
the precomposition of a \DGC map to a \DGC or twisting cochain homotopy,
the precomposition of a \DGA map to a \DGA homotopy, or
the postcomposition of a \DGC map to a \DGC homotopy
all result in another homotopy of the same type.
\end{lemma}

Moreover, the three notions are interchangeable under the adjunctions.

\begin{lemma}[{\cite[\SS1.11; Thm.~5.4, pf.]{munkholm1974emss}}]\label{thm:homotopy-adjunction}
	Suppose given a \DGC $C$ and a \DGA $A$.
	Then there are bijections of homotopies of maps
	\quation{\label{eq:homotopy-adjunction}
		\begin{aligned}
			\begin{adjunctions}
				\W C 	& A\\
				C 		& A\\
				C 		& \B A
			\end{adjunctions}
		\end{aligned}
	}
	The adjoint functors $\W \adj \B$
	thus also preserve the relation of homotopy.\footnote{\ 
		That the relation of homotopy is preserved
		is not to say that,
		for instance, if $h\: f_0 \hmt f_1\: A' \lt A$
		is a \DGA homotopy, then $\B h$ is a \DGC homotopy
		from $\B f_0$ to $\B f_1$,
		but that there exists a \DGC homotopy.
Early drafts of this document
addressed this
in much more detail, as
some of it is never explained in the primary sources.
	}
\end{lemma} 

Twisting cochain homotopies, 
despite being maps between different types of objects,
are in a way more flexible than \DGA or \DGC homotopies, 
because they are composable.

\begin{lemma}[{\cite[\SS1.12]{munkholm1974emss}}]\label{thm:homotopy-compose}
	Let a \DGC $C$ and a \DGA $A$ be given.
	A homotopy $h_{0,1}\: t_0 \hmt t_1\: C \lt A$
	of twisting cochains admits a two-sided cup-inverse
	$h_{0,1}^{\cup \,-1}$ which is a homotopy $t_1 \hmt t_0$.
	Given another homotopy $h_{-1,0}\: t_{-1} \hmt t_0\: C \lt A$,
	the cup product $h_{-1,0} \cup h_{0,1}$ is a homotopy $t_{-1} \hmt t_1$. 	
\end{lemma}

\brmk\label{rmk:eq-rel}
A suggestive phrasing is that the twisting cochains in $\GM_0(C,A)$ 
are the objects of a groupoid 
whose morphisms are the homotopies.
Particularly, homotopy is an equivalence relation.
The same then holds for the equivalent hom-sets
$\Algs(\W C,A) \longbij \Coalgs(C,\B A)$,
which are thus privileged over generic hom-sets
$\Algs(A',A)$ or $\Coalgs(C,C')$, which lack this property.
Note that cocompleteness is critical for the existence of inverses.
\ermk

We can use \Cref{thm:homotopy-adjunction}
to exchange \DGA homotopies $f_{-1} \hmt f_0 \hmt f_1$
we wish to concatenate for twisting cochains, 
take the cup product of these as in \Cref{thm:homotopy-compose},
and then move the resulting composite homotopy back to $\Algs$
using \Cref{thm:homotopy-adjunction}
to get a homotopy $f_{-1} \hmt f_1$.
The next subsections attempt to describe this process internally to $\Algs$.

\subsection{Path objects}\label{sec:path}

It is well known that the data of 
a homotopy $j\: g_0 \hmt g_1\: C \lt C'$ of maps of chain complexes 
can be realized as single map $C \ox I \lt C'$,
where $\defm I$ is the complex 
$\kk\{u_{[0,1]}\} \to \kk\{u_{[0]},u_{[1]}\}$
of nondegenerate chains
in the standard simplicial structure on the interval $[0,1]$.
Moreover, there is a natural coproduct endowing $I$
with a \DGC structure
so that \DGC homotopies can be realized in the same way.
Dually~\cite[Thm.~5.4, pf.]{munkholm1974emss}, 
the algebra of normalized cochains on the simplicial interval,
with the cup product, defines a \DGA $\defm{\I} = \kk\{v_0,v_1,e\}$
such that the data of a \DGA homotopy $h\: f_0 \hmt f_1\: A' \lt A$
can be realized by a \DGA map

\quation{\label{eq:hP}
	\begin{aligned}
		\defm{h^P}\:
		A' &\lt \I \ox A\mathrlap,\\
		a &\lmt v_0 \ox f_0(a) - e \ox h(a) + v_1 \ox f_1(a)\mathrlap.
\end{aligned}
}
Explicitly, 
the grading on $\I$ is given by $|v_0| = |v_1| = 0$ and $|e| = 1$,
the unit $\h\: \kk \lt \I$ by $\h(1) = v_0 + v_1$,
the nonzero differentials by $d v_1 = e = -dv_0$, 
and the multiplication by 

\begin{center}
	\begin{tabular}{ R || R | R | R|}
	\cup& 	v_0 & v_1	& e	\\
								 \hline\hline
	v_0 &	v_0 & 0 	& e \\ 
								\hline
	v_1 &	0 	& v_1 	& 0 \\  
								\hline
	e	&	0 	& e 	& 0 	\rlap{\phantom{X}}.%
								\\
								\hline
	\end{tabular}
\end{center}

It is easily seen that $\I$ has trivial cohomology 
$\H(\I) = H^0(\I) \iso \kk$
generated by the class of $v_0 + v_1$,
so by the K\"unneth theorem, the projections 
\[
\defm{\pi_{\mn j}} \: \I \ox A \longepi \kk\{v_{\mn j}\} \ox A \simto A
\]
are quasi-isomorphisms such that $\pi_{\mn j} \o h^P = f_{\mn j}$.
Thus tensoring with $\I$ functorially yields what 
we will call a \defd{naive path object} for \DGAs.

Unfortunately, none of the natural augmentations on $\I \ox A$
are such that both ``endpoint'' maps $\I \ox A \lt A$
are augmentation-preserving~\cite[\SS4.3]{munkholm1978dga}, 
which we need in order
to apply the natural transformation $\psi$
of \Cref{def:psi}
and to use the $\W \adj \B$ adjunction.
So we repair our path object by separating out 
$\kk \ceq \kk\{v_0 + v_1\} \ox \im \h_A$ as the image
of our unit and augmentation 
and defining the augmentation ideal
to be $\I \ox \bar A$:
\quation{\label{eq:P}
	\defm P \mnn A \ \mnn =\ \mnn \kk \,\oplus\, (\I \otimes \bar A)
	\mathrlap.\footnote{\ 
		Munkholm's description of this substitution%
		~\cite[p.~229, last line]{munkholm1978dga}
		seems to suffer from a typo.
		What is there does not quite parse as written, 
		and the most natural reading 
		yields $\kk\{v_0+v_1,e\} \+ (\I \ox \bar A)$,
		which is not quasi-isomorphic to $A$
		because $e$ generates a new $\kk$ summand in $H^1$
		now that $v_0$ and $v_1$ have been removed.
	}
}
The inclusion into $\I \ox A$ is a quasi-isomorphism,
and
the condition $\e h = 0$ on homotopies $h$ and  
unitality condition $f_{\mn j}(1) = 1$ on \DGA homomorphisms
guarantee 
the map $h^P\: A' \lt \I \ox A$ of \eqref{eq:hP}
lands in $P \mnn A$.

\begin{definition}\label{def:P}
	Given a \DGA $A$, we denote by $P \mnn A$ the augmented \DGA of \eqref{eq:P},
	equipped with the projections $\defm{\pi_0},\defm{\pi_1}\: P \mnn A \lt A$
	restricted from those of $\I \ox A$,
	and refer to it as the \defd{standard path object} of $A$.
	Given a homotopy $h\: f_0 \hmt f_1\: A' \lt A$ of \DGA maps,
	we refer to the associated \DGA map 
	$h^P\: A' \lt P \mnn A$
	of \eqref{eq:hP}
	as a \defd{right homotopy}
	and the composites $\pi_{\mn j} h^P = f_{\mn j}$ as its \defd{endpoint maps}.
	We write $\defm\zeta\: A \lt P \mnn A$ 
	for unital map defined by $a \lmt (v_0 + v_1) \ox a$
	on $\bar A$.
\end{definition}

	In later proofs,
	we will encounter many 
	right homotopies witnessed by \emph{nonstandard} path objects,
	which is to say \DGA quasi-isomorphisms
	$A' \to P' A \rightrightarrows A$ 
	representing homotopies
	$h\: f_0 \hmt f_1\: A' \lt A$
	through some other \DGA $P'A$ equipped with 
	two surjective quasi-isomorphisms to~$A$.
	Much of the material in this section is aimed at
	allowing us to convert these back to standard 
	right homotopies when needed
	in \Cref{sec:CGA,sec:ring-map}.

We now make an elementary observation about $P \mnn A$
that looks like it should follow purely model-categorically
but seems to be a fact about the category $\Algs$
(see \Cref{rmk:right-homotopy}).

\begin{proposition}\label{thm:hmt-P}
	Given a \DGA $X$, the standard path object $PX$
	is right-homotopy--equivalent to $X$ with respect to the notion 
	of right homotopy determined by $PX$ itself.
	In particular, for any other \DGA $A$,
	the mappings $\zeta$ and $\pi_{\mn j}$
	induce bijections $[A,X] \longbij [A,PX]$
	of right-homotopy classes of \DGA maps.
\end{proposition}
\begin{proof}
It would be enough to find a right homotopy between
$\zeta \o \pi_{\mn j}$ and $\id_{PX}$, 
	but it is psychologically more convenient to maintain $A$.
	Since $\pi_{\mn j} \o \zeta = \id_X$, 
	the map $\zeta_*\: [A,X] \lt [A,PX]$ is an injection,
	so it is enough to see it is surjective too. 
	For this, we note an arbitrary map 
	$A \lt PX$
	is a right homotopy $h^P$  
	between maps $f = \pi_0 \o h^P\: A \lt X$ and $g = \pi_1 \o h^P$,
	and show $h^P$ itself is right homotopic to the map $\zeta \o f$	
	representing the constant homotopy $f \hmt f$.
%
	The intuition for why this should be is given by the square
	\[
	\xymatrix@C=1em@R=.625em{
		f \ar@{-}[rr]^h \ar@{=}[dd]_0&&g\\ 
		&0&\\
		f\ar@{=}[rr]_0&&f \ar@{-}[uu]_h\mathrlap,
		}
	\]
	where we think of the left edge as $\zeta \o f$ and the right edge as $h^P$, 
	the labels $0$ and $h$ on the edges 
	representing the homotopies,
	which is to say the degree-$(-1)$ maps $A \lt X$
	which are the ``$e$-components'' of the right homotopies $\zeta \o f$
	and $h^P$.
	Explicitly, the degree-$(-1)$ map $A \lt PX$
	given by $a \lmt v_1 \ox h(a)$
	can be checked to be a \DGA homotopy $\zeta \o f \hmt h^P$
	in the sense of \Cref{def:homotopy}.
\end{proof}
\begin{corollary}\label{thm:homotopies-homotopic}
	Given a \DGC $C$ and a \DGA $X$,
	any two right homotopies $\W C \lt PX$
	both representing homotopies $f \hmt g$ of \DGA maps $\W C \lt X$
	are themselves homotopic as \DGA maps. 
\end{corollary}
\begin{proof}
	Note that $\zeta \o f$ is homotopic to	both right homotopies 
	by \Cref{thm:hmt-P}
	and recall from \Cref{rmk:eq-rel}
	that homotopy is an equivalence relation on $\Algs(\W C, X)$.
\end{proof}

\brmk\label{rmk:models}
There is a standard cofibrantly 
generated model structure on $\Algs$,
with weak equivalences quasi-isomorphisms
and surjections fibrations,
due when $\kk$ is a field 
to later work of Munkholm~\cite{munkholm1978dga},
and more generally to Jardine~\cite{jardine1997dga}.
The counit quasi-isomorphism $\e \:\W\B A \lt A$
is only a cofibrant replacement if $\kk$ is a field,
essentially because only projective modules will lift 
against surjections and $\W\B A$ is a projective $\kk$-module
if and only if $A$ is.
We will nevertheless be able to use $\W\B$ 
much in the manner of a functorial cofibrant replacement.\footnote{\ 
	The author at one point hoped to make real use of this model structure,
	but inexplicitly summoning liftings does not seem 
	to give enough control over composition;
	one wants something as near functorial as possible
	to describe composition of \DGA homotopies.
	Another thought was to use \DGCs instead, 
	but when $\kk$ is not a field, the standard proof for the model structure 
	breaks down and for some time the author 
	was under the impression 
	that for $\kk$ an arbitrary ring,
	there is no model structure
	(in fact, there is~\cite[Cor.~6.3.5]{HKRS2017necessary}).
}
\ermk

\brmk\label{rmk:right-homotopy}
The standard model-categorical notion of 
right homotopy uses generic path objects,
not specifically the standard path object of \Cref{def:P}
witnessing the classical notion of \DGA homotopy, 
which as we have pointed out
in \Cref{rmk:eq-rel}
is not typically an equivalence relation on $\Algs(A,X)$
unless $A = \W C$ for some \DGC $C$.
The natural witness $P \mnn A \x_{A} \mnn P \mnn A$ 
for a composite of homotopies
is not again $P \mnn A$, 
but \emph{is} another path object,
and model-categorical right homotopy on $\Algs$
is the transitive closure of classical \DGA homotopy.
Because of the functoriality and accessibility of the path objects
and the cobar--bar adjunction,
and to maintain back-compatibility with the classics,
we prefer to keep our discussion in terms of the classical notion.

\Cref{thm:hmt-P} looks,
if we took $RX = PX$,
a bit like the standard lemma
that the right-homotopy class of a(ny choice of) 
fibrant replacement $RA \lt RX$
of a map $A \lt X$ 
is determined by the composition $A \to X \to RX$.
But since all objects are fibrant,
this statement is distinct and actually trivial,
and $\z\: X \lt PX$ is very rarely a cofibration anyway.
The homotopy category of $\Algs$ is determined by 
right homotopy after cofibrant
replacement of domains, $QA \lt RX = X$,
so right-homotopy classes in $\Algs$$(\W\B A, X)$ give
the correct homotopy notion for $\kk$ a field,
but in general the relation is not clear.
\ermk

\subsection{Double- and triple-path objects 
	and concatenation}\label{sec:concatenation}

We have seen we can compose homotopies 
with compatible endpoint maps 
in $\Algs(\W C,A)$, for $A$ a \DGA and $C$ a \DGC, 
but our procedure passes through twisting cochains $C \lt A$
and gives no explicit description of the result 
in terms internal to $\Algs$. 
In this subsection, whose content we could find no other reference for, 
we rectify (or more honestly, circumnavigate) this 
shortcoming.

\begin{definition}\label{def:D}
	Given a \DGA $A$, we write 
	\[
		\defm D \mnn A \ceq \xu A {P \mnn A} {P \mnn A}
	\]
	for the pullback of the diagram $P \mnn A \os{\pi_1}\to A \os{\pi_0}\from P \mnn A$
	and refer to it as the \defd{double-path object}. 
	By definition it comes equipped with two projections to $P \mnn A$
	and three maps 
	\[\defm{p_0} \ceq \Dpra,
		\qquad\qquad
	 {\Dprb},
	 \qquad\qquad
	 \defm{p_1} \ceq \Dprc
	 \] 
	 to $A$, 
	all quasi-isomorphisms.
	As $P \mnn A$ is a subalgebra of $\I \ox A$,
	we may apply distributivity of $\ox$ over $\+ = \mnn\x$ 
	to identify $D \mn A$ with a subalgebra of $(\I \x \I) \ox A$,
	and then,
	in terms of the naive, unaugmented double path-object 
	$\I\x_\kk\I < \I \x \I$ of $\kk$,
	we have a decomposition
	\[
	D \mn A 
	  \ \iso\  
	\kk\big\{\mn(v_0,0) + (v_1,v_0) + (0,v_1)\big\} \,\+\,\,\mnn 
	(\xu\kk\I\I) \ox \bar A
	\]
	identifying $D \mn A$ as an augmented \DGA.
\end{definition}

The \emph{raison d'\^etre} 
of the double-path object, of course, is to represent pairs
of composable homotopies $f_{-1} \hmt f_0 \hmt f_1\: A' \lt A$
of \DGA maps,
which it achieves tautologically 
since a pair $h_{-1,0}\: f_{-1} \hmt f_0$
and $h_{0,1}\: f_0 \hmt f_1$ of homotopies 
induces representatives $h_{-1,0}^P$, $h_{0,1}^P\: A' \lt P \mnn A$
such that $\pi_1 h_{-1,0}^P = f_0 = \pi_0 h_{0,1}^P\: A' \lt A$,
and thus the map $(h_{-1,0}^P,h_{0,1}^P)\: A' \lt P \mnn A \x P \mnn A$
factors through the fiber product.
Evidently $p_0(h_{-1,0}^P,h_{0,1}^P) = f_{-1}$
and $p_1(h_{-1,0}^P,h_{0,1}^P) = f_1$.

If the desired composition of homotopies were realized by a map $D \mn A \lt P \mnn A$,
then the concatenation of any pair of compatible homotopies $A' \lt A$,
would be represented by the composite of the associated map 
$A' \lt D \mn A$ and the unattested $D \mn A \lt P \mnn A$,
but we know this is only possible when~$A'$
is the cobar construction $\W C$ on some \DGC $C$.
The composition of homotopies is, nevertheless, a natural transformation 
$\defm{\underline\Y}\:\Algs\big(\W(-),D \mn A\big) \lt \Algs\big(\W(-),P \mnn A\big)$.\footnote{\ 
	The intended visual mnemonic is that $\Y$ takes two things
	and combines them into one.
}
In particular, plugging $\B D \mn A$ in as the variable,
$\underline\Y$ takes the counit $\e\:\W\B D \mn A \lt D \mn A$
to a map $\defm\Y \ceq \underline\Y(\e)\: \W\B D \mn A \lt  P \mnn A$,
and a Yoneda-style argument yields the following.

\begin{lemma}\label{def:Y}
	Let a \DGC $C$ and \DGA $A$ 
	and homotopies $f_{-1} \hmt f_0 \hmt f_1$ of \DGA maps $\W C \lt A$
	be given. 
	If the \DGA map $\defm{h^D}\: \W C \lt D \mn A$
	represents this pair of homotopies,
	then the composite homotopy $f_{-1} \hmt f_1$
	is represented by
\quation{\label{eq:Y-comp}
	\W C \xtoo{\W \h} \W\B\W C \xtoo{\W\B h^D} \W\B D \mn A \xtoo{\Y} P \mnn A\mathrlap,
}
	where the map $\Y$ implementing the concatenation is a quasi-isomorphism.
\end{lemma}

We note that $\defm\w\:\W C \lt \W \B D \mn A$ 
is the map $(h^D)^\#$ induced up from $h^D$
as in \Cref{thm:induction}, 
so composition of homotopies in $\Algs$ is implemented 
by a diagram of the shape

\quation{\label{eq:Y-diagram}
	\begin{aligned}
		\xymatrix@C=1.5em@R=2.875em{
		& \W\B D \mn A\ar[d]^{\e}\ar[r]^(.56){\Y}&P \mnn A\mathrlap.\\
		\W C\ar[ur]^{\defm\w}\ar[r]_{h^D}& D \mn A
	}
	\end{aligned}
}

\begin{proof}
	Since $h^D = \w^* \e$ in \eqref{eq:Y-diagram},
	by naturality of $\underline\Y$ we have  
	$
		\underline\Y(h^D) 
		= \underline\Y \big(\w^* \e\big) 
		= \w^*\underline\Y(\e)
		= \w^* \Y
	$%
	.

	To see $\Y$ is a quasi-isomorphism, fix $C = \B A$
	and let $h^D\: \W\B A \lt D \mn A$
	represent a pair of constant (= trivial) homotopies of 
	$f_{-1} = f_0 = f_1 \ceq \e\: \W\B A \lt A$,
	so that explicitly
	$h^D\: x \lmt \big((v_0,0) + (v_1,v_0) + (0,v_1)\big) \ox \e(x)$.
	Then $h^D$ is a quasi-isomorphism by 2-of-3, 
	as $p_0 h^D = \e$ and $p_0$
	itself both are quasi-isomorphisms.
	If we follow the composite \eqref{eq:Y-comp}
	with the quasi-isomorphism $\pi_0$, 
	we recover the quasi-isomorphism $\e$,
	so by 2-of-3 again, 
	the composition in \eqref{eq:Y-comp} is a quasi-isomorphism.
	But $\W\h$ is always a quasi-isomorphism, 
	and $\W\B h^D$ is a quasi-isomorphism since $h^D$ is,
	so by 2-of-3 yet again, 
	so too must be $\Y$.
\end{proof}

The same trick works equally for composable triples of homotopies.

\begin{definition}\label{def:T}
	Given a \DGA $A$, its \defd{triple-path object} is the pullback
	\[
		\defm{T}\mnn A \ \ceq\  P \mnn A\, \us{A}\x\, P \mnn A\, \us{A}\x \, P \mnn A\mathrlap,
	\]
	equipped with 
	the expected three projections $TA \lt P \mnn A$
	and four projections $TA \lt A$.
\end{definition}

The same proof as for $D \mn A$ yields the following.

\begin{lemma}\label{def:Sha}
Let a \DGC $C$ and \DGA $A$ 
and homotopies $f_0 \hmt f_1 \hmt f_2 \hmt f_3$ of \DGA maps $\W C \lt A$
be given. 
There is a natural map $\defm\Sha\: \W\B TA \lt P \mnn A$
such that if the \DGA map $\defm{h^T}\: \W C \lt TA$
represents this triple of homotopies,
then the composite homotopy $f_{0} \hmt f_3$
is represented by
\[
	\W C \xtoo{\W \h} \W\B\W C \xtoo{\W\B h^T} \W\B TA \xtoo{\Sha} P \mnn A
	\mathrlap.
\]
The map $\Sha$ representing composition is a quasi-isomorphism.
\end{lemma}

In one instance we will encounter, 
the composable pair we apply \Cref{def:Y}
to comprises two adjacent sides of a 
square of homotopies $\W\C \lt (\I)\ot \ox A\ot$.

\begin{lemma}\label{thm:II-D}
	Let $A$ be a \DGA.
	Then there is a natural \DGA map $\defm r\: P \mnn A \ox P \mnn A \lt D(A \ox A)$
	such that $p_0 r = \pi_0 \ox \pi_0$ and $p_1 r = \pi_1 \ox \pi_1$.
\end{lemma}
\begin{proof}
	First note there exists a natural quotient map 
		\[
		\I \ox \I \lt 
		\kk\{v_0 \ox v_0,\, 
		e \ox v_0,\, 
		v_1 \ox v_0,\, 
		v_1 \ox e,\, 
		v_1 \ox v_1\}
		\]	
	modding out from the naive square object
	the ideal spanned by 
	$v_0 \ox e$,\, 
	$v_0 \ox v_1$,\, 
	$e \ox e$, and $e \ox v_1$.
	It is not hard to see this quotient is 
	isomorphic to the naive double-path object $\I\x_\kk\I$
	under the
	assignment
	\[
	v_0 \ox v_0 \lmt (v_0,0),\quad
	e \ox v_0 \lmt (e,0),\quad
	v_1 \ox v_0 \lmt (v_1,v_0),\quad
	v_1 \ox e \lmt (0,e),\quad
	v_1 \ox v_1 \lmt (0,v_1)\mathrlap.
	\]
	Write $\defm{r^I}\mn\: \I \ox \I \lt \I \x_\kk \I$
	for the composition
	and note it sends the idempotent
	$(v_0+v_1)\ot$ to the unity of $\I \x_\kk \I$.

	Recalling from \Cref{def:P} that $P \mnn A$ is a subalgebra of the tensor product
	$\I \ox A$,
	we may permute tensor-factors 
	to identify $P \mnn A \ox P \mnn A$ with a subalgebra of $(\I)\ot \ox A\ot$.
	The unity of $P \mnn A \ox P \mnn A$ is sent to 
	$(v_0+v_1)\ot \ox 1\ot$
	under this identification, and the augmentation ideal
	$\ol{P \mnn A \ox P \mnn A}$ into $(I^*)\ot \ox \ol{A\ot}$.
	Thus $r^I \ox \id_A\ot$ sends the unity of $P \mnn A \ox P \mnn A$
	to the unity of $D\mn A$
	and $\ol{P \mnn A \ox P \mnn A}$ 
	into $(\I \x_\kk \I) \ox \ol{A \ox A} = \ol{D(A\ox A)}$,
	so we may define the intended $r$ by corestriction.
\end{proof}

The last formal trick we will perform with homotopies
is to move $P$ past $\W\B$.
For this, starting with a right homotopy 
$h^P\: \W C \lt P \mnn A$ representing a homotopy $h\: f_0 \hmt f_1\: \W C \lt A$,
note the string of conversions
\[
\begin{adjunctions}
	h^P\: \W C 			& P \mnn A\\
	h\: \W C	 		& A\\
	{\wt h}\: \W C 		& \W\B A\\
	\wt h^P\: \W C 	& P \W B A
\end{adjunctions}
\]
afforded by \Cref{thm:homotopy-adjunction},
\Cref{thm:induction}, 
and \Cref{def:P},
amounting to a natural transformation
\[
\defm{\underline Z}\: \Algs\big(\W(-),P \mnn A\big) \lt \Algs\big(\W(-),P\,\mn\W\B A\big)\mathrlap.
\]
%
Following through the construction,
the induced $\wt h$ is a homotopy $f_0^\# \hmt f_1^\#$,
so that $\pi_{\mn j} \o \wt h^P = f_{\mn j}^\#$ for $j \in \{0,1\}$.

%
%
%

As with $\underline\Y$, a Yoneda-esque argument
shows the natural transformation $\underline Z$ is represented by one map.
Plugging $\B P \mnn A$ into the hole and applying $\underline Z$ to 
the counit $\e\: \W\B P \mnn A \lt P \mnn A$
yields a \DGA map $\defm Z \ceq \underline Z(\e)\: \W\B  P \mnn A \lt P \W \B A$.

\begin{lemma}\label{thm:P-W}
	Given a \DGA $A$,
	there is a natural \DGA quasi-isomorphism $Z\: \W\B P \mnn A \lt P\,\mn\W\B A$
	such that $\pi_{\mn j} \o Z = \W\B\pi_{\mn j}\: \W\B  P \mnn A \lt \W\B A$.
\end{lemma}
\begin{proof}
	The counit $\e\: \W\B  P \mnn A \lt  P \mnn A$ itself 
	represents a homotopy $h$
	between the two composites 
	$\pi_{\mn j} \o \e\: \W\B P \mnn A \lt A$ for $j \in \{0,1\}$,
	and the induced $\wt h$ is a homotopy between the maps
	${(\pi_{\mn j} \o \e)^\#} \: \W\B P \mnn A \lt \W\B A$.
	By the naturality property \Cref{thm:induction} of induction,
	$(\pi_{\mn j} \o \e)^\# = \W\B \pi_{\mn j} \o \e^\#$
	and $\e^\# = \id_{\W \B A}$,
	so $\smash{\wt h}$ is a homotopy $\W\B\pi_0 \hmt \W\B\pi_1$,
	and thus by definition $Z = \smash{\wt h^P}$
	satisfies $\pi_{\mn j} \o \wt h^P = \W\B\pi_{\mn j}$ 
	for $j \in \{0,1\}$.
	It follows $Z$ is a quasi-isomorphism,
	for $\pi_{\mn j}$ and $\W\B\pi_{\mn j}$ both are.
\end{proof}

\brmk
Though we will not use these facts, 
we should mention that the transformations $\underline\Y$, $\underline\Sha$, $r$, and $\underline Z$
(from \Cref{def:Y,def:Sha,thm:II-D,thm:P-W} respectively)
are also natural in the variable $A$
and $r$ is a quasi-isomorphism.
\ermk

\counterwithin{figure}{section}
\counterwithin{definition}{section}
\counterwithin{theorem}{section}
\counterwithin{lemma}{section}
\counterwithin{corollary}{section}
\counterwithin{proposition}{section}
\counterwithin{equation}{section}
\counterwithin{notation}{section}
\numberwithin{remark}{section}

\section{Maps on Tor}\label{sec:Tor}
One of the goals 
of this work is to determine the cohomology ring of the homotopy
pullback of a span $X \to B \from E$ of spaces from the cohomology 
of the input spaces.
The tool of choice here is the \EMSS,
which is a special instance of 
a so-called \emph{algebraic \EMSS},
converging~\cite[XI.3.2]{maclane} 
to differential Tor of a triple of \DGAs
from classical Tor of their cohomology
and functorial in all three variables.


\begin{lemma}[{\cite[Cor.~1.8]{gugenheimmay}%
		\cite[Theorem~5.4]{munkholm1974emss}}]
	\label{thm:Tor-quism}
	Given a \DGA map $f\: R' \lt R$,
	a right $R'$-module~$M'$, a left $R'$-module $N'$,
	a right $R$-module $M$, a left $R$-module $N$,
	and \DG module maps $u\:M' \lt M$ and $v\:N' \lt N$ making the expected squares
	\quation{\label{eq:Tor-functoriality-squares}
		\begin{aligned}
			\xymatrix@C=1.5em{
				M' \ar[d]
				&	\ar[l] 
				M' \ox R\mathrlap'
				\ar[d]
				&
				\qquad
				&	
				R' \ox N'\,
				\ar[r]
				\ar[d]
				& 	N\mathrlap'
				\ar[d]\\
				M
				&	\ar[l]
				M \ox R
				&
				\qquad
				&	R \ox N
				\ar[r]	
				& N
			}
		\end{aligned}
	}
	commute, there is induced a map of algebraic {\EMSS}s
	from that of $(M',R',N')$ to that of $(M,R,N)$,
	converging
	to the functorial map
	\[
	\Tor_f(u,v)\:	\Tor_{R'}(M',N') \lt \Tor_{R}(M,N)
	\]
	of graded modules.
	Moreover, if the maps $f$, $u$, $v$
	are quasi-isomorphisms,
	then the map of spectral sequences is an isomorphism from the $E_2$ page on 
	and $\Tor_f(u,v)$ is an isomorphism.
\end{lemma}

We will only apply these considerations in the 
special case that $M$ and $N$ are \DGAs and the $R$-module structure maps
are induced by \DGA homomorphisms $M \from R \to N$,
so that we have the condensed compatibility diagram
\quation{\label{eq:Tor-DGA-functoriality-squares}
	\begin{aligned}
		\xymatrix@C=2em{
			M' \ar[d]_(.45)u
			&	\ar[l]_{\phi_{M'}} \ar[r]^{\phi_{N'}}
			R'
			\ar[d]|(.375)\hole|(.45)f|(.525)\hole
			& 	N\mathrlap'
			\ar[d]^(.45)v\\
			M
			& R
			\ar[l]^{\phi_M} \ar[r]_{\phi_N}	
			& N
		}
	\end{aligned}
}
of \DGA maps. 
In later sections,
where we will produce diagrams of Tors comprising
mainly isomorphisms, we will ceaselessly apply this result.
We will also need to expand the notion of a map of Tors 
to include squares which commute only up to homotopy.


\begin{lemma}[{\cite[Thm.~5.4]{munkholm1974emss}}]%
	\label{thm:homotopy-Tor-map}
	Let \DGAs and \DGA maps as in \eqref{eq:Tor-DGA-functoriality-squares}
	be given such that the squares commute up to homotopies
	$h_M\: u \o \phi_{M'} \hmt \phi_M \o f$ and
	$h_N\: v \o \phi_{N'} \hmt \phi_N \o f$.
	Then there is induced a map 
	\[
	\defm{\Tor_{f}(u,v;h_M,h_N)}\:
	\Tor_{R'}(M',N')
	\lt
	\Tor_R(M,N)
	\]
	of graded modules
	which is a quasi-isomorphism if each of $u$, $f$, and $v$ is.
\end{lemma}
\begin{proof}
	Letting $h_M^P\: R \lt P \mnn M'$
	and $h_N^P\: R \lt PN'$
	be the respective right homotopies
	representing the homotopies $h_M$, $h_N$,
	as described in \Cref{def:P},
	the following equivalent diagrams commute by definition:
	\quation{\label{eq:Tor-DGA-homotopy-squares}
		\begin{aligned}
\xymatrix@C=1.90em{
  M' \ar[r]^u							&M					&P \mnn M \ar[l]_{\pi_0}\ar[r]^(.525){\pi_1}	&M\ar@{=}[r]		&M
  \\R' \ar[u]^{\phi_{M'}} \ar[d]_{\phi_{N'}}\ar@{=}[r]		&R'\ar[u] \ar[d]\ar@{=}[r]		&R'\ar[u]|(.425)\hole|{h^P_M}|(.575)\hole \ar[d]|(.425)\hole|{h^P_N}|(.575)\hole \ar@{=}[r]			&R'\ar[u] \ar[d] 	\ar[r]|(.475)\hole|f|(.525)\hole	&R\ar[u]_{\phi_M} \ar[d]^{\phi_N}
\\
  N' \ar[r]_v			&N				&PN \ar[l]^{\pi_0}\ar[r]_(.525){\pi_1}	&N\ar@{=}[r]		&N\mathrlap,
			}
			\qquad
			\quad
			\qquad
\xymatrix@C=1.90em{
  M' \ar[r]^u			&M				&P \mnn M \ar[l]_{\pi_0}\ar[r]^(.525){\pi_1}  		&M
  \\R' \ar[u]^{\phi_{M'}} \ar[d]_{\phi_{N'}}\ar@{=}[r]		
  &R'\ar[u] \ar[d] 		\ar@{=}[r]			
  &R'\ar[u]|(.425)\hole|{h^P_M}|(.575)\hole \ar[d]|(.425)\hole|{h^P_N}|(.575)\hole \ar[r]|(.475)\hole|f|(.525)\hole
  &R\ar[u]_{\phi_M} \ar[d]^{\phi_N}
\\
  N' \ar[r]_v			&N				&PN \ar[l]^{\pi_0}\ar[r]_(.525){\pi_1}			&N\mathrlap.
			}
		\end{aligned}
	}
	Since the $\pi_{\mn j}$ are quasi-isomorphisms,
	three applications of \Cref{thm:Tor-quism}
	let us set
	\[
	\Tor_{f}(u,v;h_M,h_N) 
	\ceq
	\Tor_{f}(\pi_1,\pi_1) \o \Tor_{\id}(\pi_0,\pi_0)\- \o \Tor_{\id}(u,v)\mathrlap.
	\qedhere
	\]
\end{proof}

This result accounts for our fixation 
on representing homotopies by \DGA maps in \Cref{sec:homotopy}.
\wrong{It seems to be a new observation that the choice of homotopies doesn't matter.}

\begin{proposition}\label{thm:homotopy-independence}
\wrong{
	The map $\defm{\Tor_{f}(u,v)} \ceq \Tor_{f}(u,v;h_M,h_N)$
	is independent of the homotopies $h_M$, $h_N$.
}
\end{proposition}
\begin{proof}
\wrong{
	First assume $R' = \W C$ is the cobar construction of some \DGC $C$.
	By \Cref{thm:homotopies-homotopic}, then
	the right homotopy $h^P_M\: R' \lt P \mnn M$
	and any other right homotopy $\wt h^P_M\: R' \lt P \mnn M$
	representing some other homotopy $\wt h_M\: u\o \phi_{M'} \hmt \phi_M \o f$
	are homotopic,
	and hence there is some right homotopy
	$H^P_M\: R' \lt  P \mnn P M$ with $\pi_0 \o H^P_M = h^P_M$ 
	and $\pi_1 \o H^P_M = \wt h^P_M$.
	Similarly, given a second homotopy  $\wt h_N\: v\o \phi_{N'} \hmt \phi_N \o f$,
	we have a homotopy $h^P_N \hmt \wt h^P_N$
    and witnessing right homotopy $H^P_N\: R' \lt P \mnn PN\vphantom{x_{A_{A_A}}}$.
	By construction, one has $\pi_0 \o \pi_0 \o H^P_M = \pi_0 \o \pi_1 \o H^P_M$
	and $\pi_1 \o \pi_0 \o H^P_M = \pi_1 \o \pi_1 \o H^P_M$,
	and similarly for $H^P_N$.
	This gives us the commutative diagram of maps below, 
	where the Tor in the central object
	is defined by $H^P_M$ and $H^P_N$.
}
\wrong{
	\begin{equation}\label{eq:diamonds}
	\begin{aligned}
	\xymatrix@R=2.5em@C=1.5em{
		&	
		&	
		&
		\Tor_{R'}(P \mnn M,PN)
		\ar[ld]^{\sim}_{\Tor_{\id}(\pi_0,\pi_0)\ }
		\ar[rd]^{\ \Tor_{\id}(\pi_1,\pi_1)}
		&
		\\
		\Tor_{R'}(M',N')	
		\ar[rr]^{\Tor_{\id}(u,v)}
		&	
		&
		\Tor_{R'}(M,N)	&
		\Tor_{R'}(P \mnn P \mnn M,P \mnn PN)
		\ar[u]_(.35){\Tor_{\id}(\pi_0,\pi_0)}^\vertsim
		\ar[d]^(.25){\Tor_{\id}(\pi_1,\pi_1)}_\vertsim
		&
		\Tor_{R'}(M,N)
		\ar[rr]^{\Tor_{f}(\id,\id)}
		&
		&
		\Tor_R(M,N)
		\\
		&
		&	
		&
		\Tor_{R'}(P \mnn M,PN)
		\ar[lu]_{\sim}^{\Tor_{\id}(\pi_0,\pi_0)\ \ }
		\ar[ru]_{\ \ {\vphantom{X^X}}\Tor_{\id}(\pi_1,\pi_1)}
	}
	\end{aligned}
	\end{equation}
	The horizontal composite along the top,
	${\raisebox{-5pt}{$\to$} \!\mnn\nearrow \searrow \!\mnn \raisebox{-5pt}{$\to$}}$,
	is $\Tor_f(u,v;h_M,h_N)$
	while the bottom composite
	${\raisebox{5pt}{$\to$} \!\mnn \searrow \nearrow \!\mnn \raisebox{5pt}{$\to$} }$
	is $\Tor_f(u,v;\wt h_M,\wt h_N)$.
	Both can be factored 
	through $\Tor_{R'}(P \mnn P \mnn M,P \mnn PN)$
	using the vertical isomorphisms,
	the top as 
	${\raisebox{-5pt}{$\to$}  \!\mnn \nearrow \downarrow \mn \uparrow \searrow \! \raisebox{-5pt}{$\to$}}$,
	the bottom as ${\raisebox{5pt}{$\to$}}  \!\mnn \searrow \uparrow \mn \downarrow \nearrow \! \raisebox{5pt}{$\to$}$,
	but by commutativity of the diagram,
	one has ${\nearrow \downarrow} \,=\, {\searrow \uparrow}$
	and ${\uparrow \searrow} \,=\, {\downarrow \nearrow}$,
	so the composites are equal.
}
\wrong{
	In the general case, we have the diagram \eqref{eq:diamonds}
	but with the central $\Tor_{R'}(P \mnn P \mnn M,P \mnn PN)$ omitted,
	and we do not know \emph{a priori} the top and bottom composites
	are the same.
	If we precompose all maps from $R'$ in the discussion above
	with the counit $\e\: \W\B R' \lt R'$,
	then the preceding discussion with $C = \B R'$ does apply,
	and we obtain the commutative diagram \eqref{eq:diamonds}
	but with $R'$ everywhere replaced by $\W\B R'$.
	In that diagram, we do know the top and bottom composites are equal
	by the preceding paragraph.
	Now, deleting the central $\Tor_{\W\B R'}(P \mnn P \mnn M,P \mnn PN)$
	from that diagram, the rest admits the map
	$\Tor_{\e}(\id,\id)$ down to the first diagram, 
	so that	by the naturality of $\e$,
	all new squares commute.
	As $\e$ is a quasi-isomorphism,
	we know
	$\Tor_{\e}(\id,\id)$ is a natural isomorphism
	by \Cref{thm:Tor-quism},
	and thus
	equality of the two composites in the $R'$
	diagram follows from 
	the equality of the corresponding composites in the $\W\B R'$	diagram.
}
\end{proof}

\begin{manualrmk}{4.6b}
\new{The issue with the preceding proof is that, though we have correctly noted 
that $\pi_0 \o \pi_0 \o H^P_M = \pi_0 \o \pi_1 \o H^P_M$
	and $\pi_1 \o \pi_0 \o H^P_M = \pi_1 \o \pi_1 \o H^P_M$,
	we still have $\pi_0 \o \pi_0 \neq \pi_0 \o \pi_1$ and
			$\pi_1 \o \pi_0 \neq \pi_1 \o \pi_1$,
	and it is this that we would need for the central diamond of Tors to commute.
	Morally, what we need is for $H_M$ and $H_N$ to be an endpoint-fixing rather than free 
	homotopies.}
\end{manualrmk}

\begin{manualprop}{4.6c}\label{thm:homotopy-independence-mk2}
\new{The map 
	$ \Tor_{f}(u,v;h_M,h_N)$
depends only on the endpoint-fixing homotopy classes 
of the homotopies $h_M$ and $h_N$,
in the following sense:
given other homotopies 
	$\ell_M\: u \o \phi_{M'} \hmt \phi_M \o f$ and
	$\ell_N\: v \o \phi_{N'} \hmt \phi_N \o f$,
if there exist \DGA homotopies
$H_M\: h_M^P \hmt \ell_M^P$ and
$H_N\: h_N^P \hmt \ell_N^P$
such that 
\[
P\pi_0 \o H_M^P = \z \o u \o \phi_{M'},
\quad
P\pi_0 \o H_N^P = \z \o v \o \phi_{N'},
\quad
P\pi_1 \o H_M^P = \z \o \phi_M \o f,
\quad
P\pi_1 \o H_N^P = \z \o \phi_N \o f,
\]
then 
	$ \Tor_{f}(u,v;h_M,h_N)
	= \Tor_{f}(u,v;\ell_M,\ell_N)
	$.
}
\end{manualprop}
\begin{proof}
  \new{Both the maps of Tors under discussion are composites of 
$ 
\Tor_{\id}(u,v)\:
		\Tor_{R'}(M',N')	
		\lt
		\Tor_{R'}(M,N)	
$
preceding and
$
		{\Tor_{f}(\id,\id)}\:
		\Tor_{R'}(M,N)	
		\lt
		\Tor_{R}(M,N)	
$
following some map 
		$\Tor^{(u\phi_{M'},v\phi_{N'})}_{R'}(M,N)
			\lt
		\Tor^{(\phi_{M}f,\phi_{N}f)}_{R'}(M,N)$,
where we have included the maps $M \from R' \to N$
defining these distinct but isomorphic Tors as superscripts.
We now need to compare these maps as defined by $h_M^P$ and $h_N^P$
on the one hand and by  $\ell_M^P$ and $\ell_N^P$
on the other.
The postulated right homotopies $h_M^P$, $\ell^M_P$, and $H_M^P$
fit into a commuting cone of maps from $R'$ to the following diagram,
where we tag each codomain with the map it receives from $R'$:
\[
\xymatrix@C=4em{
&  \mbox{$\overbrace{M}^{{u\phi_{M'}}}$} 
& \mbox{$\overbrace{PM}^{{h_M^P}}$}  
\ar[l]_{\pi_0}\ar[r]^{\pi_1}
&  \mbox{$\overbrace{M}^{{\phi_{M}f}}$}  
& \\
\smash{\mathllap{\mbox{\footnotesize $u\phi_{M'}$}} \Big\{ } M 
\ar[r]_\z \ar[ur]^\id \ar[dr]_\id 
& 
\smash{\mathllap{{ } ^{\z u\phi_{M'}\!\!\!\!\!\rotatebox{-45}{$\Big\{$}}\!\!\!}}
PM  
\ar[u]_{\pi_0}\ar[d]^{\pi_1}
& 
\smash{\mathllap{{ } ^{H^P\!\!\!\!\!\rotatebox{-45}{$\Big\{$}}\!\!\!}}
PPM \ar[u]_{\pi_0}\ar[d]^{\pi_1}
\ar[l]^{P\pi_0}\ar[r]_{P\pi_1}
& 
\smash{\mathllap{{ } ^{\z \phi_{M}f\!\!\!\!\!\rotatebox{-45}{$\Big\{$}}\!\!\!}}
PM \ar[u]_{\pi_0}\ar[d]^{\pi_1}
& M \smash{\Big\}  }
\mathrlap{\mbox{\footnotesize $\phi_{M}f$};} 
\ar[l]^\z \ar[ul]_\id \ar[dl]^\id 
&\\
&  \mbox{$\underbrace{M}_{{u\phi_{M'}}}$}  
& \mbox{$\underbrace{PM}_{{\ell_M^P}}$}  
\ar[l]^{\pi_0}\ar[r]_{\pi_1}
&  \mbox{$\underbrace{M}_{{\phi_{M}f}}$}  & 
}
\]
and similarly for $N$.
The commutativity of the triangles falls out of the definition
and the commutativity of the cones over the squares come from the fact
that $H_M^P$ and $H_N^P$ are homotopies and the displays.
There is an induced commutative diagram of Tors and isomorphisms,
in which we have compactified the subscripts and superscripts to fit the page
width:
\[
\xymatrix{
&  
\Tr {R'}\os{\mathclap{(u\phi_{M'},v\phi_{N'})}}{(M,N)}
&  
\Tr{R'}\os{\mathclap{(h_M^P,h_N^P)}}{(PM,PN)}
\ar[l]\ar[r]
&  
\Tr{R'}\os{\mathclap{(\phi_{M}f,\phi_{N}f)}}{(M,N)}
& \\
\Tr{R'}\os{\mathclap{(u\phi_{M'},v\phi_{N'})}}{(M,N)}
\ar[r] \ar[ur]^\id \ar[dr]_\id 
& 
\Tr{R'}\os{\mathclap{(\z u\phi_{M'},\z v\phi_{N'})}}{(PM,PN)}
\ar[u]
\ar[d]
& 
\Tr{R'}\os{\mathclap{(H_M^P,H_N^P)}}{(PPM,PPN)}
\ar[u]\ar[d]
\ar[l]\ar[r]
& 
\Tr{R'}\os{\mathclap{(\z\phi_{M}f,\z\phi_{N}f)}}{(PM,PN)}
\ar[u]\ar[d]
& 
\Tr{R'}\os{\mathclap{(\phi_{M}f,\phi_{N}f)}}{(M,N)}
\ar[l] \ar[ul]_\id \ar[dl]^\id 
\mathrlap.
&\\
&  
\Tr{R'}\os{\mathclap{(u\phi_{M'},v\phi_{N'})}}{(M,N)}
& 
\Tr{R'}\os{\mathclap{(\ell_M^P,\ell_N^P)}}{(PM,PN)}
\ar[l]\ar[r]
&  
\Tr{R'}\os{\mathclap{(\phi_{M}f,\phi_{N}f)}}{(M,N)}
& 
}
\]
The path along the top is the middle factor in
	$ \Tor_{f}(u,v;h_M,h_N)$
	and that along the bottom is
	$ \Tor_{f}(u,v;\ell_M,\ell_N)$,
	so the two are equal.}
\end{proof}

\wrong{Not only do these maps not depend on the homotopies involved,}
\new{So far as these maps are well-defined,
we will show} 
they compose functorially.
We already see the diagrams  
in the preceding proof stretching
the limits of what can be fit on normal page,
and things are only going to get worse from here.
To make what follows more legible, 
we introduce a convention to save on the mental and physical space 
required for repeating triples of operations involved 
in defining maps of Tors.
We will make increasing use of this convention as the diagrams evolve,
as we will eventually arrive at a point where there is no other choice.

\begin{notation}\label{def:suppression}
	Given \DGA maps $X \from A \to Y$,
	functors $F,G,F',G'\: \Algs \lt \Algs$,
	and natural transformations
	$F \lt G$,\, $F' \lt G'$,\, 
	$\phi\:F \lt F'$,
	and $\chi\:G \lt G'$
	such that the two composites $F \lt G'$ are equal,
	we make the abbreviations
	\[
	\defm{\Tor_{FA}} \ceq \Tor_{FA}(FX,FY),
	\qquad\mnn\qquad\mnn\qquad\mnn
	\defm{\Tor_{FA}(GX)} \ceq \Tor_{FA}(GX,GY)\mathrlap,
	\]
	\vspace{-1.5em}
	\[
	\defm{\Tor_\phi} \ceq \Tor_\phi(\phi,\phi)\:
	\Tor_{FA} \lt \Tor_{F'A},
	\qquad\quad
	\defm{\Tor_\phi(\chi)} \ceq \Tor_\phi(\chi,\chi)\:
	\Tor_{FA}(GX) \lt \Tor_{F'A}(G'X)\mathrlap.
	\]
\end{notation}


The following, apparently original,\footnote{\ 
but \emph{cf}. work of Gugenheim--Munkholm~\cite{gugenheimmunkholm1974}
which achieves something similar for a version of Tor
defined as the cohomology of the two-sided bar construction
}
shows the functoriality 
of Tor with respect to this extended mapping notion.
\begin{theorem}\label{thm:functoriality}
Assume given a diagram of \DGA maps
\[
\xymatrix{
\ar[d]_{u'}M''& R'' \ar[d]|(.48)\hole|{f'}	\ar[l]\ar[r]&	N''	\ar[d]^{v'}	\\
\ar[d]_{u}M'& R' \ar[d]|(.48)\hole|{f'}\ar[l]\ar[r]	& N' \ar[d]^{v}		\\
M& R\ar[l]\ar[r] &	N 	\\
}
\]
in which there are \DGA homotopies 
\new{%
$h'_M$, 
$h'_N$, 
$h_M$,
$h_N$} 
making each square commute.
Then \new{a} $\Tor_{f\!\! f'}(uu',vv'\new{;\ell_M,\ell_N})$ 
\wrong{is well} \new{can be} defined 
\new{in such a way as to} \wrong{and} 
equal\wrong{s}
$\Tor_{f}(u,v;\new{h_M,h_N}) \o \Tor_{f'}(u',v';\new{h'_M,h'_N})$.
\end{theorem}

\begin{proof}
If $R''$ is a cobar construction,
we define $\Tor_{\mn f\!\! f'}(uu',vv')$
using the composite \DGA homotopies~
\new{
$\ell_M$ associated to the right homotopies $h_M^P \o f'$ and $Pu \o h'_M$
and
$\ell_N$ associated to the right homotopies $h_N^P \o f'$ and $Pv \o h'_N$
and
}
guaranteed by \Cref{thm:homotopy-compose}
 making the vertical rectangles
 $R''M''M R$ and $R''N''NR$ commute.
\wrong{We have seen in \Cref{thm:homotopy-independence} that this 
choice of homotopies does not matter so long as any exist in the first place.}

If $R''$ is not a cobar construction, such composite homotopies do not necessarily exist.
However, since Tor is functorial in its arguments
and $\e\: \W\B R'' \lt R''$ is a quasi-isomorphism,
the induced map $\Tor_{\e}(\id,\id) \:\Tor_{\W\B R''}({-},{-}) \lt \Tor_{R''}({-},{-})$
is a natural isomorphism, and conjugating by these isomorphisms,
we may define $\Tor_{\mn f\!\! f'}(uu',vv')$ as 
$\Tor_{\e}(\id,\id) \o \Tor_{\W\B \mn f\, \W\B \mn f'}(uu',vv') \o
\Tor_{\e}(\id,\id)\-$.
By naturality of the isomorphism $\Tor_{\e}(\id,\id)$,
the composite 
 $\Tor_{f}(u,v) \o \Tor_{f'}(u',v')$
 can be written equally well as
\[
\big[\mn\Tor_{\e}(\id,\id) \o \Tor_{\W\B f}(u,v) \o \Tor_{\e}(\id,\id)\-\big]\o 
\big[\mn\Tor_{\e}(\id,\id) \o \Tor_{\W\B f'}(u',v') \o \Tor_{\e}(\id,\id)\-\big]\mathrlap,
\]
so it is enough to compare 
$\Tor_{\W\B \mn f \, \W\B \mn f'}(uu',vv')$ with
$\Tor_{\W\B \mn f}(u,v) \o  \Tor_{\W\B \mn f'}(u,v)$,
and thus we may assume $R'' = \W C$ for some \DGC $C$ to begin with.

Now, omitting the $N$ arguments, taking for granted that the diagrams are symmetric
in $M$'s and $N$'s, 
the four-square diagram \emph{\`a la} \eqref{eq:Tor-DGA-homotopy-squares} giving $\Tor_{\mn f\!\! f'}(uu',vv')$ can be rewritten as
\quation{\label{eq:one-shot}
	\begin{aligned}
\xymatrix{
R'' \ar[d]\ar@{=}[r]&
\ar@{=}[r] R''  \ar@{=}[r]\ar[d]&
R'' \ar@{=}[r]\ar[d]&
R'' \ar@{=}[r]\ar[d]^{\new{\ell_M^P}}&
R''\ar[r]^{f'}\ar[d]&R'\ar[r]^f\ar[d]&R\ar[d]\\
M''\ar[r]_{u'}&M'\ar[r]_{u}&M&P \mnn M\ar[r]_{\pi_1}\ar[l]^{\pi_0}&M \ar@{=}[r]&M \ar@{=}[r]&M
}
	\end{aligned}
}
while the diagram yielding $\Tor_{\mn f}(u,v) \o \Tor_{\mn f'}(u',v')$ is
\quation{\label{eq:two-shot}
	\begin{aligned}
\xymatrix{
	R'' \ar@{=}[r]\ar[d]& R'' \ar@{=}[r]\ar[d]&
	R'' \ar@{=}[r]\ar[d]^{\new{h'_M}}& R'' \ar[r]^{f'}\ar[d]& 
	R' \ar@{=}[r]\ar[d]& R' \ar@{=}[r]\ar[d]& R' \ar@{=}[r]\ar[d]^{\new{h_M}}& 
	R' \ar[r]^{f}\ar[d]& R\ar[d]\\
	M'' \ar[r]_{u'}& M' &\ar[l]^{\pi_0}P \mnn M'\ar[r]_{\pi_1}&\ar@{=}[r]M'&M'\ar[r]_{u}&M&\ar[l]^{\pi_0}P \mnn M\ar[r]_{\pi_1}&M \ar@{=}[r]&M\mathrlap.
}
	\end{aligned}
}
By \Cref{thm:homotopy-compose},
we can concatenate the right homotopies $R'' \xtoo{\new{h'_M}} P \mnn M' \xtoo{Pu} P \mnn M$
and $R'' \os{f'}\to R' \xtoo{\new{h_M}} P \mnn M$ implied by \eqref{eq:two-shot},
and by \Cref{def:Y} this composite homotopy can be implemented via a map
with codomain $\W\B D \mnn M$
followed by $\Y\: \W\B D \mnn M \lt P \mnn M$.
We can combine these ingredients into \Cref{fig:functoriality-diagram},
abusively labelling maps on Tor by the \DGA maps inducing them.
\begin{figure}
\centering
$
\xymatrix@C=1em{
\smash{\Tr{R''}}(M'')\ar@[jviolet][d]_(.525){u'}		\\
\smash{\Tr{R''}}(M')\ar@[jred][d]_(.5){u}&
\smash{\Tr{R''}}(P \mnn M')\ar[d]_{Pu}\ar@[RoyalBlue][r]^\sim\ar@[RoyalBlue][l]_\sim&
\smash{\Tr{R''}}(M')\ar[d]_{u}\ar@[RoyalBlue][r]^{f'}
			&
			\smash{\Tr{R'}}(M')\ar@[RoyalBlue][d]_u			\\
\smash{\Tr{R''}}(M) &
\smash{\Tr{R''}}(P \mnn M) \ar[r]^\sim\ar[l]_\sim&
\smash{\Tr{R''}}(M) \ar[r]^{f'}&
\smash{\Tr{R'}}(M) &
\smash{\Tr{R'}}(P \mnn M) \ar@[RoyalBlue][r]^\sim\ar@[RoyalBlue][l]_\sim
			&
\smash{\Tr{R'}}(M) \ar@[jviolet][r]^f&
\Tr{R}(M) \mathrlap.	\\
			&
\smash{\Tr{R''}}(\W\B DM) \ar[u]^\sim
							\ar[rrr]_\sim
							\ar[dl]^\sim_{\Y(\e)}
							\ar[lu]^\sim_{p_0}
							\ar@/_1.5pc/[rrrr]_(.6)\sim_(.4){p_1}
			&				&				&
\smash{\Tr{R''}}(P \mnn M) \ar[u]_{f'}
					\ar@/^1pc/[ull]^(.65)\sim
					\ar[r]^\sim
			&
{\phantom{\Tor}}{\mathllap{\Tr{R''}}}(M) \ar@[jred][u]_{f'}				\\
\smash{\Tr{R''}}(P \mnn M) \ar@[jred][uu]^{\vertsim}\ar@/_3pc/@[jred][urrrrr]_\sim		
}
$
\caption{The diagram for functoriality of Tor}
\label{fig:functoriality-diagram}
\end{figure}
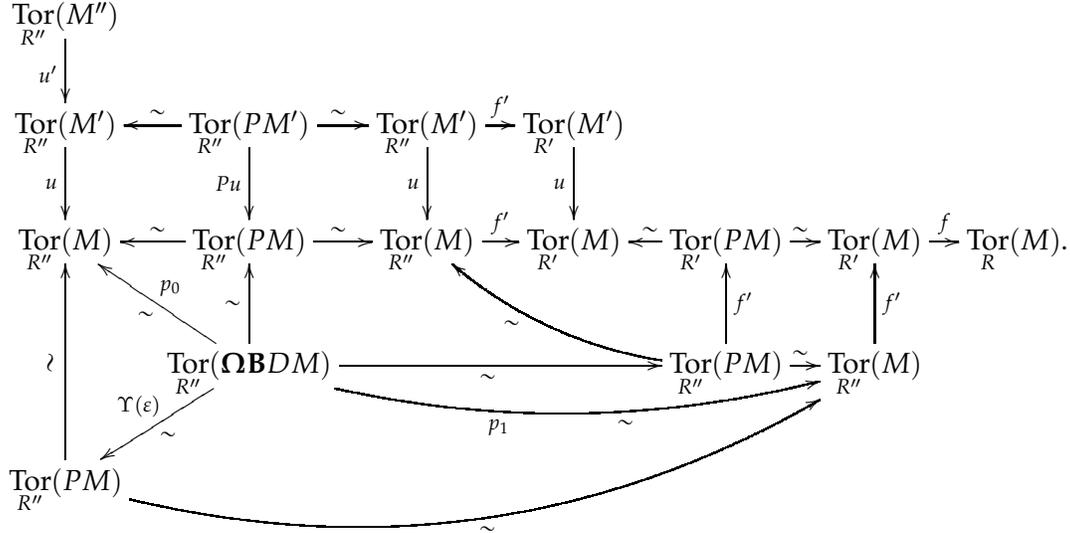
Here the maps induced by \eqref{eq:one-shot},
along the left and bottom, are in red 
and those induced by \eqref{eq:two-shot}, along the top, are in blue;
the first and last arrows, labeled $u'$ and $f$ respectively,
are violet because they are both.
Using the commutativity of the various squares
and triangles,
one sees $\Tor_{\mn f\!\! f'}(uu',vv')$
and $\Tor_{\mn f}(u,v) \o \Tor_{\mn f'}(u',v')$ are indeed equal.
%
%
\end{proof}

\section{SHC-algebras}\label{sec:SHC}

A commutative \DGA $A$ is one for which the multiplication 
$\mu\: A \ox A \lt A$ is itself a \DGA homomorphism.
Cohomology rings are of this sort, 
and a large part 
of why homotopy theory is so much more tractable
over a field $\kk$ of characteristic $0$ is that
there are functorial \CDGA models for cochains.
For other characteristics
this is not the case~\cite[Thm.~7.1]{borel1951leray}, 
but we can weaken the requirement
by asking only that $\mu$ extend to an \Ai-algebra map.
Munkholm's product is defined in terms of such a structure,
as first considered by Stasheff and Halperin.

To make sense of the following definition, 
recall from \Cref{sec:CDGA} 
that the canonical twisting cochain 
$t^A\: \B A \lt A$ of a \DGA $A$ 
factors through the projection onto a \DG direct summand $\desusp \bar A$,
on which it restricts to a cochain isomorphism
of degree $1$,
whose inverse we dub $\defm{\desusp_A}\: \bar A \lt \B A$ 
and call the \defd{desuspension}.

\begin{definition}[Stasheff--Halperin~{\cite[Def.~8]{halperinstasheff1970}}]\label{def:WHC}
	We refer to a \DGA $A$
	equipped with a \DGC map $\defm\Phi_A\: \B(A \ox A) \lt \B A$ 
	such that the composition $t_A \o \Phi \o \desusp_{A \ox A}\: 
	\ol{A \tensor A} \lt \bar A$
	is 
	the multiplication $\muA\: A \ox A \lt A$
	as a \defd{weakly homotopy commutative} (\textcolor{RoyalBlue}{\WHC}\defd{-})\defd{algebra}.\footnote{\
		Stasheff--Halperin call
		this a \emph{strongly homotopy commutative algebra} 
		structure, but we will meet a stronger notion momentarily,
		so we rename their notion ``weak.''
	}
	Given two {\WHCA}s $A$ and $Z$,
	a \textcolor{RoyalBlue}{\WHC}\defd{-algebra map}
	from $A$ to $Z$ is 
	a \DGC map $g\: \B A \lt \B Z$
	such that there exists a \DGC homotopy between 
	the two paths around
	the square
	\quation{\label{eq:SHC-map}
		\begin{aligned}
	\xymatrix@C=1.75em@R=3em{
		\mathllap\B(A\ox A)  
		\ar[d]_{g \,\,\mn\T\,\mnn g}		
		\ar[r]^(.565){\Phi_A}
		&
		\B A
		\ar[d]^{g}								\\
		\mathllap\B(Z \ox Z)
		\ar[r]_(.565){\Phi_Z}			
		\ar@{}[r]_{{\vphantom{g}}}				&
		\B Z\mathrlap.
	}
		\end{aligned}
}
The particular homotopy is not prescribed as part of the data of a \WHCA map.
\end{definition}

W$\mspace{-1.5mu}$\textsc{hc}-algebra
structures thus enable us to upgrade the non-\DGA map $\mu\: A \ox A \lt A$
to a legitimate \DGA map $\W\Phi\: \W\B(A\ox A) \lt \W\B A$,
which is more tractable categorically even if less intuitive on the point-set level,
and which by the naturality of $\e\: \W\B \lt \id$ 
carries the same information as $\mu$ up to quasi-isomorphism.

Munkholm's definition
adds to Stasheff--Halperin's a weakening of the standard \CGA axioms.

\begin{definition}%
{\cite[Def.~4.1]{munkholm1974emss}}
\label{def:SHC}
	A \defd{strongly homotopy commutative} 
	(henceforth \textcolor{RoyalBlue}{\SHC}\defd{-})
	\defd{algebra}
	is a \WHCA $A$ whose structure map
	$\Phi_A\mn\: \B(A \tensor A) \lt \B A$,
	satisfies the following conditions:
	\begin{enumerate}
		\item\label{def:SHC-unit} 
		It is strictly unital:
		$\Phi \o \B(\id\mn_A \ox\ \h_A) = \id_{\B A} = \Phi \o \B(\h_A \ox \id\mn_A)$.
		
		\item\label{def:SHC-comm} 
		It is homotopy-commutative:
		there is a \DGC homotopy 
		from $\Phi$
		to $\Phi \o \B\chi\: \B(A \tensor A) \lt \B A$, 
		where $\chi\: A \tensor A \isoto A \tensor A$
		is the factor transposition $a \ox b \lmt (-1)^{|a||b|} b \ox a$.	
		
		\item\label{def:SHC-assoc} 
		It is homotopy-associative:
		there is a \DGC homotopy
		between the maps
		$\Phi(\Phi \T \id\mn_{\B A})$
		and
		$\Phi(\id\mn_{\B A} \T \Phi)\:
		\B(A \tensor A \tensor A) \lt \B A$.
	\end{enumerate}
%
	The associativity and commutativity homotopies postulated 
	are again not themselves specified in the data of an \SHCA,
	only the fact of their existence. 
	An \textcolor{RoyalBlue}{\SHC}\defd{-}\defd{algebra map}
	is a \WHCA map between {\SHCA}s.
\end{definition}


The canonical example is that of an authentically commutative algebra.

\bex\label{thm:CGA-SHC}
If $A$ is a \CDGA, 
then the morphism $\Phi = \B\muA\: \B(A \ox A) \lt \B A$ makes $A$ 
an \SHCA.
The cohomology ring $\H(X;\kk)$ of a simplicial set is of this type,
and will always come considered with this \SHCA structure.
If $\rho\: A \lt B$ is a map of \CDGAs,
then $\B\rho$ is an \SHCA map.
\eex

\bex\label{thm:H-SHC}
If $A$ is any \DGA,
there is a unique \DG-module section $\defm{i}\: A \lt \W\B A$ of 
$\e\: \W\B A \lt A$
that is unital and restricts to $t_{\B A}\o \desusp$ on $\bar A$~%
\cite[Prop.~2.14]{munkholm1974emss}.
If $A$ is an \SHCA,
then the homotopy-commutativity of $\Phi$ 
implies $\mu$ and $\mu \o \chi$ 
are homotopic cochain maps,
and $\e \o \W(\Phi \o \B\chi) \o i = \mu \o \chi$
and $\e \o \W\Phi \o i = \mu$,
and so $\H(A)$ is a \CDGA (with trivial differential).
Thus $\Phi = \B \mu_{\H(A)}$ gives an \SHCA structure on $\H(A)$
by \Cref{thm:CGA-SHC}.
\emph{\textbf{The cohomology ring of an}} 
{\SHC}\emph{\textbf{-algebra will always be endowed with this}} 
{\SHC}\emph{\textbf{-algebra structure}}.
\eex
%
%

\begin{theorem}[{\cite[Prop.~4.7]{munkholm1974emss}}]\label{thm:SHC-cochain}
	Let $X$ be a simplicial set and $\kk$ any ring.
	Then the normalized cochain algebra $\defm{\C}(X) = \C(X;\kk)$,
	augmented by restriction to $\C({*};\kk) \iso \kk$ for
	some basepoint ${*} \in X_0$,
	admits an \SHCA structure $\Phi_{\C(X)}$,
	and this structure is \emph{strictly natural} in the sense
	that given a basepoint-preserving map $f\: Y \lt X$ of simplicial sets,
	the induced \DGC map $\B \C(f)\: \B \C(X) \lt \B \C(Y)$
	renders the square \eqref{eq:SHC-map}
	commutative on the nose.
\end{theorem}

\nd This natural \SHCA structure on cochains 
is a reinterpretation of the classical Eilenberg--Zilber theorem,
and only verifying the homotopy-associativity axiom requires
substantial work. 

The class of known \SHCAs has recently expanded significantly:

\begin{theorem}[Franz~{\cite{franz2019shc}}]\label{thm:Franz-SHC}
A homotopy Gerstenhaber algebra $A$ admits a 
strictly natural \WHCA structure $\Phi_A$
satifsying the axioms \Cref{def:SHC}.\ref{def:SHC-unit} and \Cref{def:SHC}.\ref{def:SHC-assoc}.
If $A$ is an \emph{extended} homotopy Gerstenhaber algebra,
then $\Phi_A$ is in fact an \SHCA structure.
\end{theorem}

\brmk\label{rmk:operads}
A homotopy Gerstenhaber algebra is an algebra over a certain $E_2$-operad
$F_2 \ms X$
and similarly, an extended homotopy Gerstenhaber algebra is an algebra over 
a certain suboperad of an $E_3$-operad $F_3 \ms X$,
accounting for the phrasing we employ in the abstract.\footnote{\ 
	Both of these are 
	filtrands of the so-called \emph{surjection operad} $\ms X$ 
	of interval-cut
	operations on cochains,
	which is a quotient of the 
	\DG-operad $\ms E$ associated to the classical Barratt--Eccles simplicial operad~\cite{mccluresmith2003,bergerfresse2004operad}.
}
Since $E_2$ is not very far along the road to $E_\infty$,
morally speaking we require \emph{some} commutativity to obtain 
the product on Tor, but not very much.

On the other hand, since the readiest source of 
homotopy Gerstenhaber algebras is algebras over this particular operad, 
rather than just any $E_2$-operad, 
an $E_2$-algebra is not necessarily a \WHCA,
so the notions are not strictly comparable.
\ermk

\section{The product}\label{sec:product}
Munkholm's product
can be motivated as a sort of least common generalization of
the classical products on $\Tor_{\C B}(\C X, \C E)$
and $\Tor_{\H B}(\H X, \H E)$,
rephrased
in terms of the canonical \SHCA structures.
We choose not to use his definition of the product,
but an equivalent definition of our own.\footnote{\ 
	In a sequel,
	we will elaborate on both definitions and prove their equivalence,
	but it is not necessary to show our product agrees with Munkholm's
	in order to use it.
}

Given \DGAs $R_0$, $R_1$ and right and left \DG $R_i$-modules 
$M_i$ and $N_i$ respectively,  
there is a classically defined exterior product~%
\cite[p.~206]{cartaneilenberg}
\[
\Tor_{R_0}(M_0, N_0) \ox\Tor_{R_1}(M_1, N_1)
\lt 
\Tor_{R_0 \ox R_1}(M_0 \ox M_1, N_0 \ox N_1)\mathrlap,
\]
functorial in all six variables in the sense that
given similarly defined $R'_i$, $M'_i$, $N'_i$
such that the squares \eqref{eq:Tor-functoriality-squares} commute,
so does the square
\[
\xymatrix{
\Tor_{R_0'}(M_0',N_0') \ox \Tor_{R_1'}(M_1',N_1') \ar[r] \ar[d]
&\Tor_{R_0' \ox R_1'}(M_0' \ox M_1',N_0' \ox N_1') \ar[d] 
\\
\Tor_{R_0}(M_0,N_0) \ox \Tor_{R_1}(M_1,N_1) \ar[r]
&\Tor_{R_0 \ox R_1}(M_0 \ox M_1,N_0 \ox N_1)\mathrlap,
}
\]
and given further $R''_i,M''_i,N''_i$, such squares glue.
If $R = R_0 = R_1$ is a \emph{commutative} \DGA,
then $\mu\: R' = R \ox R \lt R$ is a \DGA map,
and if $M = M_0 = M_1$ and $N = N_0 = N_1$ are themselves \DGAs,
then $\mu\: M' = M \ox M \lt M$
and $\mu\: N' = N \ox N \lt N$
make the diagrams \eqref{eq:Tor-functoriality-squares} commute,
so we may follow the external product with the map
\[
\Tor_\mu(\mu,\mu)\:	\Tor_{R \ox R}(M \ox M, N \ox N)
\lt
\Tor_R(M, N)
\]
to obtain the classical product on Tor.
This particularly applies to $R = \H\mn A$, \ $M = \H X$, \ $N = \H Y$
for $X \from A \to Y$ maps of spaces.

If $A$ fails to be commutative, this fails to give a product,
but taking $\C(X) \from \C(B) \to \C(E)$ as $X \from A \to Y$,
one can 
use 
the natural \DGA maps
\[
\C(B) \ox \C(B) 
\lt 
(C_* B \ox C_* B)^*
\xleftarrow{\nabla^*}
{\C(B \x B)}
\xtoo{\C(\D)} 
\C(B)
\]
inducing the cup product (the first two are quasi-isomorphisms,
so the direction of $\nabla^*$ is not an issue)
to obtain a map
\[
	\Tor_{\C B \ox \C B} (\C X \ox \C X,\C E \ox \C E) \lt
	\Tor_{C^*B}(C^*X,C^*E)
\]
which we apply after the exterior product, and this yields the product on Tor.
In the situation of the Eilenberg--Moore theorem,
this product can be shown to be preserved by the isomorphism with $\H(X \x_B E)$%
~\cite[Corollary 7.18]{mcclearyspectral}%
\cite[Cor.~3.5]{gugenheimmay}%
\cite[Prop.~3.4]{smith1967emss}%
\cite[Thm.~A.27]{carlsonfranzlong}.%
\footnote{\ 
	No source the author knows actually shows the product is preserved,
	but McCleary at least reduces it to an exercise,
	and Carlson--Franz~\cite[A.27]{carlsonfranzlong}
	spell out some of the steps to this exercise.
}

Munkholm is able to describe 
both these products as instances of another product.\footnote{\ 
It takes a little work to see these products as instances 
of Munkholm's, or that the product we give here agrees with Munkholm's,
and we will spell out the details
in a sequel article comparing definitions of products,
but here we will take the specialization as given.}
We assume given the following homotopy-commutative squares of \DGC maps.
\begin{equation}\label{eq:SHC-map-2}
\begin{aligned}
\xymatrix@C=1.75em@R=4em{
	\B(X \ox X)\ar[d]_{\Phi_X}& 
	\B(A \ox A)\ar[d]|(.45)\hole|{\Phi_A}|(.55)\hole
	\ar[l]_{\xi \,\T\, \xi}
	\ar[r]^{\upsilon \,\T\, \upsilon}	
	& \B(Y \ox Y)\ar[d]^{\Phi_Y} \\
	\B X	& \B A			\ar[l]^{\defm\xi}\ar[r]_{\defm\upsilon}	
	& \B Y
}
\end{aligned}
\end{equation}
Applying $\W$,
which preserves the relation of homotopy by \Cref{thm:homotopy-adjunction},
this induces a map
$\Tor_{\W\B(A \ox A)} \lt \Tor_{\W\B A}$ by \Cref{thm:homotopy-Tor-map},
where we have utilized the abbreviation convention in \Cref{def:suppression}.
The exterior product is a map $(\Tor_{\W\B A})^{\otimes 2} \lt \Tor_{(\W\B A)\ot}$,
in this notation,
so to define a candidate product
we must connect $\Tor_{(\W\B A)\ot}$
with $\Tor_{\W\B (A\ot)}$.
Munkholm does this using $\psi$ and $\e$,
but we can do it in a somewhat simpler way 
applying the natural quasi-isomorphisms
\[
{\W\B{Z}} \ox {\W\B{Z}} \xleftarrow{\gamma} \W(\B {Z}\ox\B {Z}) \xtoo{\W\nabla} \W\B({Z} \ox {Z})
\]
of \Cref{def:shuffle}
to the span $X \from A \to Y$.
The map $\Tor_\g$ goes in the wrong direction,
but this is no issue by \Cref{thm:Tor-quism},
since $\g$ is a \quism.
All told, one gets the following composite.

\begin{definition}\label{def:product}
Given \WHCA maps and homotopies as
in \eqref{eq:SHC-map-2},
the \defd{product on Tor} is 
\[
\xymatrix@C=6em{
		\us{\W\B A}\Tor(\W\B X,\W\B Y)
		\ox
		\us{\W\B A}\Tor(\W\B X,\W\B Y)
	\ar[r]^(.475){\mathrm{ext}}
&
		\us{\smash{\W\B A \ox \W\B A}}{\Tor}(\W\B X \ox \W\B X,\W\B Y \ox \W\B Y)
		\ar@{<-}[d]^(.525){\Tor_\g(\g,\g)}_(.525){\vertsim}
\\
		\us{\W\B( A \ox A)}{\Tor}\big(
		\W\B( X \ox X),\W\B( Y \ox Y)
		\big)
	\ar[d]^{\Tor_{\id}(\W\Phi,\W\Phi)}
&
	\ar[l]_{\Tor_{\W\nabla}(\W\nabla,\W\nabla)}^\sim
		\us{\W(\B A \ox \B A)}{\Tor}\big(\W(\B X \ox \B X),
						 \W(\B Y \ox \B Y)\big)
\\
		\us{\W\B(A \ox A)}{\Tor}(\W\B X, \W\B Y)
&
		\us{\W\B(A \ox A)}{\Tor}(P\,\mn\W\B X, P\,\mn\W\B Y)
	\ar[l]_{\Tor_{\id}(\pi_0,\pi_0)}^\sim
	\ar[d]_{\Tor_{\W\Phi}(\pi_1,\pi_1)}
\\
&
		\us{\W\B A}{\Tor}(\W\B X, \W\B Y)
	\mathrlap.
}
\]
Compactifying notation as per \Cref{def:suppression},
this is
	\quation{\label{eq:full-product-v2}
		\Big(\,\us{\W\B A}\Tor\,\Big){}\ot
		\os{\mr{ext}}\to			
		\us{(\W\B A)\ot}{\Tor}
		\us{\sim}{\os{\Tr\g}\from}
		\us{\W(\B A)\ot}{\Tor}
		\xtoo[\sim]{\!\Tr{\W\nabla}\!}
		\us{\W\B(A\ot)}{\Tor}
		\!\!\!\!
		\os{\us{\id}\Tor(\W\Phi)\!\!\!\!\!\!}\lt
		\!\!\!
		\us{\W\B(A\ot)}\Tor\!(\W\B X)
		\!\!
		\us{\sim}{\os{\us {\id}\Tor(\pi_0)\!\!\!\!}\lf}
		\!\!\!\!
		\us{\W\B(A\ot)}\Tor\!(P\,\mn\W\B X) 
		\!\!
		{\,\,\os{\us{\W\Phi}\Tor(\pi_1)\!\!\!\!}\lt}
		\,
		\us{\W\B A}{\Tor}\mathrlap.
	}
\end{definition}

\counterwithin{figure}{subsection}
\counterwithin{theorem}{subsection}
\counterwithin{lemma}{subsection}
\counterwithin{proposition}{subsection}
\counterwithin{corollary}{subsection}
\counterwithin{equation}{subsection}
\counterwithin{notation}{subsection}
\numberwithin{remark}{subsection}

\section{The algebra structure on Tor}\label{sec:CGA}

As noted in the introduction,
Munkholm's product depends for its definition on 
a choice of homotopies making 
\eqref{eq:SHC-map-2} commute, 
and he conjectured that its properties might therefore be bad. 
\wrong{We aim to show this is untrue,
proving the \CGA structure from the statement of \Cref{thm:Tor-CGA}
and discussing some consequences.}
\new{In general, he is right,
but with some additional assumptions---which are not apparently 
easy to verify in general---%
desirable properties can still be proven.}

\wrong{First, we show the product map does not depend on our choices.
\begin{theorem}\label{thm;well-defined}
	Given \WHCA maps and homotopies as
	in \eqref{eq:SHC-map-2}
	the product of \Cref{def:product}
	does not depend on the given homotopies
	of maps $\W\B(A \ox A) \lt \W\B X$
	and $\W\B(A \ox A) \lt \W\B Y$.
\end{theorem}
\begin{proof}
Recall that the use of the homotopy in \Cref{eq:full-product-v2}
is only in the last two steps,
as an application of \Cref{thm:homotopy-Tor-map}
with 
\[
R' = \W\B(A \ox A)\quad
R = \W\B A,			\quad
M' = \W\B(X \ox X),\quad
M = \W\B X,			\quad
N' = \W\B(Y \ox Y),\quad
N = \W\B Y\mathrlap.
\]
But then independence is immediate from \Cref{thm:homotopy-independence}.
\end{proof}
}

The desired properties making Munkholm's product as described 
in \Cref{sec:product} a \CGA follow in bijection with the 
defining properties of an \SHCA in \Cref{def:SHC}.
We subdivide the proof accordingly.

\begin{theorem}\label{thm:Tor-omnibus}
	Let {\WHCA}s $A$, $X$, $Y$ and \WHCA maps 
	$\B X \os\xi\from \B A \to \B Y$ be given.
	Suppose each \WHCA structure satisfies 
\benum

\item\label{thm:Tor-unit}
	the unitality condition \ref{def:SHC}.\ref{def:SHC-unit}.
	Then the product \eqref{eq:full-product-v2} 
	on $\Tor_A(X,Y)$
	is unital, with unit\footnote{\ 
		The $\h$ here are the units of the $\kk$-algebra structures on
		$\W\B Z$ for $Z \in \{A,X,Y\}$, unusually,
		not to be confused with the unit $\id \lt \B\W$
		of the bar--cobar adjunction.
		}
	\[
	\kk \simto \Tor_\kk(\kk,\kk) \xtoo{\Tor_{\h}(\h,\h)} 
	\Tor_{\W\B A}(\W\B X,\W\B Y)
 	\mathrlap.
	\]
\item\label{thm:Tor-comm}
	the commutativity condition \ref{def:SHC}.\ref{def:SHC-comm}.
	Then the product is commutative\new{,
	assuming an additional, not particularly transparent 
	or easily verified compatibility condition,
	to be stated in the course of the proof,
	on the various defining homotopies}.
\item\label{thm:Tor-assoc}
	the associativity condition \ref{def:SHC}.\ref{def:SHC-assoc}.
	Then the product is associative\new{,
	assuming an additional, not particularly transparent 
	or easily verified compatibility condition,
	to be stated in the course of the proof,
	on the various defining homotopies}.
	\eenum
\end{theorem}

\subsection{Unitality}
Unitality is easiest.

\begin{proof}[Proof of \Cref{thm:Tor-omnibus}.\ref{thm:Tor-unit}]
We prove that this map is a left unit, the right proof being symmetric.
Using the identifications $\kk \iso \B\kk \iso \W\B\kk$,
and following through our modified definition of the product,
we get the diagram
\begin{equation}\label{eq:unit-scheme}
\begin{aligned}
\xymatrix@C=1.25em@R=1.75em{
	\us{\W\B\kk}{\Tor}{}\, \ox \us{\W\B A}{\Tor}{}	\ar[rr]^{\mr{ext}}_\sim
									\ar[dd]_{\Tor_{\W\B\h} \ox \Tor_{\W\B{\id}}}&&
	\us{\W\B\kk\, \ox \W\B A}\Tor			\ar[dd]_{\Tor_{\W\B\h \ox \W\B{\id}}}
									&&	
	\us{\W(\B\kk \,\ox \B A)}\Tor				\ar[ll]^\sim_{\Tor_\g}
									\ar[rr]_\sim^{\Tor_{\W\nabla}} 
									\ar[dd]|{\Tor_{\W(\B\h \ox \B{\id})}}&&
	\us{\W\B(\kk \,\ox A)}\Tor				\ar[rr]_\sim
									\ar[dd]^{\Tor_{\W\B(\h \ox {\id})}}&&
	\us{\W\B A}\Tor						\ar@{=}[dd]				
														\\
	& &&&&
	&& 														\\	
	\us{\W\B A}{\Tor}{} \ox \us{\W\B A}{\Tor}{}	\ar[rr]_{\mr{ext}}&& 
	\us{\W\B A \ox \W\B A}\Tor						&& 
	\us{\W(\B A \ox \B A)}\Tor				\ar[ll]_\sim^{\Tor_\g}
									\ar[rr]^\sim_{\Tor_{\W\nabla}} 
&&
	\us{\W\B(A \ox A)}\Tor					\ar[rr]&&
	\us{\W\B A}\Tor	\mathrlap.							
}
\end{aligned}
\end{equation}
Commutativity of the first three squares follows from naturality,
in brief.
At length, the external product is functorial in all six of its entries,
giving the first square~\cite[p.~206]{cartaneilenberg}.
For the second and third,
by \Cref{def:shuffle}, note $\g$ 
and $\W\nabla$ are respectively
natural transformations
$\W({-} \ox {-}) \lt \W(-) \ox \W(-)$ and
$\W(\B{-} \ox \B{-}) \lt \W\B({-} \ox {-})$.

The last square 
obviously should commute as a result of the assumed condition
$\Phi\o \B(\h \ox \id) = \id$,
but to formally verify it
we require the six-square diagram \eqref{eq:Tor-DGA-homotopy-squares},
of which we display only the $A$-$X$ portion,
appending the vertical map $\Tor_{\W\B(\h \T {\id})}$:
\[
\xymatrix@R=3em{
	\W\B(\kk\,\ox A) 
		\ar[rr]^(.48){\W\B(\h \ox {\id})}
		\ar[d]_{\W(\id_{\B\kk} \T \xi)}
	&&
	\W\B(A \ox A) 
		\ar@{=}[r]
		\ar[d]_{\W(\xi \T \xi)}
	&
	\W\B(A \ox A) 
		\ar@{=}[r]
		\ar[d]
	&
	\W\B(A \ox A) 
		\ar[r]^(.6){\W\Phi}
		\ar[d]
	&	
	\W\B A	
		\ar[d]^{\W\xi}
	\\	
	\W\B(\kk\,\ox X) 
		\ar[rr]_(.48){\W\B(\h \ox {\id})}
	&&
	\W\B(X \ox X) 
		\ar[r]_(.55){\W\Phi}
	&
	\W\B X 
	&
	P \W\B X 
		\ar[l]\ar[r]
	&	
	\W\B X\mathrlap.
}
\]
That the leftmost square commutes is the same as stating 
the vertical map we call $\Tor_{\W\B(\h \T {\id})}$
exists in the first place, and follows from \Cref{thm:T-compose},
using the fact that $\B\h_A = \h_{\B A}$ and
$\B\h_X = \h_{\B X}$ are the coaugmentations.
Using $\Phi_X\o \B(\h_X \ox \id\mn_X) = \id_{\B X}$,
we may merge the first two squares.
Then we extend the commutative diagram to include the composite right homotopy 
$\W\B(\kk \ox A) \to \W\B(A \ox A) \to P\,\mn\W\B X$:
\[
\xymatrix@R=1.5em{
&&\W\B(\kk \,\ox A)
	\ar@{=}[lld]
	\ar[ld]
	\ar[rd]|{\W\B(\h \ox {\id})}
	\ar[dd]
	\ar@/^.25pc/[rrd]^{\sim}
\\
	\W\B(\kk\,\ox A) 
		\ar[r]_\sim
		\ar[dd]	&
	\W\B(A \ox A) 
		\ar@{=}[rr]|\hole
		\ar[dd]|(.75)\hole
	&&
	\W\B(A \ox A) 
		\ar[r]_(.57){\W\Phi}
		\ar[dd]|(.8)\hole
	&	
	\W\B A	
		\ar[dd]
	\\
&&P\,\mn\W\B X	
	\ar[dll]_(.4)\sim
	\ar[drr]^(.45){\pi_1}
	\ar[dl]^(.45){\pi_0}
	\ar@{=}[dr]
	\\	
	\W\B(\kk\,\ox X) 
		\ar[r]_(.55)\sim
	&
	\W\B X 
	&&
	P \W\B X 
		\ar[ll]^{\pi_0}\ar[r]_{\pi_1}
	&	
	\W\B X\mathrlap,
}
\]
where the map $\W\B(\kk\,\ox A) \lt \W\B A$
is the one we have been calling $\id_{\W\B A}$
under the identification 
using the condition $\Phi_A(\h_A \T \id\mn_A) = \id_{\B A}$.
Converting this to a commutative diagram of Tors, we get
\[
\xymatrix@C=.725em{
	&
	&
\us{\W\B(\kk \ox A)}\Tor (P\,\mn\W\B X)
	\ar[dll]_\sim
	\ar[dl]
	\ar[dr]
	\ar[drr]^\sim
	\\
\us{\W\B(\kk\,\ox A)}\Tor\big(\W\B(\kk\,\ox X)\big)
	\ar[r]
	&
\us{\W\B(A \ox A)}\Tor(\W\B X) 
	&&
\us{\W\B(A \ox A)}\Tor(P\,\mn\W\B X) 
	\ar[ll]^\sim
	\ar[r]
	&
\us{\W\B A}\Tor(\W\B X) \mathrlap,
}
\]
where the isomorphisms along the top become the identity
under the standard identifications.
\end{proof}

\subsection{Commutativity}

The proof of commutativity is more involved.

\begin{proof}[Proof of \Cref{thm:Tor-omnibus}.\ref{thm:Tor-comm}]
Assume also given  
{\WHCA}s $A'$, $X'$, $Y'$ and \WHCA maps 
$\B X' \from \B A' \to \B Y'$.
We write $\chi\: A \ox A' \lt A' \ox A$.
As in the proof of \Cref{thm:Tor-unit},
we trace through the definition of the product;
in the last square, we will assume
$A' = A$, $X' = X$, and $Y' = Y$
and finally use the homotopy-commutativity assumption
on the \WHCA structures.
\begin{equation}\label{eq:comm-scheme}
\begin{aligned}
\xymatrix@C=1.25em@R=1.75em{
	\us{\W\B A}{\Tor}{}\, \ox \us{\W\B A'}{\Tor}{}	\ar[rr]^{\mr{ext}}_\sim
									\ar[dd]_\chi&&
	\us{\W\B A\, \ox \W\B A'}\Tor			\ar[dd]_{\Tor_\chi}
									&&	
	\us{\W(\B A \,\ox \B A')}\Tor				\ar[ll]^\sim_{\Tor_\g}
									\ar[rr]_\sim^{\Tor_{\W\nabla}} 
									\ar[dd]|{\Tor_{\W\chi}}&&
	\us{\W\B(A \,\ox A')}\Tor				\ar[rr]_\sim
									\ar[dd]^{\Tor_{\W\B\chi}}&&
	\us{\W\B A'}\Tor						\ar@{=}[dd]				
														\\
	& &&&&
	&& 														\\	
	\us{\W\B A'}{\Tor}{} \ox \us{\W\B A}{\Tor}{}	\ar[rr]_{\mr{ext}}&& 
	\us{\W\B A' \ox \W\B A}\Tor						&& 
	\us{\W(\B A' \ox \B A)}\Tor				\ar[ll]_\sim^{\Tor_\g}
									\ar[rr]^\sim_{\Tor_{\W\nabla}} 
&&
	\us{\W\B(A' \ox A)}\Tor					\ar[rr]&&
	\us{\W\B A}\Tor	\mathrlap.							
}
\end{aligned}
\end{equation}
The vertical maps in the first two squares make sense 
by the naturality of the external product and $\g$.
The map we have written as $\Tor_{\W\B\chi}$
makes sense because Munkholm~\cite[Prop.~3.5]{munkholm1974emss} shows the following 
square commutes up to homotopy:
\eqnl{\label{eq:chi-nabla-square}
\xymatrix@C=4em@R=4.25em{
	\W\B(A \ox A') 
		\ar[r]^{\W(\xi \T \xi')}
		\ar[d]_{\W\B\chi}					&
	\W\B (X \ox X')
		\ar[d]	^{\W\B\chi}			
								\\
	\W\B(A' \ox A) 
		\ar[r]_{\W(\xi' \T \xi)}		
						&
	\W\B (X' \ox X)
	\mathrlap.
}
}
Thus we are forced to induce the map using the technique of
\Cref{thm:homotopy-Tor-map}.

The commutativity of the external product square 
appears in Cartan--Eilenberg~\cite[Prop.~X.2.1]{cartaneilenberg};
the extra sign in their expression is implicit in our definition of the interchange map.
The next two squares say, roughly, that the shuffle maps $\gamma$ and $\nabla$
of \Cref{def:shuffle} are commutative. 
On prepending the tautological twisting cochain, for $\g$,
and postpending it, for $\nabla$, this boils down to the equation
chains
\eqn{
\chi\gamma t_{\B A \ox \B A'} &= 
\chi (t_{\B A} \ox \eta_{\W\B A'}\e_{\B A'} + \eta_{\W\B A}\e_{\B A} \ox t_{\B A'}) 
\\	&=
(t_{\B A'} \ox \eta_{\W\B A}\e_{\B A} + \eta_{\W\B A'}\e_{\B A'} \ox t_{\B A}) \chi
	=
\gamma t_{\B A' \ox \B A} \chi
	=
\gamma \, \W\chi\, t_{\B A \ox \B A'}\mathrlap,\\
t^{A' \ox A}\nabla \chi &=
(t^{A'} \ox \eta_{A}\e_{\B A} +  \eta_{A'}\e_{\B A'} \ox t^A) \,\chi 
\\	&=
\chi (t^A \ox \eta_{A'}\e_{\B A'} +  \eta_{A}\e_{\B A} \ox t^{A'})	
	=
\chi t^{A \ox A'}\nabla
	=
t^{A' \ox A}\,\B\chi \, \nabla\mathrlap,
}
where we have used naturality of the tautological twisting cochains. 
This works fully for the $\gamma$ square, 
but there is something to check for the $\W\nabla$ square because
the Tor is defined in steps.
Let $h\: \W\B(A \ox A') \lt P\,\mn\W\B(X' \ox X)$
be a right homotopy witnessing 
the homotopy making \eqref{eq:chi-nabla-square} commute.
Note that this square commutes on the nose if we precompose
$\W\nabla\: \W(\B A \ox \B A') \lt \W\B(A \ox A')$,
for by the preceding equations 
and \Cref{thm:tensor-psi},
we have
\[
	\B\chi
	\,(\xi \T \xi')
	\,\nabla
		=
	\B\chi
	\,\nabla
	\,(\xi \ox \xi')
		=
	\nabla
	\,\chi
	\,(\xi \ox \xi')
		=
	\nabla
	\,(\xi' \ox \xi)
	\,\chi
		=
	(\xi' \T \xi)
	\,\nabla 
	\,\chi
		=
	(\xi' \T \xi)
	\,\B\chi
	\,\nabla
		\mathrlap.
\]
Thus the following diagram commutes:
\[
\xymatrix@R=1em@C=-1em{
	\W(\B A \ox \B A')
\ar@{=}[dd] 
\ar@[jred][rd]^(.675){\textcolor{jblue}{\W(\xi \ox \xi')\ \ }}
\ar@{=}[rrr]
	&&&	
	\W(\B A \ox \B A')
\ar@{=}[dd]|\hole  
\ar@{-->}@[jcmp][rd]
\ar@{=}[rrr]
	&&&	
	\W(\B A \ox \B A')
\ar@{=}[dd]|\hole  
\ar@{-->}@[jcmp][rd]
\ar[rrr]^{\textcolor{jblue}{\W\nabla}}
	&&&	\W\B (A \ox A')
\ar@{=}[dd]|\hole  
\ar@[jred][rd]^(.65){\textcolor{jblue}{\W(\xi' \T \xi')}}
	\\
	&\textcolor{jX2}{\W(\B X \ox \B X')}
\ar@[jred][dd]_(.275){\textcolor{jblue}{\W\chi}} 
\ar@[jX2][rrr]^(.375){\textcolor{jblue}{\W\nabla}}
	&&&\textcolor{jX2}{\W\B(X \ox X')}
\ar@[jred][dd]
\ar@[jX2]@{=}[rrr]
	&&&\textcolor{jX2}{\W\B(X \ox X')}
\ar@[jred][dd]
\ar@[jX2]@{=}[rrr]
	&&&\textcolor{jX2}{\W\B(X \ox X')}
\ar@[jred][dd]^{\textcolor{jblue}{\W\B\chi}}
	\\
	\W(\B A \ox \B A')
\ar@{=}[rrr]|!{[ur];[dr]}\hole
\ar@{=}[dd]
\ar@{-->}@[jcmp][rd]  
	&&&\W(\B A \ox \B A')
\ar@{=}[rrr]|!{[ur];[dr]}\hole
\ar@{=}[dd]|\hole  
\ar@{-->}@[jcmp][rd]  
	&&&\W(\B A \ox \B A')
\ar[rrr]|!{[ur];[dr]}\hole
\ar@{=}[dd]|\hole  
\ar@{-->}@[jcmp][rd]  
	&&&\W\B (A \ox A')
\ar@{=}[dd]|\hole 
\ar@{-->}@[jcmp][rd] 
	\\
	&\textcolor{jX2p}{\W(\B X' \ox \B X)}
\ar@[jX2p][rrr]
	&&&\textcolor{jX2p}{\W\B(X' \ox X)}
\ar@[jX2p]@{=}[rrr]
	&&&\textcolor{jX2p}{\W\B(X' \ox X)}
\ar@[jX2p]@{=}[rrr]
	&&&\textcolor{jX2p}{\W\B(X' \ox X)}
	\\
	\W(\B A \ox \B A')
\ar@[jred][dd]_{\textcolor{jblue}{\W\chi}}
\ar@{-->}@[jcmp][rd]
\ar@{=}[rrr]|!{[ur];[dr]}\hole 
	&&&\W(\B A \ox \B A')
\ar@[jred][dd]|\hole 
\ar@{-->}@[jcmp][rd]
\ar@{=}[rrr]|!{[ur];[dr]}\hole 
	&&&\W(\B A \ox \B A')
\ar@[jred][dd]|\hole 
\ar@{-->}@[jcmp][rd]_(.375){\textcolor{jblue}{h\,\W\nabla}}
\ar[rrr]|!{[ur];[dr]}\hole 
	&&&\W\B (A \ox A')
\ar@{~>}@[jhmt][rd]_(.35){\textcolor{jblue}{h}}
\ar@[jred][dd]|\hole 
	\\
	&\textcolor{jX2p}{\W(\B X' \ox \B X)}
\ar@[jX2p]@{=}[uu]
\ar@[jX2p]@{=}[dd]
\ar@[jX2p][rrr]
	&&&\textcolor{jX2p}{\W\B(X' \ox X) }
\ar@[jX2p]@{=}[uu]
\ar@[jX2p]@{=}[dd]
\ar@[jX2p]@{<-}[rrr]
	&&&\textcolor{jX2p}{P\,\mn\W\B(X' \ox X) }
\ar@[jX2p][uu]
\ar@[jX2p]@{=}[dd]
\ar@[jX2p]@{=}[rrr]
	&&&	\textcolor{jX2p}{P\,\mn\W\B(X' \ox X)}
\ar@[jX2p][uu]_{\textcolor{jblue}{\pi_0}}
\ar@[jX2p][dd]^{\textcolor{jblue}{\pi_1}}
	\\
		\textcolor{jQ2}{\W(\B A' \ox \B A)}
\ar@[jred][dr]_(.35){\textcolor{jblue}{\W(\xi' \ox \xi)\ } }
\ar@[jQ2]@{=}[rrr]|!{[dr];[ur]}\hole 
	&&&	\textcolor{jQ2}{\W(\B A' \ox \B A)}
\ar@{-->}@[jcmp][dr]
\ar@[jQ2]@{=}[rrr]|!{[dr];[ur]}\hole 
	&&&	\textcolor{jQ2}{\W(\B A' \ox \B A)}
\ar@{-->}@[jcmp][dr]_(.325){\textcolor{jblue}{h\,\W\nabla\,\W\chi}}
\ar@[jQ2][rrr]|!{[dr];[ur]}\hole 
	&&&	\textcolor{jQ2}{\W\B (A' \ox A)}
\ar@[jred][rd]_(.3){\textcolor{jblue}{\W(\xi' \T \xi)} }
	\\
	&\textcolor{jX2p}{\W(\B X' \ox \B X)}
\ar@[jX2p][rrr]_{\textcolor{jblue}{\W\nabla}}
	&&&\textcolor{jX2p}{\W\B (X' \ox X)}
\ar@[jX2p]@{<-}[rrr]_{\textcolor{jblue}{\pi_0}}
	&&&\textcolor{jX2p}{P\,\mn\W\B (X' \ox X)}
\ar@[jX2p][rrr]_{\textcolor{jblue}{\pi_1}}
	&&&\textcolor{jX2p}{\W\B (X' \ox X)}
\mathrlap.
}
\]

There is less going on here than meets the eye,
and the color-coding of objects by quasi-isomorphism class 
is hoped to make the few changes somewhat easier to follow.
Arrows defined as composites are grey and dashed,
other non-quasi-isomorphisms are red, 
and the one right homotopy that is not a composite is gold and wavy.
The left face just expresses in a long-winded way
that $\W\chi\, \W(\xi \ox \xi') = \W(\xi' \ox \xi)\,\W\chi$,
and the next vertical face comes from postcomposing $\W\nabla$.
The right face is the three-square diagram inducing the map
we abusively called $\Tor_{\W\B\chi}$.
The third vertical face is the only interesting one.
Its top square is inherited from the second vertical face,
its middle square expresses that 
$\pi_0\,h\,\W\nabla$ is the function in the preceding long display,
and the bottom square commutes since $\chi \o \chi = \id$.
The map from the third face back to the second is trivial
except for the bottom two maps in front, which are $\pi_0$,
which is possible because 
$\pi_0\,h\,\W\nabla$ 
and 
$\pi_1\,h\,\W\nabla$ 
are equal by the long display
and again because $\chi \o \chi = \id$.
The map from the third face to the right face
makes sense entirely by definition in the top two cubes;
in the bottom cube,
the back commutes since $\B\chi \o \nabla = \nabla \o \chi$,
the front commutes trivially,
and the bottom commutes from the long display and yet again 
because $\chi \o \chi = \id$.

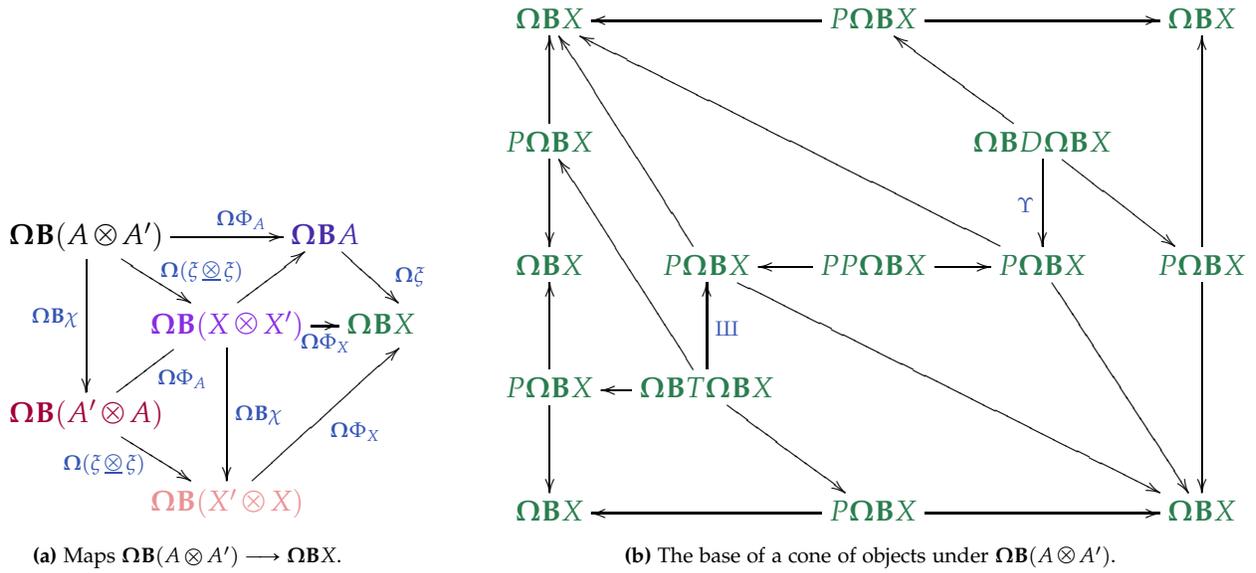
\begin{figure}
	\centering
	\begin{subfigure}{0.3\textwidth}
		\centering
$
\begin{aligned}
\xymatrix@R=1.5em@C=-1em{	
	\W\B(A \ox A')
	\ar@[jred][rr]^(.65){\textcolor{jblue}{\W\Phi_A}}
	\ar@[jred][dr]^(.6){\textcolor{jblue}{\W(\xi \T \xi)}}
	\ar@[jred][dd]_(.45){\textcolor{jblue}{\W\B\chi}}
	&&
	\textcolor{jA}{\W\B A}		
	\ar@[jred][dr]^(.65){\textcolor{jblue}{\W\xi}}	
	\\
	&
	\textcolor{jX2}{\W\B(X \ox X')}
	\ar@[jred][rr]_(.5){\textcolor{jblue}{\W\Phi_X}}
	\ar@[jred][dd]^(.5){\textcolor{jblue}{\W\B\chi}}	
	&&
	\textcolor{jgreen}{\W\B X}	
	{\phantom{({}\ox X')}}
	\\
	\textcolor{jQ2}{\W\B(A' \ox A)}
	\ar@[jred][uurr]	_(.3){\textcolor{jblue}{\W\Phi_A}}
	|(.4)\hole|(.45)\hole|(.5)\hole|(.55)\hole|(.6)\hole
	\ar@[jred][dr]_(.35){\textcolor{jblue}{\W(\xi \T \xi)}}
	&&&&&&	
	\\
	&
	\textcolor{jX2p}{\W\B(X' \ox X)}
	\ar@[jred][uurr]_{\textcolor{jblue}{\W\Phi_X}}
}
\end{aligned}
$
\caption{Maps $\W\B(A \otimes A') \lt \W\B X$.}
\label{fig:comm-homotopy-prism}
	\end{subfigure}
	\hfill%
\begin{subfigure}{0.6\textwidth}
	\centering$
		\xymatrix@C=0.25em@R=3em{
			\textcolor{jgreen}{	\W\B X } &&&& 
			\textcolor{jgreen}{	P\,\mn\W\B X } \ar@[jgreen][llll]\ar@[jgreen][rrrr]&&&&
			\textcolor{jgreen}{	\W\B X }
			\\
			\textcolor{jgreen}{	P\,\mn\W\B X } \ar@[jgreen][d]\ar@[jgreen][u]	&&&&&&
			\textcolor{jgreen}{	\W\B D \W\B X } \ar@[jgreen][llu]\ar@[jgreen][d]_{\textcolor{jblue}{\Upsilon}}\ar@[jgreen][drr]
			\\
			\textcolor{jgreen}{	\W\B X } &&
			\textcolor{jgreen}{	P\,\mn\W\B X }	\ar@[jgreen][lluu]\ar@[jgreen][ddrrrrrr]&& 
			\textcolor{jgreen}{	P \mnn P\,\mn\W\B X }\ar@[jgreen][ll]\ar@[jgreen][rr]&&
			\textcolor{jgreen}{	P\,\mn\W\B X }
			\ar@[jgreen][uullllll]\ar@[jgreen][ddrr]&&
			\textcolor{jgreen}{	P\,\mn\W\B X } \ar@[jgreen][dd]\ar@[jgreen][uu]
			\\
			\textcolor{jgreen}{	P\,\mn\W\B X } \ar@[jgreen][d]\ar@[jgreen][u]		&&
			\textcolor{jgreen}{	\W\B T\W\BX } \ar@[jgreen][u]_{\textcolor{jblue}{\Sha}}\ar@[jgreen][uull]\ar@[jgreen][ll]\ar@[jgreen][drr]
			\\
			\textcolor{jgreen}{	\W\B X } &&&& 
			\textcolor{jgreen}{	P\,\mn\W\B X } \ar@[jgreen][rrrr]\ar@[jgreen][llll]&&&&
			\textcolor{jgreen}{	\W\B X }
		}
	$
	\caption{The base of a cone of objects under $\W\B(A \ox A')$.}
	\label{fig:comm-plug}
\end{subfigure}
	\caption{Auxiliary diagrams for the commutativity argument.}
	\label{fig:comm-extra-figs}
\end{figure}

We can finally consider the $\Phi$ triangle at the right of \eqref{eq:comm-scheme}.
The maps of Tors are induced by the homotopy-commutative 
squares of the prism in \Cref{fig:comm-homotopy-prism}.
There are five edge-paths from $\W\B(A \ox A')$ to $\W\B X$,
all of lengths one or two, and if we say 
two paths \emph{neighbor} one another if they
together bound a face of the prism,
then each path has two neighbors. 
The top and bottom face homotopies,
which are the same, are prescribed 
by the fact $\xi\: \B A \lt \B X$ is a \WHCA map,
\wrong{although we have seen we have some flexibility
in which homotopies we use to define the product,}
and the left face homotopy is $h$ 
from the preceding argument.
The right homotopies representing these homotopies making 
\Cref{fig:comm-homotopy-prism} commute
can be expanded
to give the cubical \Cref{fig:commutativity-cube}.
\begin{figure}
\centering
$
\xymatrix@C=-.25em@R=1em{
\W\B(A \ox A')		\ar@[jred][dr]^(.6){\textcolor{jblue}{\ \W(\xi \T \xi)}}
				\ar@{=}[rr]
				\ar@{=}[dddd]
				\ar@/_2.5pc/@{=}[ddrr]|(.59)\hole|(.75)\hole|(.775)\hole
&& \W\B(A \ox A')	\ar@{=}[rr]\ar@{-->}@[jcmp][dr]
&& \W\B(A \ox A')	\ar@[jred][rr]^(.5775){\textcolor{jblue}{\W\Phi_{\!\smash{A}}\!}}
&& \textcolor{jA}{\W\B A} 			\ar@[jred][dr]^(.55){\mn\!\!\textcolor{jblue}{\W\xi}}
\\
& \textcolor{jX2}{\W\B(X \ox X')}		\ar@[jred][dddd]_(.45){\textcolor{jblue}{\W\B\chi}}
				\ar@[jred][rr]^(.35){\textcolor{jblue}{\W\Phi_{X}}} 
				\ar@/_2pc/@{~>}@[jhmt][ddrr]
&& \textcolor{jgreen}{\W\B X} 		
&& \textcolor{jgreen}{P\,\mn\W\B X} 		\ar@{<~}@[jhmt][ul]
				\ar@[jgreen][ll]
				\ar@[jgreen][rr]
&& \textcolor{jgreen}{\W\B X}		
\\
&&\W\B(A \ox A')		\ar@{=}[uu]|\hole
				\ar@{=}[dd]|\hole
&&
&&
\\
&
&&\textcolor{jgreen}{P\,\mn\W\B X} 		
				\ar@{<~}@[jhmt][ul]
				\ar@[jgreen][uu]
				\ar@[jgreen][dd]
&&
&& 
\\
\W\B(A \ox A')		\ar@{=}[dd]
				\ar@{=}[rr]|\hole
				\ar@{-->}@[jcmp][dr]
&&\W\B(A \ox A')	
				\ar@{=}[dd]|\hole
				\ar@{~>}@[jhmt][dr]
&&
&&\textcolor{jA}{P\,\mn\W\B A}	\ar@[jred][dr]_(.45){\ \textcolor{jblue}{P\,\mn\W\xi}\!\!\!}
			\ar@[jA][dddd]
			\ar@[jA][uuuu]|(.75)\hole
\\
&\textcolor{jX2p}{\W\B(X' \ox X)}		\ar@[jred][rr]^(.42){\textcolor{jblue}{\W\Phi_X}}
&&\textcolor{jgreen}{\W\B X}
&&
&&\textcolor{jgreen}{P\,\mn\W\B X}		
				\ar@[jgreen][uuuu]
				\ar@[jgreen][dddd]
\\
\W\B(A \ox A')		\ar@[jred][dd]_{\textcolor{jblue}{\W\B\chi}} 
				\ar@{=}[rr]|\hole
&&\W\B(A \ox A')		\ar@{-->}@[jcmp][dd]|\hole
&&   
&&
\\
&\textcolor{jX2p}{P\,\mn\W\B(X' \ox X)}	
				\ar@{<~}@[jhmt][ul]
				\ar@[jX2p][uu]
				\ar@[jX2p][dd]
				\ar@[jred][rr]^(.44){\textcolor{jblue}{P\,\mn\W\Phi_{\mn X}}}
&&\textcolor{jgreen}{P\,\mn\W\B X} 		\ar@{<--}@[jcmp][ul]
				\ar@[jgreen][uu]
				\ar@[jgreen][dd]
&& 
&& 
\\
\textcolor{jQ2}{\W\B(A' \ox A)}		
				\ar@[jred][dr]_(.41){\mathllap{\textcolor{jblue}{\W(\xi \T \xi)}}}
				\ar@{=}@[jQ2][rr]|\hole
&&\textcolor{jQ2}{\W\B(A' \ox A)}		
				\ar@{=}@[jQ2][rr]|\hole
				\ar@{-->}@[jcmp][dr]
&& \textcolor{jQ2}{\W\B(A' \ox A)}	\ar@[jred][rr]^(.575){\textcolor{jblue}{\W\Phi_A}}
&& \textcolor{jA}{\W\B A}			\ar@[jred][dr]_(.3){\textcolor{jblue}{\W\xi}\!} 
\\
&\textcolor{jX2p}{\W\B(X' \ox X)}		
				\ar@[jred][rr]_{\textcolor{jblue}{\W\Phi_X}}
&&\textcolor{jgreen}{\W\B X} 
&& \textcolor{jgreen}{P\,\mn\W\B X} 		\ar@{<~}@[jhmt][ul]
				\ar@[jgreen][ll]
				\ar@[jgreen][rr]
&& \textcolor{jgreen}{\W\B X}
}
$
\caption{The cube diagram for commutativity.}
\label{fig:commutativity-cube}
\end{figure}
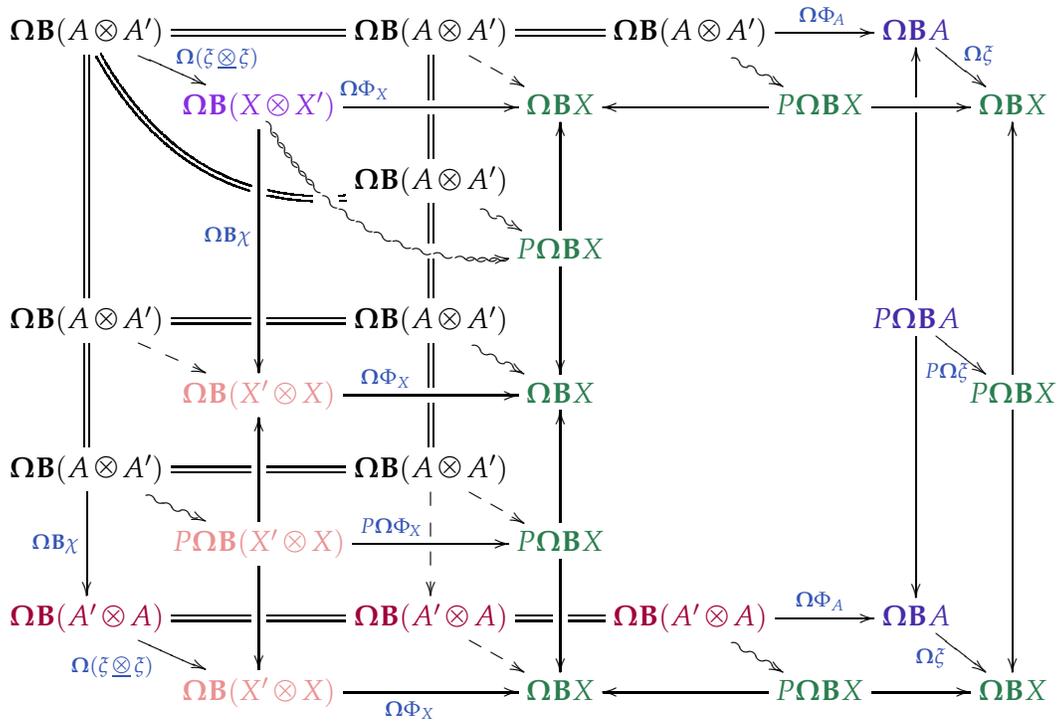
In this figure, the objects are color-coded by homotopy type,
primitive arrows are red,
right homotopies wavy gold, 
and composites defined
so as to make squares commute grey.
The right face will induce the identity isomorphism on $\Tor_{\W\B A}$.
The three small interior cubes and triangular prism on the left
commute by definition.

By assumption, there is a \DGC homotopy between $\Phi_A$ and $\Phi_A \o \B\chi$,
inducing a homotopy from $\W \Phi_A$ to $\W\Phi_A \o \W\B\chi$,
which is represented by a right homotopy
$\W\B(A \ox A') \lt P\,\mn\W\B A$. 
Postcomposing $P\,\mn\W\xi$ gives a right homotopy
$\W\B(A \ox A') \lt P\,\mn\W\B X$.
The composition of the homotopy represented by the 
top face in the large cube and this homotopy
is a homotopy from $\W\Phi_X \o \W(\xi\T\xi)$ to $\W\xi \o \W\Phi_A \o \W\B\chi$,
as is the composition of the three homotopies
represented by the right homotopies along the lower left.
The right homotopies representing these two compositions
are themselves homotopic
by \Cref{thm:homotopies-homotopic} 
\new{(not that this is helpful or relevant)},
so one finds a \wrong{commutative} diagram of objects
receiving maps from $\W\B(A \ox A')$
as in \Cref{fig:comm-plug}. 
\new{The faces of the cone including the vertex $\W\B(A \ox A')$
are indeed commutative by assumption, but the base \Cref{fig:comm-plug}
is not; this is is precisely the same issue that the purported proof of 
\Cref{thm:homotopy-independence} faced,
which is repaired in \Cref{thm:homotopy-independence-mk2},
and so the necessary condition to repair this proof is that 
the right homotopy 
$H\: \W\B(A \ox A') \lt PP\W\B X$
in the middle of the cone be \emph{endpoint-fixing}
in the sense that the composites 
$P\pi_0 \o H$ and $P\pi_1 \o H$ 
factor respectively as $\z \o \W\Phi_X \o \W(\xi \T \xi)$
and $\z \o \W\xi \o \W\Phi_A$,
where $\z\: \W\B X \lt P \W\B X$ is the natural map
defined in \Cref{def:P}.
The cone we find we need is thus not one whose base 
is \Cref{fig:comm-plug},
but instead \Cref{fig:true-base-comm}.
}

\new{
We observe that there is no obvious criterion to determine when
this is possible.
}

\begin{figure}
	\centering
$
\xymatrix@C=-.25em@R=2.5em{
	\textcolor{jgreen}{\W\B X} 
	\ar@[jgreen][drrrr]_{\id}
	\ar@[jgreen][ddrrr]|{\z}
	\ar@[jgreen][dddrr]^{\id}
	&&
	&&
	&&
	\textcolor{jgreen}{P\,\mn\W\B X} 
	\ar@[jgreen][llllll]_{\pi_0}
	\ar@[jgreen][rr]^{\pi_1}
	&&
	\textcolor{jgreen}{\W\B X}
	\\
	&&&&
	\textcolor{jgreen}{\W\B X}
	&&
	\textcolor{jgreen}{\W\B D\W\B X} 
	\ar@[jgreen][u]
	\ar@[jgreen][dl]_(.45){\textcolor{jblue}{\Upsilon}}
	\ar@[jgreen][ddrr]
	\\
	&&&
	\textcolor{jgreen}{P\W\B X}
	\ar@[jgreen][dl]^{\pi_0}
	\ar@[jgreen][ur]_{\pi_1}
	&&
	\textcolor{jgreen}{P\W\B X} 
	\ar@[jgreen][ul]^{\pi_0}
	\ar@[jgreen][dr]_{\pi_1}
	\\
	\textcolor{jgreen}{P\,\mn\W\B X}  
	\ar@[jgreen][ddd]_{\pi_1}
	\ar@[jgreen][uuu]^{\pi_0}
	&&
	\textcolor{jgreen}{\mn\W\B X}
	&& 
	\textcolor{jgreen}{P \mnn P\,\mn\W\B X}
	\ar@[jgreen][dl]|{\pi_0}
	\ar@[jgreen][ur]|{\pi_1}
	\ar@[jgreen][ul]|{P\pi_0}
	\ar@[jgreen][dr]|{P\pi_1}
	&&
	\textcolor{jgreen}{\W\B X} 
	&&
	\textcolor{jgreen}{P\,\mn\W\B X}
	\ar@[jgreen][ddd]^{\pi_1}
	\ar@[jgreen][uuu]_{\pi_0}
	\\
	&&&
	\textcolor{jgreen}{P\,\mn\W\B X} 
	\ar@[jgreen][ul]_{\pi_0}
	\ar@[jgreen][dr]^{\pi_1}
	&&
	\textcolor{jgreen}{P\,\mn\W\B X}
	\ar@[jgreen][dl]_{\pi_0}
	\ar@[jgreen][ur]^{\pi_1}
	\\
	&&
	\textcolor{jgreen}{\W\B T\W\B X} 
	\ar@[jgreen][ur]_(.45){\textcolor{jblue}{\Sha}}
	\ar@[jgreen][d]
	\ar@[jgreen][uull]
	\ar@[jgreen][drrrr]
	&&
	\textcolor{jgreen}{\W\B X} 
	\\
	\textcolor{jgreen}{\W\B X} 
	&&
	\textcolor{jgreen}{P\,\mn\W\B X} 
	\ar@[jgreen][ll]^{\pi_0}
	\ar@[jgreen][rr]_{\pi_1}
	&&
	\textcolor{jgreen}{\W\B X} 
	&& 
	\textcolor{jgreen}{P\,\mn\W\B X} 
	\ar@[jgreen][ll]^{\pi_0}
	\ar@[jgreen][rr]_{\pi_1}
	&&
	\textcolor{jgreen}{\W\B X}
	\ar@[jgreen][ullll]_{\id}
	\ar@[jgreen][uulll]|{\z}
	\ar@[jgreen][uuull]^{\id}
}
$
\caption{\new{The true base for the cone at the end of the commutativity proof.}}
	\label{fig:true-base-comm}
\end{figure}

Noting that right homotopy $\W\B(A \ox A') \lt P\,\mn \W\B X$ on the right edge
coming from $\W \Phi_A \hmt \W\Phi_A \o \W\B\chi$
factors through $P\,\mn\W\B A$
and the right homotopy on the bottom factors through $\W\B(A' \ox A)$,
we may plug this diagram into the large rectangle in the front face
of the previous cube and take Tor to obtain a large commutative diagram.
The map induced on the right edge is the identity map of $\Tor_{\W\B A}(\W B X)$
since the projections $\pi_1$ and $\pi_0$ induce the same map in cohomology,
so this completes the final square in \eqref{eq:comm-scheme}
and with it the proof.
\end{proof}

\subsection{Associativity}

The associativity proof is again more involved.\footnote{\ 
	To follow this proof carefully,
	it may be helpful to first absorb \Cref{sec:ring-map},
	as broadly similar arguments are presented
	more expansively there.
	}

\begin{proof}[Proof of \Cref{thm:Tor-omnibus}.\ref{thm:Tor-assoc}]
	The template is Figure \ref{fig:assoc},
	in which we show each square commutes.
%
\begin{figure}
\resizebox{149mm}{!}{
\xymatrix@C=0.5em{
	\Tr{\W\B A}	\ox \Tr{\W\B A} \ox \Tr{\W\B A}
		\ar[rr]^(.5){\id \ox \ext}
		\ar[dd]_{\ext \ox \id}
&&
	\Tr{\W\B A} \ox \Tr{(\W\B A)\ot}
		\ar[dd]_{\ext}
		\ar@{<-}[rr]^{\id \ox \Tr{\g}}
&&
	\Tr{\W\B A} \ox \Tr{\W(\B A)\ot}
		\ar[dd]_{\ext}
		\ar[rr]^{\id \ox \Tr{\W\nabla}}
&&
	\Tr{\W\B A} \ox \Tr{\W\B (A\ot)}
		\ar[dd]_{\ext}
		\ar[rr]^{\id \ox \Tr{\W\Phi}}
&&
	\Tr{\W\B A} \ox \Tr{\W\B A}
		\ar[dd]_{\ext}
\\
&
	\mr{I}
&&
	\mr{II}
&&
	\mr{III}
&&
	\mr{IV}
&
\\
	\Tr{(\W\B A)\ot} \ox \Tr{\W\B A}
		\ar[rr]^(.55){\ext}
		\ar@{<-}[dd]_{\Tr{\g} \ox \id}
&&
	\Tr{\W\B A \ox \W\B A \ox \W\B A}
		\ar@{<-}[rr]^{\Tr{\id \ox \g}}		
		\ar@{<-}[dd]_{\Tr{\g \ox \id}}
&&
	\Tr{\W\B A \ox \W(\B A\ot)}
		\ar@{<-}[dd]_{\Tr{\g}}
		\ar[rr]^{\Tr{\id \ox \W\nabla}}
&&
	\Tr{\W\B A \ox \W\B (A\ot)}
		\ar@{<-}[dd]_{\Tr{\g}}
		\ar[rr]^{\Tr{\id \ox \W\Phi}}
&&
	\Tr{(\W\B A)\ot}
		\ar@{<-}[dd]_{\Tr{\g}}
\\
&
\mr{II}'
&&
\mr{V}
&&
\mr{VI}
&&
\mr{VII}
&
\\
	\Tr{\W(\B A)\ot} \ox \Tr{\W\B A} 
			\ar[rr]^(.5){\ext}
			\ar[dd]_{\Tr{\W\nabla} \ox \id}
&&
	\Tr{\W(\B A\ot) \ox \W\B A}
		\ar@{<-}[rr]^{\Tr{\g}}
		\ar[dd]_{\Tr{\W\nabla \ox \id}}
&&
	\Tr{\W(\B A)\ott}
		\ar[rr]^{\Tr{\W(\id \ox \nabla)}}
		\ar[dd]_{\Tr{\W(\nabla \ox \id)}}
&&
	\Tr{\W(\B A \ox\B (A\ot))}
		\ar[rr]^{\Tr{\W(\id \ox \Phi)}}
		\ar[dd]_{\Tr{\W\nabla}}
&&
	\Tr{\W(\B A)\ot}
		\ar[dd]_{\Tr{\W\nabla}}
\\
&
\mr{III}'
&&
\mr{VI}'
&&
\mr{VIII}
&&
\mr{IX}
&
\\
	\Tr{\W\B(A\ot)} \ox \Tr{\W\B A} 
		\ar[rr]^(.5){\ext}
		\ar[dd]_{\Tr{\W\Phi} \ox \id}
&&
	\Tr{\W\B(A \ot) \ox \W\B A}
		\ar@{<-}[rr]^{\Tr{\g}}
		\ar[dd]_{\Tr{\W\Phi \ox \id}}
&&
	\Tr{\W(\B(A\ot) \ox \B A)}
		\ar[rr]^{\Tr{\W\nabla}}
		\ar[dd]_{\Tr{\W(\Phi \ox \id)}}
&&
	\Tr{\W\B(A\ott)}
		\ar[rr]^{\Tr{\W(\id \T \Phi)}}
		\ar[dd]_{\Tr{\W(\Phi \T \id)}}
&&
	\Tr{\W\B(A\ot)}
		\ar[dd]_{\Tr{\W\Phi}}
\\
&
\mr{IV}'
&&
\mr{VII}'
&&
\mr{IX}'
&&
\mr{X}{\phantom{\mr{I}}}
&
\\
	\Tr{\W\B A} \ox \Tr{\W\B A} 
			\ar[rr]^(.5){\ext}
&&
	\Tr{(\W\B A)\ot}
		\ar@{<-}[rr]^{\Tr{\g}}
&&
	\Tr{\W(\B A)\ot}
		\ar[rr]^{\Tr{\W\nabla}}
&&
	\Tr{\W\B(A\ot)}
		\ar[rr]^{\Tr{\W\Phi}}
&&
	\Tr{\W\B A}
    }
}
\caption{The associativity diagram.}\label{fig:assoc}
\end{figure}
We can be brief about the proofs the squares
from the first three columns commute, 
which mostly involve only naturality and functoriality, 
and only deal with the unprimed labels.
\benum
\item[I:]
	The associativity of the external product is classical; \emph{cf.}~Cartan--Eilenberg%
	~\cite[p.~206]{cartaneilenberg}.
\item[II:]
	This is the naturality of the external product in the second three variables~\cite[XI.2.1]{cartaneilenberg}.
\item[III:]	
	This too is the naturality of the external product in the second three variables.
\item[V:]
	This follows from naturality of $\gamma$ and 
	the equation $(\id \ox \g)\g = (\g \ox \id)\g$.
	To see this, it is enough to precompose the tautological twisting cochain
	$t\: (\B A)\ott \lt \W(\B A)\ott$
	and expand using the definition in \Cref{def:shuffle}.
\item[VI:]
	This follows because $\g$ is a natural transformation 
	$\W({-} \ox {-}) \lt \W({-}) \ox \W({-})$.
\item[VIII:]
		This follows from functoriality of $\W$,
		then naturality of $\nabla$ and 
	the equation $\nabla(\id \ox \nabla) = \nabla(\nabla \ox \id)$,
	whose proof is dual to that in [V].
\eenum

The other squares require their own diagrams.
\benum
\item[IV:]
Commutativity of the square 
follows from naturality of the external product in the last
three variables (the first three variables fixed as $\W\B X \from \W\B A \to \W\B Y$), 
applied to the diagram
\begin{equation}\label{eq:IV-diagram}
\begin{aligned}
	\xymatrix@C=1em@R=2em{
		\W\B(A \ox A)	\ar[r]\ar[d]
		&	\W\B(A \ox A)	\ar[r]\ar[d]
		&	\W\B(A \ox A)	\ar[d]^{\defm{\wt h}}\ar[r]
		&\W\B A					\ar[d]
		\\
		\W\B (X \ox X)	\ar[r]
		&	\W\B X					\ar[r]
		& P\,\mn\W\B X 				\ar[r]
		& \W\B X	
	}
\end{aligned}
\end{equation}
and the symmetric diagram with $X$ replaced with $Y$. 
The map we call $\wt h$
can be obtained
in various (equivalent) ways from 
the assumed \DGC homotopy
between \DGC maps $\B(A\ox A) \lt \B X$.
For our purposes it will be most convenient
to transpose this to a \DGA homotopy between two \DGA maps 
$\B(A\ox A) \lt X$
using \Cref{thm:homotopy-adjunction}, 
then represent that as a right homotopy
$\defm{h^P}\: \W\B(A \ox A) \lt PX$.
We write $\defm {h^\dagger}$ for the composite
$\B h^P \o \h\: \B(A\ox A) \to \B\W\B(A \ox A) \to \B P X$.
Then $\W h^\dagger$ is the map 
$(h^P)^\#\: \W\B(A \ox A) \lt \W\B PX$ induced up from $h^P$ 
by \Cref{thm:induction},
and we finally define $\wt h$ to be the composite
$Z \o \W h^\dagger$ with the natural map $Z\: \W\B PX \lt P\,\mn\W\B X$
of \Cref{thm:P-W}.
\eenum

We simplify life by referring only to the $A$-$X$ side
in the remaining squares.

\bs

\nd\emph{A replacement.}
To proceed in the diagram, we will need to replace
$\Tor_{\W\B A \ox \W\B(A\ot)}(\W B X \ox \textcolor{red}{P}\,\mn \W \B X)$
with 
$\Tor_{\W\B A \ox \W\B(A\ot)}(\W \B X \ox \W \B \textcolor{red}PX )$.
Recall that we defined the right homotopy $\wt h$
as the composite of the
\DGA quasi-isomorphism $Z\: \W\B PX \lt P \W \B X$
and a map $\W h^\dagger\:\W\B(A\ot) \lt \W \B  P X$.
%
Factoring $\wt h$ in this way in \eqref{eq:IV-diagram}, 
we get 
\[
\xymatrix@C=1em@R=2em{
	\W\B(A \ox A)	\ar[r]\ar[d]
	&	\W\B(A \ox A)	\ar[r]\ar[d]
	&	\W\B(A \ox A)	\ar[d]^{{\W h^\dagger}}\ar[r]
	&\W\B A					\ar[d]
	\\
	\W\B (X \ox X)	\ar[r]
	&	\W\B X					\ar[r]
	& P\,\mn\W\B X 				\ar[r]
	& \W\B X	
}
\]
(and symmetrically on the $A$-$Y$ side).
Tensoring this diagram with $\W\B X \from \W\B A \to \W\B Y$
yields another diagram inducing the map $\Tor_{\id \ox \W\Phi}$
at the top of [VII],
connected to the diagram 
inducing the same map $\Tor_{\id \ox \W\Phi}$ at the bottom of [IV]
by the triple $(\id_{\W\B X} \ox Z,\id_{\W\B A\ot},\id_{\W\B Y} \ox Z)$
in the third column, 
and by the identity elsewhere.

\benum
\item[VII:]
This follows from the functoriality of Tor
on applying $\g$, which is a natural transformation
$\g\: \W({-} \ox {-}) \lt \W({-}) \ox \W({-})$,
to the diagram   
\quation{\label{eq:VII-IV}
\xymatrix{
	\W\big(\B A \ox \B(A\ot)\big)	\ar[r]\ar[d]
&\W\big(\B A \ox \B(A\ot)\big)		\ar[r]\ar[d]
&\W\big(\B A \ox \B(A\ot)\big)		\ar[r]\ar[d]^{\W(\xi \ox h^\dagger)}
&\W(\B A \ox \B A)						\ar[d]
\\
	\W\big(\B X \ox \B (X\ot)\big)		\ar[r]
&	\W\big(\B X \ox \B X)					\ar[r]
& \W\big(\B X \ox \B PX)				\ar[r]
& \W\big(\B X \ox \B X)\mathrlap,
}
}
once we observe that
$\g \o \W(\xi \ox h^\dagger) = (\W\xi \ox \W h^\dagger)\o \g$.

\item[IX:]
This follows from functoriality of Tor 
on applying $\W\nabla$, 
which is a natural transformation
$\W\big(\B({-}) \ox \B({-})\big) \lt \W\B({-} \ox {-})$
to the diagram \eqref{eq:VII-IV}.
The right homotopy is now witnessed by
$\W(\xi \T h^\dagger)\: \W\B(A \ox A\ot) \lt \W\B(X \ox PX)$.
\eenum

\bs

\nd\emph{Another replacement.}
We now want to free the $P$ trapped inside the $\W\B$.
To this end, we note that there is a quasi-isomorphic embedding
$\iota\: X \ox PX \lt P(X \ox X)$
probably most easily understood by identifying each
as a subalgebra of $X \ox \I \ox X \iso \I \ox X \ox X$.
Postcomposing this $\iota$, we may replace 
$\Tor_{\W\B(A\ot)}\mn\big(\W\B(X \ox PX)\big)$ with
$\Tor_{\W\B(A\ot)}\mn\big(\W\B P(X\ox X)\big)$ on the bottom of [IX].
We may now further postcompose $Z$
and replace $\Tor_{\W\B(A\ot)}\mn\big(\W\B P(X\ox X)\big)$
with $\Tor_{\W\B(A\ot)}\mn\big(P\,\mn\W\B (X\ox X)\big)$.

\bs

\benum
\item[X:]
Consider the cube in Figure \ref{fig:homotopy-associativity-cube},
giving six maps
$\W\B(A\ott) \lt \W\B X$ 
and six homotopies between them,
each the composite of a 
map and a homotopy across one of the faces.

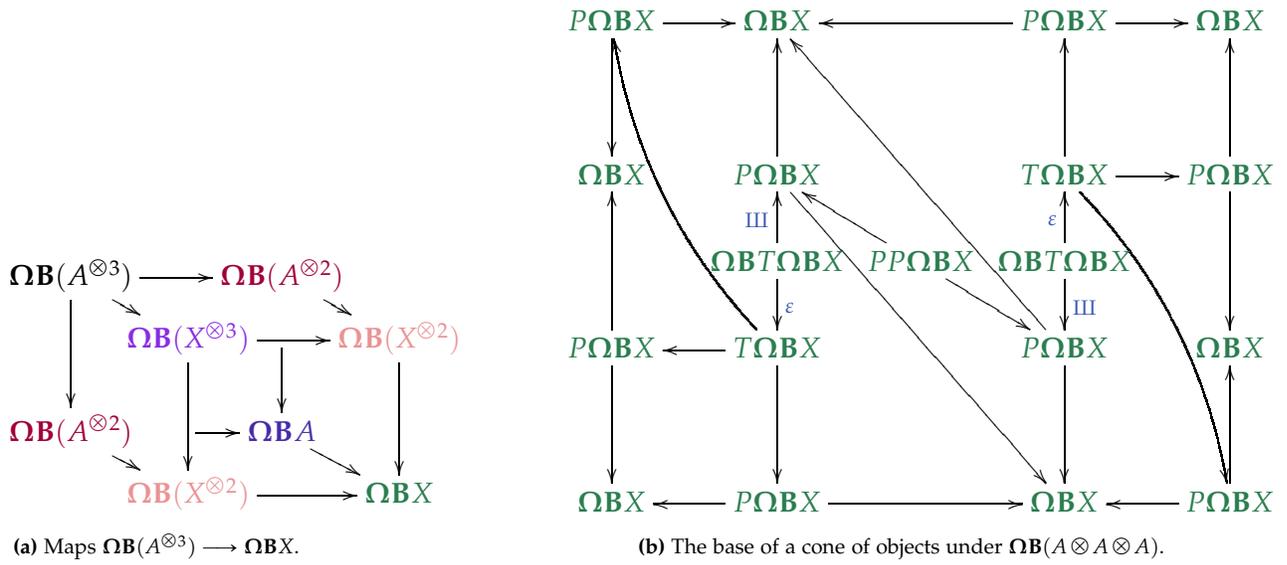
\begin{figure}
	\begin{subfigure}{0.3\textwidth}
		\centering
$
\begin{aligned}	
	\xymatrix@R=.5em@C=-.75em{	
			\W\B(A\ott)
		\ar@[jred][rrrrrr]
		\ar@[jred][dr]
		\ar@[jred][ddd]
		&&&&&&	
			\textcolor{jQ2}{\W\B(A\ot)}	
		\ar@[jred][dr]
		\ar@[jred][ddd]|!{[dlllll];[dr]}\hole
		&\\
		&
			\textcolor{jX2}{\W\B(X\ott)}
		\ar@[jred][rrrrrr]
		\ar@[jred][ddd]
		&&&&&	&
			\textcolor{jX2p}{\W\B (X\ot)}
		\ar@[jred][ddd]
		\\&&&&&&	
		\\
			\textcolor{jQ2}{\W\B(A\ot)}
		\ar@[jred][rrrrrr]|!{[uur];[dr]}\hole
		\ar@[jred][dr]
		&&&&&&	
			\textcolor{jA}{\W\B A}
		\ar@[jred][dr]
		&\\
		&
			\textcolor{jX2p}{\W\B(X\ot)}
		\ar@[jred][rrrrrr]
		&&&&&&			
			\textcolor{jgreen}{\W\B X}			
	}
\end{aligned}
$
		\caption{Maps $\W\B(A\ott) \lt \W\B X$.}
		\label{fig:homotopy-associativity-cube}
	\end{subfigure}%
	\hfill%
	\begin{subfigure}{0.6\textwidth}
		\centering
$
\begin{aligned}
	\xymatrix@C=.25em@R=1.75em{
		\ar@[jgreen][dd]
		\ar@[jgreen][rr]	
\textcolor{jgreen}{		P\,\mn\W\B X}					&&
\textcolor{jgreen}{		\W\B X}					&&
		\ar@[jgreen][ll]
		\ar@[jgreen][rr]	
\textcolor{jgreen}{		P\,\mn\W\B X}					&&
\textcolor{jgreen}{		\W\B X}									\\
		\\
\textcolor{jgreen}{		\W\B X}					&&
		\ar@[jgreen][uu]
		\ar@[jgreen][ddddrr]
\textcolor{jgreen}{		P\,\mn\W\B X}						&&
		\ar@[jgreen][rr]	
		\ar@[jgreen][uu]
		\ar@[jgreen]@/^1pc/[rrdddd]	
\textcolor{jgreen}{		T\W\B X}			&&
		\ar@[jgreen][dd]
		\ar@[jgreen][uu]	
\textcolor{jgreen}{		P\,\mn\W\B X}									\\
		&\ &		
		\ar@[jgreen][u]^{\textcolor{jblue}{	\Sha}}
		\ar@[jgreen][d]^{\textcolor{jblue}{\e}}
\textcolor{jgreen}{		\W\B T\W\B X}							
		&
		\ar@[jgreen][dr]
		\ar@[jgreen][ul]	
\textcolor{jgreen}{		P \mnn P\,\mn\W\B X}							
		&
		\ar@[jgreen][d]^{\textcolor{jblue}{	\Sha}}
		\ar@[jgreen][u]^{\textcolor{jblue}{\e}}
\textcolor{jgreen}{		\W\B T\W\B X}									
		&\ &
		\\
		\ar@[jgreen][dd]
		\ar@[jgreen][uu]	
\textcolor{jgreen}{		P\,\mn\W\B X}					&&
		\ar@[jgreen][dd]
		\ar@[jgreen][ll]	
		\ar@[jgreen]@/^1pc/[lluuuu]	
\textcolor{jgreen}{		T\W\B X}		&&
		\ar@[jgreen][dd]
		\ar@[jgreen][uuuull]
\textcolor{jgreen}{		P\,\mn\W\B X}					&&
\textcolor{jgreen}{		\W\B X}									\\
		\\
\textcolor{jgreen}{		\W\B X}					&&
		\ar@[jgreen][ll]
		\ar@[jgreen][rr]	
\textcolor{jgreen}{		P\,\mn\W\B X}					&&
\textcolor{jgreen}{		\W\B X}					&&
		\ar@[jgreen][uu]
		\ar@[jgreen][ll]	
\textcolor{jgreen}{		P\,\mn\W\B X}					
	}
\end{aligned}
$
		\caption{The base of a cone of objects under $\W\B(A \ox A \ox A)$.}
		\label{fig:assoc-path-plug}
	\end{subfigure}
	\caption{Auxiliary diagrams for the associativity argument.}
	\label{fig:assoc-extra-figs}
\end{figure}

\begin{figure}
	\centering{
		\xymatrix@C=-0.5em@R=1em{
			\ar@{=}[rrrr]
			\ar@{=}[dddd]\ar@/^2.5pc/@{=}[rrdd]|(.415)\hole|(.47)\hole |(.8)\hole
			\ar@[jred][rd]
			\W\B(A\ott)								&&
			&&
			\ar@{=}[rr]
			\ar@{=}[dd]|\hole
			\ar@{-->}@[jcmp][rd]
			\W\B(A\ott)								&&
			\ar@{=}[dd]|\hole
			\ar@[jred][rr]
			\ar@{~>}@[jhmt][rd]
			\W\B(A\ott)								&&
			\ar@[jred][rd]
			\textcolor{jQ2}{\W\B(A\ot)} 	
			\\&															
			\ar@[jred][rrrr]
			\ar@[jred][dddd]
			\ar@/_2.5pc/@{~>}@[jhmt][rrdd]
			\textcolor{jX2}{\W\B(X\ott)}	&&
			&&
			\ar@[jred][dd]
			\textcolor{jX2p}{\W\B(X\ot)}	&&	
			\ar@[jX2p][rr]
			\ar@[jX2p][ll]
			\ar@[jred][dd]
			\textcolor{jX2p}{P\,\mn\W\B(X\ot)}	&&	
			\ar@[jred][dd]
			\textcolor{jX2p}{\W\B(X\ot)}	&&									
			\\																
			&&
			\ar@{=}[rr]
			\ar@{=}[dd]|(.6)\hole
			\ar@{-->}@[jcmp][rd]
			\ar@/^2.5pc/@{~>}@[jhmt][rrrrrrdddddd]|(.275)\hole
			\W\B(A\ott)								&&
			\ar@{-->}@[jcmp][rd]
			\ar@{=}[rr]|\hole
			\W\B(A\ott)								&&
			\ar@[jred][rr]|\hole
			\ar@{-->}@[jcmp][rd]
			\W\B(A\ott)	 							&&
			\ar@[jQ2]@{=}[uu]|\hole
			\ar@[jQ2]@{=}[dd]|\hole
			\ar@{-->}@[jcmp][rd]
			\textcolor{jQ2}{\W\B(A\ot)} 	
			\\&																
			&&
			\ar@[jgreen][rr]
			\ar@[jgreen][dd]
			\textcolor{jgreen}{P\,\mn\W\B X} 		&&
			\textcolor{jgreen}{\W\B X} 			&&
			\ar@[jgreen][rr]
			\ar@[jgreen][ll]
			\textcolor{jgreen}{P\,\mn\W\B X} 		&&
			\textcolor{jgreen}{\W\B X} 												
			\\																	
			\ar@{=}[rr]|\hole
			\ar@{=}[dd]
			\ar@{-->}@[jcmp][rd]
			\W\B(A\ott)									&&
			\ar@{=}[dd]|\hole
			\ar@{-->}@[jcmp][rd]
			\W\B(A\ott)									&&
			&&
			&&
			\ar@[jred][dd]
			\ar@{~>}@[jhmt][rd]
			\textcolor{jQ2}{\W\B(A\ot)} 	
			\\&																	
			\ar@[jred][rr]
			\textcolor{jX2p}{\W\B(X\ot)}		&&
			\textcolor{jgreen}{\W\B X} 			&&
			&&
			&&
			\ar@[jgreen][uu]
			\ar@[jgreen][dd]
			\textcolor{jgreen}{P\,\mn\W\B X} 	
			\\																	
			\ar@{=}[rr]|\hole
			\ar@[jred][dd]
			\ar@{~>}@[jhmt][rd]
			\W\B(A\ott)								&&
			\ar@{-->}@[jcmp][rd]
			\ar@[jred][dd]|\hole
			\W\B(A\ott)								&&
			&&
			&&
			\ar@[jred][rd]
			\textcolor{jA}{\W\B A} 	
			\\&																	
			\ar@[jX2p][uu]
			\ar@[jX2p][dd]
			\ar@[jred][rr]
			\textcolor{jX2p}{P\,\mn\W\B(X\ot)}		&&
			\ar@[jgreen][uu]
			\ar@[jgreen][dd]
			\textcolor{jgreen}{P\,\mn\W\B X} 		&&
			&&
			&&
			\textcolor{jgreen}{\W\B X} 	
			\\																	
			\ar@[jred][rd]
			\textcolor{jQ2}{\W\B(A\ot)} 	&&
			\ar@[jQ2]@{=}[ll]|\hole
			\ar@[jQ2]@{=}[rr]|\hole
			\ar@{-->}@[jcmp][rd]
			\textcolor{jQ2}{\W\B(A\ot)} 	&&
			\ar@{~>}@[jhmt][rd]
			\ar@[jred][rr]
			\textcolor{jQ2}{\W\B(A\ot)} 	&&
			\ar@[jred][rd]
			\textcolor{jA}{\W\B A} 				&&
			\ar@[jA][ll]
			\ar@[jA][uu]
			\ar@[jred][rd]
			\textcolor{jA}{P\,\mn\W\B A} 	
			\\&																	
			\ar@[jred][rr]
			\textcolor{jX2p}{\W\B(X\ot)}	&&	
			\textcolor{jgreen}{\W\B X} 			&&
			\ar@[jgreen][rr]
			\ar@[jgreen][ll]
			\textcolor{jgreen}{P\,\mn\W\B X} 		&&
			\textcolor{jgreen}{\W\B X} 			&&
			\ar@[jgreen][uu]
			\ar@[jgreen][ll]
			\textcolor{jgreen}{P\,\mn\W\B X}
		}
	}
	\caption{A part of the system of \DGAs 
		underlying Square X from \Cref{fig:assoc}.}
	\label{fig:X-DGA-partial}
\end{figure}
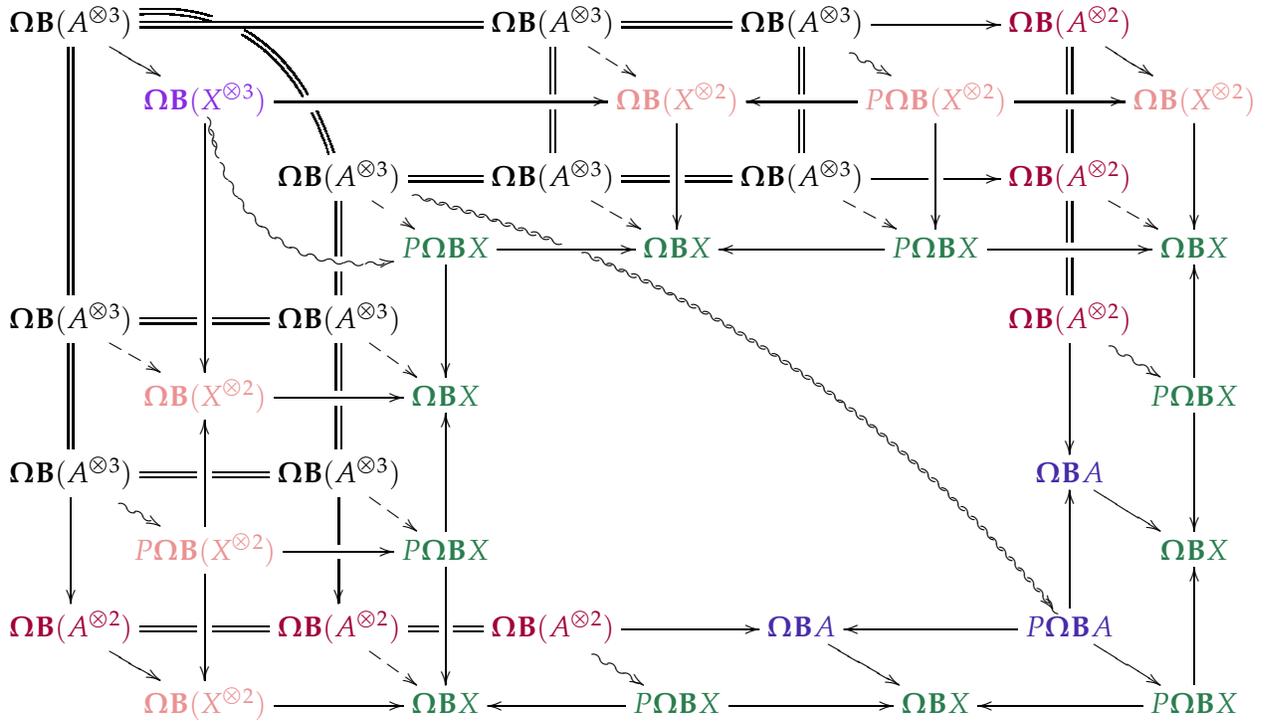

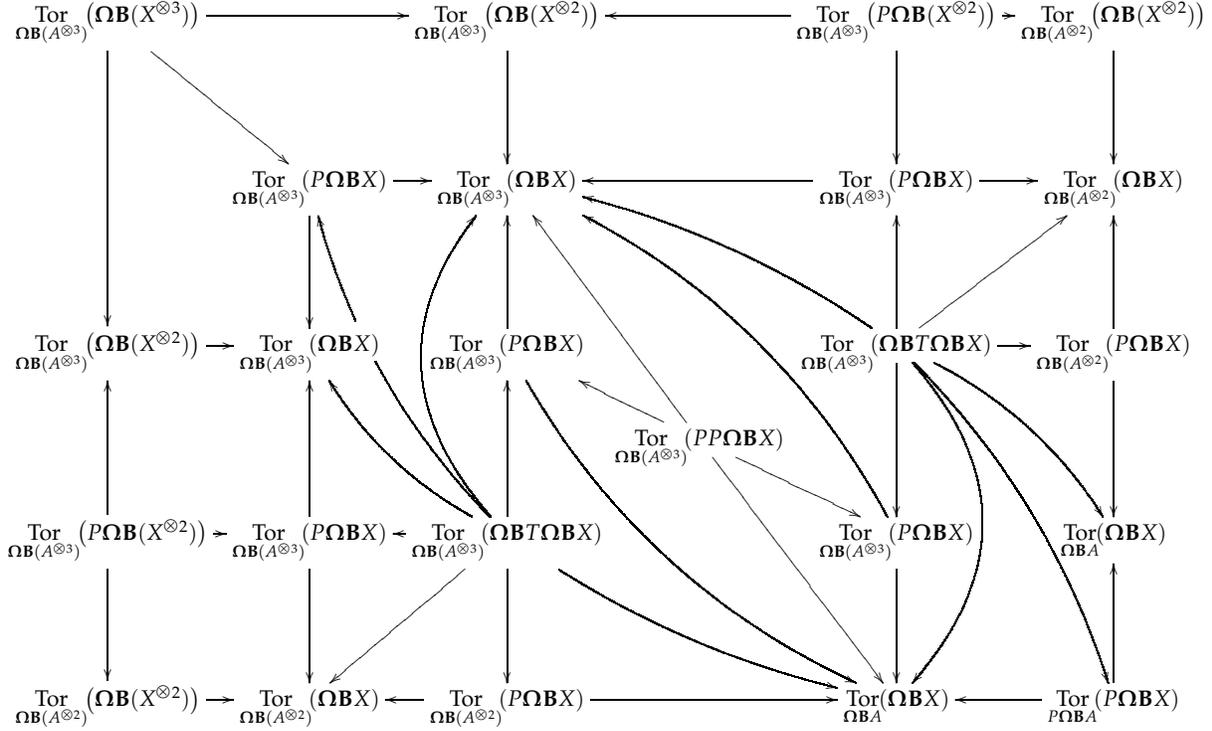
\begin{figure}
	\centering$
		\resizebox{149mm}{!}{
			\xymatrix@C=0em
			{
				\ar@[nonisored][rrdd]
				\ar@[nonisored][rrrr]
				\ar@[nonisored][dddd]
				\Tr{\W\B(A\ott)}\!\big(\W\B(X\ott)\big)
				&&&&
				\ar@[nonisored][dd]
				\Tr{\W\B(A\ott)}\!\big(\W\B(X\ot)\big)						&&
				\ar@[nonisored][rr]
				\ar[ll]
				\ar@[nonisored][dd]
				\Tr{\W\B(A\ott)}\!\big(P\,\mn\W\B(X\ot)\big)					&&
				\ar@[nonisored][dd]
				\Tr{\W\B(A\ot)}\!\big(\W\B(X\ot)\big)
				\\
				\\
				&&
				\ar[rr]
				\ar[dd]																		
				\Tr{\W\B(A\ott)}\!(P\,\mn\W\B X)									&&
				\Tr{\W\B(A\ott)}\!(\W\B X)									&&
				\ar@[nonisored][rr]
				\ar[ll]
				\Tr{\W\B(A\ott)}\!(P\,\mn\W\B X)									&&
				\Tr{\W\B(A\ot)}\!(\W\B X)	
				\\
				\\
				\ar@[nonisored][rr]						
				\Tr{\W\B(A\ott)}\!\big(\W\B(X\ot)\big)						&&
				\Tr{\W\B(A\ott)}\!(\W\B X)									&&
				\ar[uu]
				\ar@/_2pc/@[nonisored][ddddrr]
				\Tr{\W\B(A\ott)}\!(P\,\mn\W\B X)									&&
				\ar@/^1.125pc/@[nonisored][ddddrr]
				\ar@[nonisored][rr]
				\ar@/^1pc/@[nonisored][rrdd]
				\ar[uu]
				\ar[dd]
				\ar@/^3.5pc/@[nonisored][dddd]
				\ar@[nonisored][uurr]
				\ar@/_1pc/[uull]
		\smash{\Tr{{\W\B(A\ott)}}}\!
				(\W\B T\W\B X)									&&
				\ar@[nonisored][dd]
				\ar[uu]
				\Tr{\W\B(A\ot)}\!(P\,\mn\W\B X)									&&
				\\
				&&
				&&
				&
				\ar[dr]^{}
				\ar[ul]_{}
				\ar@[nonisored][dddr]
				\ar[uuul]
				\smash{\Tr{\W\B(A\ott)}}\!(P \mnn P\,\mn\W\B X)	
				&
				&
				&
				\\
				\ar@[nonisored][dd]
				\ar[uu]_{}		
				\ar@[nonisored][rr]																
				\Tr{\W\B(A\ott)}\!\big(P\,\mn\W\B (X\ot)\big)					&&
				\ar@[nonisored][dd]
				\ar[uu]_{}
				\Tr{\W\B(A\ott)}\!(P\,\mn\W\B X)									&&
				\ar@[nonisored][dd]
				\ar@/_1pc/@[nonisored][rrdd]
				\ar[ll]
				\ar@[nonisored][ddll]
				\ar@/^1pc/[uull]
				\ar@/^1.25pc/[uuuull]		|(.55)\hole	|(.625)\hole
				\ar@/^3.5pc/[uuuu]
				\ar[uu]
				\Tr{\W\B(A\ott)}\!(\W\B T\W\B X)									&&
				\ar@/_2pc/[uuuull]
				\ar@[nonisored][dd]
				\Tr{\W\B(A\ott)}\!(P\,\mn\W\B X)									&&
				\Tr{\W\B A}\!(\W\B X)									
				\\
				\\
				\ar@[nonisored][rr]_(.55){}																		
				\Tr{\W\B(A\ot)}\!\big(\W\B (X\ot)\big)					&&
				\Tr{\W\B(A\ot)}\!(\W\B X)										&&
				\ar[ll]^{}
				\ar[rr]
				\Tr{\W\B(A\ot)}\!(P\,\mn\W\B X)									&&
				\Tr{\W\B A}\!(\W\B X)											&&
				\ar[ll]^{}
				\ar[uu]_{}
				\Tr{P\,\mn\W\B A}\!(P\,\mn\W\B X)										
			}
		}
	$
	\caption{Square X from Figure \ref{fig:assoc}, filled.}
	\label{fig:X-filled}
\end{figure}

\hspace{20pt}We may use the associated right homotopies to fill out
Figure \ref{fig:X-DGA-partial}.
The edges from \eqref{fig:homotopy-associativity-cube} are red, 
the right homotopies corresponding to the faces are gold and wavy,
composite maps are grey and dashed,
and quasi-isomorphism classes of \DGAs are color-coded.
In particular all endpoint maps $\pi_0$ and $\pi_1$ are green.
The top of Figure \ref{fig:X-DGA-partial} comes from the replacement 
we have just made,
and the left from its suppressed, symmetric twin following the 
parallel square IX';
the right and bottom edges are both from
\Cref{thm:homotopy-Tor-map}
and induce $\Tor_{\W\Phi}$.
The diagram as it stands now commutes by definition,
and it remains to fill in the interior.

\hspace{20pt}By \Cref{thm:homotopy-compose}, 
the homotopies from $\W\B(A\ott)$
can be composed, and by \Cref{def:Sha}
the composite of two consecutive triples can 
be represented by a single right homotopy.
By \Cref{thm:homotopies-homotopic},
these composite right homotopies $\W\B(A\ott) \lt P\,\mn\W\B X$
are themselves homotopic,
and this is witnessed by a right homotopy $\W\B(A\ott) \lt P \mnn P\,\mn\W\B X$.
We can combine all the codomains into the 
Figure \ref{fig:assoc-path-plug}, 
to be thought of as the base of a cone under $\W\B(A\ott)$.

\new{
As with \Cref{fig:comm-plug} and \eqref{eq:diamonds}, 
this cone is not commutative, though by assumptions the faces apart from the base are.
Again, the issue can be repaired with the fix of
\Cref{thm:homotopy-independence-mk2},
at the cost of assuming
the right homotopy 
$H\: \W\B(A\ott) \lt PP\W\B X$
in the middle of the cone be {endpoint-fixing}.
That means, again,
that the two maps
$P\pi_j \o H\:\W\B(A\ott) \lt P\W\B X$
must factor through
the natural map $\z\: \W\B X \lt P \W\B X$ of \Cref{def:P}.
The relevant cone from $\W\B(A\ott)$ is thus in fact that over the base 
in \Cref{fig:true-base}.
}

\new{
As before, there does not seem to be 
an obvious criterion determining when this is actually possible.
}

\begin{figure}
	\centering
$
\xymatrix@C=-.25em@R=2.5em{
	\textcolor{jgreen}{\W\B X} 
	\ar@[jgreen][drrrr]_{\id}
	\ar@[jgreen][ddrrr]|{\z}
	\ar@[jgreen][dddrr]^{\id}
	&&
	\textcolor{jgreen}{P\,\mn\W\B X}
	\ar@[jgreen][ll]_{\pi_0}
	\ar@[jgreen][rr]^{\pi_1}
	&&
	\textcolor{jgreen}{\W\B X}&&
	\textcolor{jgreen}{P\,\mn\W\B X} 
	\ar@[jgreen][ll]_{\pi_0}
	\ar@[jgreen][rr]^{\pi_1}
	&&
	\textcolor{jgreen}{\W\B X}
	\\
	&&&&
	\textcolor{jgreen}{\W\B X}
	&&
	\textcolor{jgreen}{\W\B T\W\B X} 
	\ar@[jgreen][llllu]
	\ar@[jgreen][u]
	\ar@[jgreen][dl]_(.45){\textcolor{jblue}{\Sha}}
	\ar@[jgreen][ddrr]
	\\
	&&&
	\textcolor{jgreen}{P\W\B X}
	\ar@[jgreen][dl]^{\pi_0}
	\ar@[jgreen][ur]_{\pi_1}
	&&
	\textcolor{jgreen}{P\W\B X} 
	\ar@[jgreen][ul]^{\pi_0}
	\ar@[jgreen][dr]_{\pi_1}
	\\
	\textcolor{jgreen}{P\,\mn\W\B X}  
	\ar@[jgreen][ddd]_{\pi_1}
	\ar@[jgreen][uuu]^{\pi_0}
	&&
	\textcolor{jgreen}{\mn\W\B X}
	&& 
	\textcolor{jgreen}{P \mnn P\,\mn\W\B X}
	\ar@[jgreen][dl]|{\pi_0}
	\ar@[jgreen][ur]|{\pi_1}
	\ar@[jgreen][ul]|{P\pi_0}
	\ar@[jgreen][dr]|{P\pi_1}
	&&
	\textcolor{jgreen}{\W\B X} 
	&&
	\textcolor{jgreen}{P\,\mn\W\B X}
	\ar@[jgreen][ddd]^{\pi_1}
	\ar@[jgreen][uuu]_{\pi_0}
	\\
	&&&
	\textcolor{jgreen}{P\,\mn\W\B X} 
	\ar@[jgreen][ul]_{\pi_0}
	\ar@[jgreen][dr]^{\pi_1}
	&&
	\textcolor{jgreen}{P\,\mn\W\B X}
	\ar@[jgreen][dl]_{\pi_0}
	\ar@[jgreen][ur]^{\pi_1}
	\\
	&&
	\textcolor{jgreen}{\W\B T\W\B X} 
	\ar@[jgreen][ur]_(.45){\textcolor{jblue}{\Sha}}
	\ar@[jgreen][d]
	\ar@[jgreen][uull]
	\ar@[jgreen][drrrr]
	&&
	\textcolor{jgreen}{\W\B X} 
	\\
	\textcolor{jgreen}{\W\B X} 
	&&
	\textcolor{jgreen}{P\,\mn\W\B X} 
	\ar@[jgreen][ll]^{\pi_0}
	\ar@[jgreen][rr]_{\pi_1}
	&&
	\textcolor{jgreen}{\W\B X} 
	&& 
	\textcolor{jgreen}{P\,\mn\W\B X} 
	\ar@[jgreen][ll]^{\pi_0}
	\ar@[jgreen][rr]_{\pi_1}
	&&
	\textcolor{jgreen}{\W\B X}
	\ar@[jgreen][ullll]_{\id}
	\ar@[jgreen][uulll]|{\z}
	\ar@[jgreen][uuull]^{\id}
}
$
\caption{\new{The true base.}}
	\label{fig:true-base}
\end{figure}

Using the factorizations of the maps along the right and bottom edges
through $\W\B(A\ot)$ and $\W\B A$,
we may insert this cone into Figure \ref{fig:X-DGA-partial},
and taking Tor, obtain Figure \ref{fig:X-filled},
in which the black arrows are isomorphisms and the red are not.
This is square X of Figure \ref{fig:assoc},
and using the commutativity of all its constituent 
squares and triangles,
we see it commutes.\qedhere
\eenum
\end{proof}



\section{Functoriality of the product}\label{sec:ring-map}

Now that we have a ring structure on Tor, 
we would like also to have ring maps.
It is clear the product \eqref{eq:full-product-v2}
is functorial with respect to spans of \DGA maps, 
but our claims are more expansive.

\begin{theorem}\label{thm:main-algebraic}
Given {\WHCA}s
$A'$, $X'$, $Y'$,
$A$, $X$, $Y$ 
and \WHCA maps 
\[
	\xymatrix@C=2.8em@R=3em{
	{\vphantom{()}}
	\B X' 	\ar[d]_{{\l_{\mn X}}}		& 
	\B A'	\ar[r]^{\upsilon'}
	\ar[l]_{\xi'}
	\ar[d]|(.45)\hole|{ {\l_A}}|(.575)\hole& 
	\B Y' 	\ar[d]^{{\l_Y}{\vphantom{\T}}}	
	\ar@{}[l]_{{\vphantom{\Phi_{A'}}}}			
	\\
	{\vphantom{()}}	
	\B X												& 
	\B A	\ar[l]^{\xi}
	\ar[r]_{\upsilon}							& 
	\B Y 
	\ar@{}[r]_{{\vphantom{\Phi_{A}}}}			&
}
\]
 such that the squares commute up to \DGC homotopy,
 the $\kk$-linear map 
\[
\defm \Miso \ceq \Tor_{\W\l_A}(\W\l_X,\W\l_X)\:
\Tor_{\W\B A'}(\W\B X',\W\B Y') \lt \Tor_{\W\B A}(\W\B X,\W\B Y)
\]
defined as in \Cref{thm:homotopy-Tor-map}
is multiplicative with respect to the products
\eqn{
	\defm{\Pi'}\: &\Tor_{\W\B A'}(\W\B X',\W\B Y')\ot \lt \Tor_{\W\B A'}(\W\B X',\W\B Y')\mathrlap,\\
	\defm{\Pi}\: &\Tor_{\W\B A}(\W\B X,\W\B Y)\ot \lt \Tor_{\W\B A}(\W\B X,\W\B Y)
}
described in \eqref{eq:full-product-v2}.
That is,
$\Pi \o (\Miso \ox \Miso) = \Miso \o \Pi'$.

These algebra homomorphisms are functorial in the following sense:
given further \WHCA maps $\BX'' \from \BA'' \to \BY''$
and a triple $(\l'_X,\l'_A,\l'_Y)$ of \WHCA-maps
to $\BX' \from \BA''\to \BY'$
such that the resulting squares commute up to \DGC homotopy,
and hence an algebra map 
\[
\Xi' = \Tor_{\W\l'_A}(\W\l'_X,\W\l'_Y)\:\Tor_{\W\BA''}(\W\BX'',\W\BY'') \lt \Tor_{\W\BA'}(\W\BX',\W\BY')\mathrlap,\]
the composite $\Xi \o \Xi'$
is equal to 
$\Tor_{\W\l_A \o \,\W\l'_A}(\W\l_X \o \W\l'_X,\W\l_Y \o \W\l'_Y)$\new{,
subject to a usually intractible and unverifiable additional condition 
on the various defining homotopies,
to be gestured at briefly toward the end of the proof below}.
\end{theorem}

\begin{corollary}\label{thm:corollary-main-algebraic}
	In the situation of \Cref{thm:main-algebraic},
	suppose that 
	$A' = \H A$, 
	$X' = \H X$,   
	$Y' = \H Y$,
	with 
	$\xi_*\: \H A \lt \H X$
	and $\upsilon_*\: \H A \lt \H Y$
	the \DGA maps obtained
	by conjugating $\H\W\xi$
	and $\H\W\upsilon$ respectively by $\H\e$.
	Then the induced map 
$
\Tor_{\W\B\H A}(\W\B\H X,\W\B\H Y) 
\isoto \Tor_{A}(\W \B X,\W \B Y)$
is multiplicative with respect to the products induced as in \Cref{sec:product}
by the homotopies witnessing that 
$\xi$, $\upsilon$, $\B\xi_*$, and $\B\upsilon_*$ are \SHCA maps.
\end{corollary}

\wrong{
\begin{proof}[Proof of {\Cref{thm:main-topological}}.]
Conjugating by $\Tor_\e$, 
we may replace $\Tor_{\C(B)}$ with $\Tor_{\W\B\C(B)}$
and $\Tor_{\H(B)}$ with $\Tor_{\W\B\H(B)}$.
Now,
	assuming the spaces in the span
	$X \from B \to E$ have polynomial cohomology,
	Munkholm~\cite[7.2]{munkholm1974emss}
	uses the \SHCA structure on 
	singular cochains from 
	\Cref{thm:SHC-cochain}
	to produce a trio of \DGC quasi-isomorphisms
\[
\xymatrix@C=1em@R=3em{
	\B\H(X)\ar[d]_{\defm{\l_X}}		&
	\B\H(B)\ar[l]\ar[r]\ar[d]|{\defm{\l_B}}		&	
	\B\H(E)\ar[d]^{\defm{\l_E}}		\\
	\B\C(X)				&
	\B\C(B)\ar[l]\ar[r]				&	
	\B\C(E)	
	\mathrlap.
}
\]
If $\kk$ has characteristic $2$
then for each $Z \in \{X,B,E\}$,
if the $\cupone$-squares of a set of polynomial generators of $\H(Z)$
vanish,
$\l_Z$ is a \WHCA map;
if $\kk$ has characteristic $\neq 2$,
they are \WHCA maps no matter what~\cite[7.3]{munkholm1974emss}.
If $\l_X$ and $\l_E$ are \WHCA maps,
the squares commute up to \DGC homotopy~\cite[7.4]{munkholm1974emss},
and since we assume
additionally $\l_B$ is a \WHCA map, 
\Cref{thm:corollary-main-algebraic} applies
to show the isomorphism is multiplicative.
\end{proof}
}

%

\begin{remark}\label{rmk:exterior}
		Munkholm's
		additive predecessor of \Cref{thm:main-algebraic}
		requires $\l_A$ be a \DGC map only;
		in his intended case,
		this is the only map doing the job anyway,
		but his proof does not require it to be an \SHCA map.
		Our proof, on the other hand, \emph{does} require $\l_A$ be a \WHCA map.
		
		In the application to \Cref{thm:main-topological}, 
		when $A' = \H(B)$ is a polynomial algebra, 
		this only poses an additional
		restriction if the characteristic of $\kk$ is $2$,
		but the obstruction is genuine and not merely a defect of the proof.
		If $Z$ is topological space 
		with polynomial $\H(Z;\F_2)$ 
		but such that $\cupone$-squares of generators $z$ 
		are not all decomposable (\emph{i.e.}, 
		if it cannot be guaranteed that the $z \cupone z$ lie in 
		$\wt H^*(Z;\F_2) \. \wt H^*(Z;\F_2)$),
		then Saneblidze showed
		$\H(\Omega Z;\F_2)$ is not exterior~\cite[Cor.~1]{saneblidze2017loop}.	
		In such cases,
		$A = \C(Z;\F_2)$ and $X = Y = \C({*};\F_2) = \F_2$ 
		is a counterexample to the desired strengthening 
		of \Cref{thm:main-algebraic}.
		
		Other work of Munkholm also 
		analyzes this situation,
		and \emph{exotic} \SHCA structures
		on cochains.
\end{remark}

\brmk\label{rmk:Panov}
The result is more sensitive than it might appear,
and the author still does not know cases where 
the strategy of \Cref{thm:main-topological}
yields an isomorphism if
the cohomology of the inputs is not polynomial,
even when \DGA quasi-isomorphisms between cochains and cohomology
are known.
In what is arguably the next-best case,
that of Davis--Januszkiewicz spaces, the strategy already fails.
The failure of two separate published putative proofs of analogous results 
in this case was a major motivator for our decision
to conduct the proof as much as possible at the level of diagrams of \DGA maps.
\ermk


\subsection*{Proof}

The rest of the section constitutes the proof of \Cref{thm:main-algebraic}.
The functorial nature of the maps follows immediately from 
\Cref{thm:functoriality}, so it remains to show multiplicativity.
That is, we are to connect the products $\Pi'$ and $\Pi$ described in \Cref{def:product}
using the map $\Miso$ of \Cref{thm:main-algebraic}:
\[
\xymatrix{
	\Tor_{\W\B A'}(\W\B X',\W\B Y')\ot \ar[d]_{\Miso \ox \Miso}\ar[r]^{\Pi'} &
	\Tor_{\W\B A'}(\W\B X',\W\B Y') \ar[d]^\Miso\\
	\Tor_{\W\B A}(\W\B X,\W\B Y)\ot \ar[r]_\Pi &
	\Tor_{\W\B A}(\W\B X,\W\B Y)\mathrlap.
}
\]
Expanding out definition of the product partially
and employing the space-saving convention of 
\Cref{def:suppression},
we will fill in the following 
diagram
in such a way that commutativity of each square is manifest:
\begin{equation}\label{eq:grand-scheme}
\begin{aligned}
\xymatrix@C=1em@R=1em{
	\us{\W\B A'}{\Tor}{}\ot				\ar[rr]\ar[dd]_{\Miso \ox \Miso}&&
	\us{(\W\B A')\ot}\Tor				\ar[dd]&&	
	\us{\W(\B A')\ot}\Tor				\ar[ll]_\sim
									\ar[rr]^(.375)\sim 
									\ar[dd]&&
	\smash{\us{{\W(\BA')\ot}}\Tor}\!\big(\W\B(X')\ot\big)			\ar[rr] \ar[dd]&&
	\us{\W\B A'}\Tor					\ar[dd]^\Miso				
														\\
	& {\smash{\mathclap{\mathrm{external}}}} {\phantom{}}&&\mathclap{\g}&&
	\mathclap{\W\nabla}&& \,\, \mathllap{\Phi}\ 
														\\	
	\us{\W\B A}{\Tor}{}\ot					\ar[rr]&& 
	\us{(\W\B A)\ot}\Tor						&& 
	\us{\W(\B A)\ot}\Tor				\ar[ll]_\sim 
									\ar[rr]^(.375)\sim&&
	\us{\W(\BA)\ot}\Tor\!\big(\W\B(X\ot)\big)					\ar[rr]&&
	\us{\W\B A}\Tor	\mathrlap.							
}
\end{aligned}
\end{equation}
Expanding all squares simultaneously 
would result in a diagram
more intimidating than illuminating,
so we will consider each square separately in its own subsection.

Before we do that, we should give the reader some explanation as to why 
we might expect the thing to commute at all.
The external product, $\gamma$, and $\W\nabla$ are all natural transformations,
so one should expect the squares involving them to commute,
and they do,
transforming the objectwise tensor-square of the three-square diagram of
\eqref{eq:Tor-DGA-homotopy-squares} determining $\Miso$
into another such three-square diagram determining the left edge 
of the $\Phi$ square. 
The only casualty in this process is the right homotopy,
which is transmogrified from a standard right homotopy
to one encumbered with an increasingly ornate
witnessing path object. 
We then have to deploy the material developed in \Cref{sec:homotopy}
to recover a right homotopy witnessed by the standard path object.
The edges of the $\Phi$ square come from
three- or four-square diagrams per \eqref{eq:Tor-DGA-homotopy-squares}
determined by the new homotopy we have transported over from
the one giving $\Miso \ox \Miso$
and three of the homotopies appear as hypotheses for \Cref{thm:main-algebraic}.
Filling in the $\Phi$ square amounts to constructing \DGA maps making
these homotopies coherent, 
and for this we will use
the techniques developed in \Cref{sec:homotopy}.

In what follows, we will continuously use 
the functoriality of $\Tor$, 
viewed as a graded $\kk$-module, in triples of $\kk$-\DGA maps making
the two squares \eqref{eq:Tor-DGA-functoriality-squares} commute. 
All homotopies will be expanded in terms
of path objects so that nothing is swept under the rug.
Because the squares involving $A$ and $Y$ are notationally symmetric with
those involving $A$ and $X$, 
we truncate to the $A$-$X$ portions
what will nevertheless be a crushing overburden of diagrams.
Every argument will proceed symmetrically and silently on the $A$-$Y$ side.

\begin{notation}\label{thm:proof-notation}
Recall from the statement of \Cref{thm:main-algebraic}
that we assume a \DGC homotopy between the two paths around the 
square
\[
\xymatrix@C=1em{
\BA' \ar[r]^{\xi'}\ar[d]_{\l_A}			& \BX'\ar[d]^{\l_X}\\
\BA \ar[r]_\xi								& \BX\mathrlap.
}
\]
Via the adjunction of \Cref{thm:homotopy-adjunction}, 
we obtain a \DGA homotopy $\defm h = h_X$ between 
the two transposed maps $\W\B A' \lt X$.
Write $\defm{h^P}\: \W\B A' \lt PX$
for the \DGA map 
representing this \DGA homotopy per \Cref{def:P}.
	We write $(h^P)^\#\: \W\B A' \lt \W\B PX$
	for the map induced up as in \Cref{thm:induction},
	which encodes a homotopy between 
	{$\W(\xi \o \l_A)$ and $\W(\l_{\mn X} \o \xi')\: \W\B A' \lt \W\B X$}
	in the sense that one gets these maps
	back from $(h^P)^\#$ by postcomposing
	$\W\B\pi_0$ and $\W\B\pi_1$ respectively.
\end{notation}

\subsection{The external product square}\label{sec:external}

To express $\Miso\ox \Miso$ on the left 
of \eqref{eq:grand-scheme}
in terms of \DGA maps,
we will find it more convenient to describe 
$\Miso$ using $(h^P)^{\#}$
rather than with the {standard} right homotopy $\W\B A' \lt P\,\mn\W\B X$
{with the same endpoints}.
Fortunately, \Cref{thm:P-W} gives us 
a quasi-isomorphism $Z\: \W\B P X \lt P\,\mn\W\B X$
such that $\pi_{\mn j} \o Z = \W\B\pi_{\mn j}$,
so we can replace the expected diagram on the left below with that
on the right and have the same induced map $\Miso$ on Tor.

\quation{\label{eq:Tor-DGA-homotopy-square-substitution}
	\begin{aligned}
		\xymatrix@C=1.90em{
			\W\B A'
			\ar@{=}[d]
			\ar[r]^{\W\xi'}
			& 	\W\B X\mathrlap'
			\ar[d]^{\W\l_{\mn X}}\\
			\W\B A'\ar@{=}[d]
			\ar[r]	
			& \W\B X\\
			\W\B A' \ar[d]_{\W\l_A}
			\ar[r]	
			& P\,\mn\W\B X\ar[u]_{\pi_0}\ar[d]^{\pi_1}\\
			\W\B A \ar[r]_{\xi}
			& \W\B X 
		}
		\qquad
		\quad
		\qquad
		\xymatrix@C=1.90em{
			\W\B A'
			\ar@{=}[d]
			\ar[r]^{\W\xi'}
			& 	\W\B X\mathrlap'
			\ar[d]^{\W\l_{\mn X}}\\
			\W\B A'\ar@{=}[d]
			\ar[r]	
			& \W\B X\\
			\W\B A' \ar[d]_{\W\l_A}
			\ar[r]	
			& \W\B PX\ar[u]_{\W\B\pi_0}\ar[d]^{\W\B\pi_1}\\
			\W\B A \ar[r]_{\xi}
			& \W\B X 
		}
\end{aligned}
}

Then the external product square is the composite of subsquares
\begin{equation}\label{eq:external-Tor-squares}
\begin{aligned}
\xymatrix{
	\Tor_{\W\B  A'}(\W\B X'\Ysub)\ot \ar[r]\ar[d] & 
		\Tor_{(\W\B  A')\ot}\big((\W\B X')\ot\Ysub\big)\ar[d]\\
	\Tor_{\W\B  A'}(\W\B X\Ysub)\ot\ar[r]&
		\Tor_{(\W\B  A')\ot}\big((\W\B X)\ot\Ysub\big)\\
	\Tor_{\W\B  A'}(\W\B PX\Ysub)\2\ar[r]\ar[d]\ar[u]^\vertsim&
		\Tor_{(\W\B  A')\2}\big((\W\B PX)\2\Ysub\big)\ar[d]\ar[u]_\vertsim	\\
	\Tor_{\W\B A}(\W\B X\Ysub)\ar[r]&
		\Tor_{(\W\B A)\ot}\big((\W\B X)\ot\Ysub\big)
	\mathrlap,
	}
\end{aligned}
\end{equation}
in which each horizontal map is the  exterior product
and the vertical maps in each file are determined functorially 
by applying Tor to the right diagram
of \eqref{eq:Tor-DGA-homotopy-square-substitution}
and its tensor-square.
Here, as promised, we have suppressed the symmetric $A$-$Y$ half 
of the diagram in a bid for comprehensibility,
and the unlabeled map is the necessary composition rendering the diagram commutative.

In summary, 
the input diagram of \DGAs commutes by definition
and the output diagram of Tors commutes 
by the functoriality of the exterior product.

\subsection{The \texorpdfstring{$\g$}{gamma} square}\label{sec:gamma}

What we call the $\g$ square in \eqref{eq:grand-scheme}
arises by applying the cobar shuffle
$\g\: \W({-} \ox {-}) \lt \W({-}) \ox \W({-})$ 
of \Cref{def:shuffle}
across the board, 
landing in the right edge of the Tor diagram 
\eqref{eq:external-Tor-squares}
from the external product square of the preceding subsection:

\[
	\xymatrix{
		\Tor_{(\W\B  A')\ot}\big((\W\B X')\ot\Ysub\big)
		\ar[d] & 
		\ar[l]^\sim_(.48){\Tor_\g}
		\Tor_{\W(\B  A')\ot}\big(\W(\B X')\ot\Ysub\big)
		\ar[d]\\
		\Tor_{(\W\B  A')\ot}\big((\W\B X)\ot\Ysub\big)
		&
		\ar[l]^\sim_(.48){\Tor_\g}
		\Tor_{\W(\B  A')\ot}\big(\W(\B X)\ot\Ysub\big)
		\\
		\Tor_{(\W\B  A')\2}\big((\W\B PX)\2\big)
		\ar[d]\ar[u]^\vertsim&
		\ar[l]_\sim^(.48){\Tor_\g}
		\Tor_{\W(\B  A')\2}\big(\W (\B PX)\2\big)
		\ar[d]\ar[u]_\vertsim	\\
		\Tor_{(\W\B A)\ot}\big((\W\B X)\ot\Ysub\big)
		&
		\ar[l]_\sim^(.48){\Tor_\g}
		\Tor_{\W(\B A)\ot}\big(\W(\B X)\ot\Ysub\big)
		\mathrlap.
	}
	\]	

\bs

\nd This Tor diagram is induced by the \DGA diagram

\bs

\[
\xymatrix@R=.75em@C=3.5em{
	(\W\B A')\ot\ar@{=}[dd] \ar[rd] 
	&&&\W(\B A')\ot\ar@{=}[dd]|\hole  
				\ar[rd]^(.55){\W(\xi')\ot}
\ar[lll]	\\
	&( \W\B X' )\ot \ar[dd]
	&&&\W(\B X' )\ot \ar[dd]^{\W(\l_{\mn X}\ot)} 
\ar[lll]	\\
	(\W\B A')\ot\ar@{=}[dd]\ar[rd]
	&&&\W(\B A')\ot\ar@{=}[dd]|\hole \ar[rd] 
\ar[lll]|!{[l];[ll]}\hole 	\\
	&(\W\B X )\ot
	&&&\W(\B X )\ot
\ar[lll]	\\
	(\W\B A')\ot \ar[dd] \ar[rd]_{(\mn(h^P)^\#\mn)\ot\ \ \,}
	&&&\W(\B A')\ot\ar[rd]^{\W(h^\dagger)\ot} 
						\ar[dd]|\hole _(.275){\W(\l_{\! A}\ot)\!} 
\ar[lll]|!{[l];[ll]}\hole 	\\
	&(\W\B PX)\ot \ar[uu]\ar[dd]
	&&&\W(\B PX)\ot 	\ar[uu]_{\W(\B\pi_0)\ot} 
					\ar[dd]^{\W(\B\pi_1)\ot} 
\ar[lll]	\\
	(\W\B A)\ot \ar[dr]
	&&&\W(\B A)\ot\ar[rd]^{\W(\xi\ot)} 
\ar[lll]|!{[l];[ll]}\hole \\
	&(\W\B X)\ot
	&&&\W(\B X)\ot\mathrlap,
\ar[lll]\\
}
\]

\bs

\nd in which all horizontal arrows are $\g$
and $\defm{h^\dagger}	\: \B A' \lt \B P X$ is the transpose of 
					$h^P\: \W\B A' \lt PX$.
The left face commutes by the previous step, 
and the top two and bottom horizontal faces commute by 
naturality of $\g$, 
but the right face
and the third horizontal face 
remain to be explained.

Since the transpose is given by applying $\B$ and precomposing 
$\h\: \B A' \lt \B\W\B A'$,
and $h^P$ represents the transpose
of the original homotopy $\xi \o \l_A \hmt \l_{\mn X} \ox \xi'$
it follows that $\B \pi_{\mn j} \o h^\dagger = \B(\pi_{\mn j} \o h^P)\o \h$ 
are respectively $\xi \o \l_A$ and $\l_{\mn X} \ox \xi'$ 
again for $j \in \{0,1\}$,
giving commutativity of the right face.
The remaining horizontal face commutes
since $\W h^\dagger = \W\B h^P \o \W\h = (h^P)^\#$
by \Cref{thm:induction}.

\subsection{The \texorpdfstring{$\W\nabla$}{shuffle} square}\label{sec:nabla}


In what we called the $\W\nabla$ square in \eqref{eq:grand-scheme},
all horizontal maps are $\id$ or $\W\mnn\nabla$. 
The expanded rectangle of Tors is

\begin{equation}\label{eq:nabla-Tor-squares}
\begin{aligned}
	\xymatrix@C=4em{
		\Tor_{\W(\B  A')\ot}  \big(\W(\B X')\ot\Ysub\big)
				\ar[d] 
				\ar[r]_{\substack{\phantom{x}}{\sim}}^{\Tor_{\id}(\W\mnn\nabla)}
		& 
		\Tor_{\W(\B A')\ot}  \big(\W\B(X')\ot\Ysub\big)
				\ar[d]\\
		\Tor_{\W(\B  A')\ot}  \big(\W(\B X)\ot\Ysub\big)
				\ar[r]_{\substack{\phantom{x}}{\sim}}^{\Tor_{\id}(\W\mnn\nabla)}
		&
		\Tor_{\W(\B A')\ot}  \big(\W\B(X\ot)\Ysub\big)
		\\
		\Tor_{\W(\B  A')\2}  \big(\W (\B PX)\2\big)
				\ar[r]^\sim_{\Tor_{\id}(\W\mnn\nabla)}
				\ar[d]
				\ar[u]^\vertsim&
		\Tor_{\W(\B A')\2}\big(\W\B(PX)\2\big)
				\ar[d]
				\ar[u]_\vertsim	
		\\
		\Tor_{\W(\B A)\ot}  \big(\W(\B X)\ot\Ysub\big)
				\ar[r]^\sim_{\Tor_{\id}(\W\mnn\nabla)}
		&
		\Tor_{\W\B(A\ot)}  \big(\W\B(X\ot)\Ysub\big)
		\mathrlap.
	}
\end{aligned}
\end{equation}

The prism of inducing \DGA maps is morally (but not exactly)
the following, in which all horizontal maps are $\W\mnn\nabla$:

\bs

\[
\xymatrix@R=.75em@C=3.5em{
	\W(\B A')\ot
\ar@{=}[dd] 
\ar[rd]_{\W(\xi')\ot\ \ }
\ar[rrr]
	&&&	\W\B (A')\ot
\ar@{=}[dd]|\hole  
\ar[rd]^{\W(\xi' \T \xi')}
	\\
	&\W(\B X' )\ot
\ar[dd]^(.275){\W(\l_{\mn X}\ot)} 
\ar[rrr]
	&&&\W\B(X' )\ot
\ar[dd]^{\W(\l_{\mn X} \T \l_{\mn X})}
	\\
	\W(\B A')\ot
\ar[rrr]|!{[ur];[dr]}\hole
\ar@{=}[dd]\ar[rd]  
	&&&\W\B (A')\ot
\ar@{=}[dd]|\hole \ar[rd] 
	\\
	&\W(\B X )\ot
\ar[rrr]
	&&&\W\B(X\ot)
	\\
	\W(\B A')\ot
\ar[dd]_{\W(\l_A\ot)\!} 
\ar[rd]_{\W(h^\dagger)\ot\ \ } 
\ar[rrr]|!{[ur];[dr]}\hole 
	&&&\W\B (A')\ot
\ar[rd]^(.45){\W(h^\dagger \T h^\dagger)}
\ar[dd]_(.3){\W(\l_A\T\l_A)\!}|\hole 
	\\
	&\W(\B PX)\ot 
\ar[uu]_(.325){\W(\B\pi_0)\ot} \ar[dd]^(.325){\W(\B\pi_1)\ot} 
\ar[rrr]
	&&&	\W\B(PX)\ot
\ar[uu]_{\W\B(\pi_0\ot)}
\ar[dd]^{\W\B(\pi_1\ot)}
\\
	\W(\B A)\ot
\ar[dr]_{\W(\xi\ot)\ } 
\ar[rrr]|!{[dr];[ur]}\hole 
	&&&\W\B (A\ot)
\ar[rd]^(.45){\W(\xi \T \xi)} \\
	&\W(\B X)\ot
\ar[rrr]
	&&&
	\W\B(X\ot)
	\mathrlap.
}
\]
The left face is the right face of the $\g$ square of \Cref{sec:gamma},
and the top, bottom, front, and back, 
each containing two edges $\W\mnn\nabla$,
commute because $\nabla({-}\, \ox \,{-}) = ({-}\, \T \,{-})\nabla$
by \Cref{thm:tensor-psi}. 

The right face does not necessarily commute as stands.
The issue is the limited functoriality of $\T$ in \DGC maps:
we have $\B(\pi_0 \ox \pi_0) (h^\dagger \T h^\dagger) = 
\B \pi_0 \,h^\dagger \T \B\pi_{0}\,h^\dagger = \l_{X}\xi' \T \l_{X}\xi'$,
and similarly for $\pi_{1}$, but no guarantee that 
$ \l_{X}\xi' \T \l_{X}\xi'$ should equal $(\l_{X} \T \l_{X})(\xi' \T \xi')$.
However, because the other five faces of each cube commute, 
when one prepends 
$\W\mnn\nabla\:\W(\B A' \ox \B A') \lt \W\B(A' \ox A')$
to the composites of the two maps around any square 
of the right face, 
the resulting maps \emph{are} equal.
Hence we replace the three copies of 
$\W\B(A' \ox A')$ along the back right edge
with $\W(\B A' \ox \B A')$ and get commutative cubes.
The right face we build off of in subsequent diagrams then becomes

\begin{equation}\label{eq:nabla-right-face}
\begin{aligned}
\xymatrix@R=.25em@C=6em{
	\W(\B A'\ox \B A')
\ar@{=}[dd]  
\ar[rd]^{\ \W(\xi' \T \xi')\W\mnn\nabla}
	\\&\W\B(X'\ox X') 
\ar[dd]^{\W(\l_{\mn X} \T \l_{\mn X})}
	\\\W(\B A' \ox \B A')
\ar@{=}[dd] \ar[rd] 
	\\
	&\W\B(X\ox X)
	\\\W(\B A' \ox \B A')
\ar[rd]^(.45){\ \W(h^\dagger \T h^\dagger)\W\mnn\nabla}
\ar[dd]_{\W(\l_A \ox \l_A)}
	\\
	&\W\B(PX \ox PX)
\ar[uu]_{\W\B(\pi_0\ox\pi_0)}
\ar[dd]^{\W\B(\pi_1\ox\pi_1)}
\\\W(\B A \ox \B A)
\ar[rd]^(.45){\ \W(\xi \T \xi)\W\mnn\nabla} \\
	&
	\W\B(X \ox X)
	\mathrlap.
\\
}
\end{aligned}
\end{equation}


\subsection{Repackaging the homotopy}\label{sec:t-tilde}

The format we require for representatives of homotopies 
in the upcoming \Cref{sec:Phi} 
has only one $P$, whereas $\W\B(PX \ox PX)$ has two,
so we need to reformat this map to fit into the diagram to come.
Thus we will convert the right homotopy
$\smash{\defm{\ul {f}}}\ceq 
\W(\hd \T \hd)\o\W\nabla\: \W(\B A' \ox \B A') \lt \W\B(PX \ox PX)$
and its associated endpoint maps
$\W\B(\pi_{\mn \smash j}\ot) \o \ul f\:
\W(\B A' \ox \B A') \lt \W\B(X \ox X)$
into a standard right homotopy 
in a number of steps.

\bitem

\item
Recall the natural map $r\: \big(P(-)\big){}\ot \lt D\big((-)\ot\big)$
of \Cref{thm:II-D},
morally restricting a square of \DGA maps to two adjacent edges.
Postcomposing $\W\B r$ to $\smash{\ul f}$,
we obtain a right homotopy
$\W g
\: \W(\B A')\ot \lt  \W\B D(X\ot)$,
where
$\defm g = \B r \o (\hd \T \hd) \o \nabla$. 
By \Cref{thm:II-D}, the new endpoint maps agree with the old:
as $p\mn_j \o r = \pi_j\ot$,
we have
$\W\B p\mn_{ j} \o \W g = 
\W\B(\pi_{\mn j}\ot) \o \ul f$.

\item
Recall also the composition operation
$D(X\ot) \os \e\from \W\B D(X\ot) \os{\mn\Y\, }\lt P(X\ot)$ 
of \Cref{def:Y}.
We want to attach our existing homotopy representative
to $\W\B$ of this operation, to wit,
\[
\W\B D(X\ot) \ \xleftarrow{\W\B\e} 
\W\B\W\B D(X\ot) \xtoo{\W\B\Y} 
\W\B P(X\ot)\mathrlap.
\]
In order to accomplish this, we need 
$g\: \W\B(A')\ot \lt  \W\B D(X\ot)$
to factor through $\W\B\e$.
It indeed does, 
by \Cref{thm:W-bar-factor},
and we have the following diagram:
\[
\xymatrix@C=3em@R=4em{
	&\W\B\W\B D(X\ot) \ar[r]^(.55){\W\B\Y}\ar[d]^{\W\B\e}
	& \W\B P (X\ot)\ar[d]^{\W\B\pi_{\mn j}}
	\\ 
	 \W(\B A')\ot \ar[r]_{\W g} \ar[ur]^{
	 	\defm{\bar f} \,\ceq\, \W\B \W g\, \W\h\ \,}
		&\W\B D(X\ot) \ar[r]_(.55){\W\B p_{\mn j}} 
	& \W\B(X\ot)\mathrlap,
	}
\]
where the triangle is this factorization.
The square commutes
by the definition of $\Upsilon$ in \Cref{def:Y}
and functoriality of $\W\B$,
Letting $\defm{(H^P)^\#}$ denote the composite along the top,
$\W\B(\Y\o\W g) \o \W\h$, 
we have preservation of endpoint maps:
$\W\B\pi_{\mn j} \o (H^P)^\# = \W\B p_{\mn j} \o \W g$
for $j \in \{0,1\}$.

\medskip 

\item
Now to free $P$ from $\W\B$, 
we postcompose the natural map $Z\: \W\B P(X\ot) \lt P \W \B(X\ot)$ 
of \Cref{thm:P-W}
to get a map $\defm{\wt H^P} \ceq Z \o (H^P)^\#$
satisfying 
$\pi_{\mn j} \wt H^P = \W\B\pi_{\mn j} \o (H^P)^\#$ for $j \in \{0,1\}$.
\eitem

All told, we will be able to glue the left of the following diagram
to the lower two squares on the existing 
right face \eqref{eq:nabla-right-face} the $\nabla$ diagram:

\[
\resizebox{146mm}{!}{
\xymatrix@R=1.5em@C=-.0625em{
	&&&&&	\W(\B A')\ot
	\ar@{=}[ddd]|!{[llllld];[dr]}\hole
	\ar[dr]
	\ar[llllld]|(.45){\W(\l_{\mn X} \T \l_{\mn X})\W(\xi' \T \xi')\W\mnn\nabla										{\phantom{catcatcatcat}}
				}
	\ar[drrrrrrr]
	\ar[drrrrrrrrrrrrr]
	\ar@/^/[drrrrrrrrrrrrrrrrrrr]
																	\\
			\W\B(X\ot)	
	\ar@{=}[rrrrrr]
	&&&&&	&\W\B(X\ot)	
	\ar@{=}[rrrrrr]
	&&&&&	&\W\B(X\ot)	
	\ar@{=}[rrrrrr]
	&&&&&	&\W\B(X\ot)	
	\ar@{=}[rrrrrr]
	&&&&&	&\W\B(X\ot)											\\	\\
	&&&&&	\W(\B A')\ot	
	\ar[drrrrrrr]|!{[ruu];[rd]}\hole
						|(.625){\bar f}
	\ar[drrrrrrrrrrrrr]	
						|!{[ruu];[rd]}\hole
						|(.49)\hole
						|(.5275){(H^P)^\#}
						|!{[ruurrrrrr];[rrrrrrrd]}\hole
	\ar@/^/[drrrrrrrrrrrrrrrrrrr]|!{[ruu];[rd]}\hole
						|!{[ruurrrrrr];[rrrrrrrd]}\hole
						|!{[ruurrrrrrrrrrrr];[rrrrrrrrrrrrrd]}\hole
						^(.605){\wt H^P}
	\ar[ddd]|!{[llllld];[dr]}\hole_(.75){\W(\l_A \ox \l_A)}
	\ar[dlllll]_{\ul f}
	\ar[dr]_(.45){\W g}													&\\
			\W\B(P X)\ot 
	\ar[uuu]^{\mathllap{\W\B(\pi_0\2)}}	
	\ar[ddd]_{\mathllap{\W\B(\pi_1\2)}}
	\ar[rrrrrr]_(.425){\W\B r}	
	&&&&&	&\W\B D(X\ot)	
	\ar[uuu]_{\W\B p_0}	
	\ar[ddd]^{\W\B p_1}
	&&&&&	&\W\B\W\B D(X\ot)
	\ar[llllll]^(.52){\W\B\e}
	\ar[rrrrrr]_(.5625){\W\B\Y}
	\ar[uuu]_{\W\B p_0\,\W\B\e}	
	\ar[ddd]^{\W\B p_1\,\W\B\e}
	&&&&& 	&\W\B P(X\ot)
	\ar[uuu]_{\W\B\pi_0}	
	\ar[ddd]^{\W\B\pi_1}
	\ar[rrrrrr]_{Z}
	&&&&&	&P\,\mn\W\B(X\ot)
	\ar[uuu]_{\pi_0}	
	\ar[ddd]^{\pi_1}
															\\		\\
	&&&&&	
	\W(\B A)\ot	
	\ar[llllld]_(.5){\W(\xi \T \xi)\W\mnn\nabla}	
	\ar[dr]
	\ar[drrrrrrr]|!{[ruu];[rd]}\hole
	\ar[drrrrrrrrrrrrr]|!{[ruu];[rd]}\hole
	|!{[ruurrrrrr];[rrrrrrrd]}\hole
	\ar@/^/[drrrrrrrrrrrrrrrrrrr]|!{[ruu];[rd]}\hole
	|!{[ruurrrrrr];[rrrrrrrd]}\hole
	|!{[ruurrrrrrrrrrrr];[rrrrrrrrrrrrrd]}\hole
																\\
	\W\B(X\ot)	\ar@{=}[rrrrrr]						
	&&&&&	&\W\B(X\ot)
	\ar@{=}[rrrrrr]
	&&&&&	&\W\B(X\ot)
	\ar@{=}[rrrrrr]
	&&&&&	&\W\B(X\ot)
	\ar@{=}[rrrrrr]
	&&&&&	&\W\B(X\ot)
	\mathrlap{.} 								
}
}
\]

\bs

\nd 
It happens that all of the horizontal maps are quasi-isomorphisms,
so the diagrams induce isomorphisms in Tor,
although strictly speaking we only need to know this for 
the backward-facing $\W\B\e$.
In terms of the Tor diagram \eqref{eq:nabla-Tor-squares}
of the $\W\nabla$ square, appending the triangular prism replaces
$\Tor_{\W(\B A')\ot}\big(\W\B(\textcolor{red}P X)\ot\big)$
with 
$
\smash{\Tor_{\W(\B A')\ot}\big(\textcolor{red}P\,\W \B (X\ot) \big)}
$
at the right of the third row and otherwise leaves the diagram unchanged.

\subsection{The \texorpdfstring{$\Phi$}{Phi} square}\label{sec:Phi}

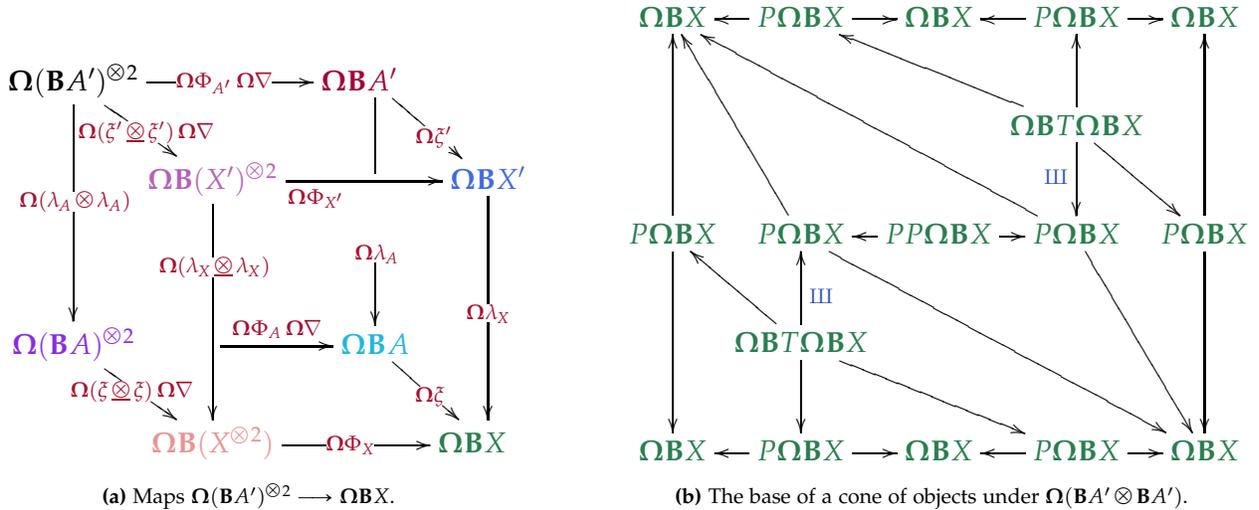
\begin{figure}
\begin{subfigure}{0.3\textwidth}
\centering
	$
	\begin{aligned}
		\xymatrix@R=1.75em@C=-.25em{	
			\W(\B A')\ot
			\ar@[jred][rrrrrr]|{\textcolor{jlabel}{\W\Phi_{A'} \, \W\mnn\nabla}}
			\ar@[jred][dr]|(.5){\textcolor{jlabel}{{\phantom{i}}\W(\xi' \T \xi') \, \W\mnn\nabla}}
			\ar@[jred][ddd]|(.45){\textcolor{jlabel}{\W(\l_{\mn A} \ox \l_{\mn A})\!}}
			&&&&&&	
		 \textcolor{jQ2}{	\W\B A'	}	{\phantom{{}\2}}
			\ar@[jred][dr]|(.55){\textcolor{jlabel}{\!\!\W\xi'}}
			\ar@[jred][ddd]|!{[dlllll];[dr]}\hole|(.65){\textcolor{jlabel}{\W \l_{\mn A}\!}}
			&\\
			&
		 \textcolor{compromise}{\W\B(X')\ot}
			\ar@[jred][rrrrrr]_(.375){\textcolor{jlabel}{\W\Phi_{\smash{X'}}}}
			\ar@[jred][ddd]|(.325){\textcolor{jlabel}{\W(\l_{\mn X} \T \l_{\mn X})\!}	}
			&&&&&	&
			 \textcolor{RoyalBlue}{\W\B X'}
			\ar@[jred][ddd]|{\textcolor{jlabel}{\W\l_{\mn X}}}
			\\&&&&&&	
			\\
			 \textcolor{BlueViolet}{\W(\B A)\ot}
			\ar@[jred][rrrrrr]|!{[uur];[dr]}\hole
							^(.675){\textcolor{jlabel}{\W\Phi_A\, \W\mnn\nabla}}
			\ar@[jred][dr]|(.45){\textcolor{jlabel}{\W(\xi \T \xi)\, \W\mnn\nabla\ }}
			&&&&&&	
		 \textcolor{jcyan}{	\W\B A}
			\ar@[jred][dr]|{\textcolor{jlabel}{\!\W\xi}}
			&\\
			&
		 \textcolor{jX2p}{	\W\B(X\ot)}
			\ar@[jred][rrrrrr]|{\textcolor{jlabel}{\W\Phi_X}}
			&&&&&&			
	 \textcolor{jgreen}	{	\W\B X}{\phantom{{}\2}} 				
		}
	\end{aligned}
	$
	\caption{Maps $\W(\B A')\ot \lt \W\B X$.}
	\label{fig:homotopy-cube}
\end{subfigure}%
\hfill%
\begin{subfigure}{0.6\textwidth}
\centering
$
\xymatrix@C=-.25em@R=2.5em{
	\textcolor{jgreen}{\W\B X} &&
	\textcolor{jgreen}{P\,\mn\W\B X}\ar@[jgreen][rr]\ar@[jgreen][ll]&&
	\textcolor{jgreen}{\W\B X}&&
	\textcolor{jgreen}{P\,\mn\W\B X} \ar@[jgreen][ll]\ar@[jgreen][rr]&&
	\textcolor{jgreen}{\W\B X}\\
	&&&&&&
	\textcolor{jgreen}{\W\B T\W\B X} \ar@[jgreen][llllu]\ar@[jgreen][u]
	\ar@[jgreen][d]_(.45){\textcolor{jblue}{\Sha}}\ar@[jgreen][rrd]\\
	\textcolor{jgreen}{P\,\mn\W\B X}  \ar@[jgreen][dd]\ar@[jgreen][uu]&&
	\textcolor{jgreen}{P\,\mn\W\B X} \ar@[jgreen][uull]\ar@[jgreen][ddrrrrrr]&& 
	\textcolor{jgreen}{P \mnn P\,\mn\W\B X}\ar@[jgreen][rr]\ar@[jgreen][ll]&&
	\textcolor{jgreen}{P\,\mn\W\B X} \ar@[jgreen][uullllll]\ar@[jgreen][rrdd]&& 
	\textcolor{jgreen}{P\,\mn\W\B X} \ar@[jgreen][dd]\ar@[jgreen][uu]\\
	&&
	\textcolor{jgreen}{\W\B T\W\B X} \ar@[jgreen][u]_(.45){\textcolor{jblue}{\Sha}}
	\ar@[jgreen][d]\ar@[jgreen][ull]\ar@[jgreen][drrrr]\\
	\textcolor{jgreen}{\W\B X} && 
	\textcolor{jgreen}{P\,\mn\W\B X} \ar@[jgreen][rr]\ar@[jgreen][ll]&&
	\textcolor{jgreen}{\W\B X} && 
	\textcolor{jgreen}{P\,\mn\W\B X} \ar@[jgreen][ll]\ar@[jgreen][rr]&&
	\textcolor{jgreen}{\W\B X}
}
$
\caption{The base of a cone of objects under $\W(\B A' \ox \B A')$.}
\label{fig:ring-map-path-plug}
\end{subfigure}
\caption{Auxiliary diagrams for the functoriality 
argument.}
	\label{fig:func-extra-figs}
\end{figure}

In filling in the $\Phi$ square of \eqref{eq:grand-scheme}, 
we are subject to a few constraints.
The isomorphisms of the previous large squares essentially 
carry $\Miso \ox \Miso$ along to the left edge,
and the top, bottom, and right of the square must
respectively describe $\Pi'$, $\Pi$, and $\Miso$.
\begin{figure}
\centering
	$
	\begin{aligned}
\xymatrix@C=1.75em@R=2em{
\smash{	\us{{\W(\B A')\ot}}
		\Tor}
	\big(\W\B(X')\ot\Ysub\big)				
											\ar[d]\ar[r]
	&\us{{\W}(\B A')\ot}
		\Tor(\W\B X') 
	&\us{\W(\B A')\ot}\Tor(P\,\mn\W\B X') \ar[l]_\sim \ar[r]
	&\smash{\us{\W\B A'}\Tor}
		\big(\W\B(X')	\Ysub\big)
										\ar[d]
	\\
	\smash{\us{\W(\B A')\ot}\Tor}
		\big(\W\B(X\ot)\Ysub\big)	
	&
	&
	&
	\smash{\us{\W\B A'}\Tor}
		(\W\B X\Ysub)
	\\
	\smash{\us{\W(\B A')\ot}\Tor}
		\big(P\,\mn\W\B (X\ot)\Ysub\big)
													\ar[u]^(.45)\vertsim
													\ar[d]
	&
	&		
	&\smash{\us{\W\B A'}\Tor}
		(P\,\mn\W\B X\Ysub)							\ar[u]_(.45)\vertsim
													\ar[d]
	\\
	\us{\W(\B A\ot)}\Tor\big(\W\B(X\ot)\Ysub\big)		\ar[r]
	&\us{\W(\B A)\ot}\Tor(\W\B X) 
	&\us{\W(\B A)\ot}\Tor(P\,\mn\W\B X) \ar[l]_\sim \ar[r]
	&\us{\W\B A}\Tor(\W\B X\Ysub)\mathrlap.
	}
	\end{aligned}
	$
	\caption{The constraints on filling the $\Phi$ square.}
	\label{fig:ring-map-Tor-ring}
\end{figure}
\nd
Thus our pre-existing commitments amount to \Cref{fig:ring-map-Tor-ring}.
It is not yet obvious this should commute, 
but we will fill it in
in such a way as to make commutativity apparent.
%
In terms of general strategy, counting $P$'s,
it is visible that we already have four homotopies of maps present,
to be accounted for as follows:
\bitem
\item 
Begin with the homotopies assumed in \Cref{thm:main-algebraic}
witnessing that $\xi',\ \xi,\ \upsilon,\ \upsilon'$
are \SHCA-maps
and 
precompose the two $\W\nabla$
maps $\W(\B A' \ox \B A') \lt \W\B(A'\ox A')$
and $ \W(\B A \ox \B A) \lt \W\B(A \ox A)$.
The maps along the top and bottom of \Cref{fig:ring-map-Tor-ring}
then follow from the associated six-square diagrams
as in \Cref{thm:homotopy-Tor-map}
\item 
The maps on the left are inherited from the previous squares of
\eqref{eq:grand-scheme}.
\item
The maps on the right come from the six-square diagram that 
\Cref{thm:homotopy-Tor-map} associates to the homotopy-commutativity
of the squares of \WHCA maps in the statement of
\Cref{thm:main-algebraic}.
\eitem
The two remaining homotopies 
are those we have assumed to
make $\l_A$, $\l_X$, $\l_Y$ \SHCA maps.
The six of these on the $A$-$X$ side 
(the $A$-$Y$ argument as usual proceeds silently in parallel)
together make up the cube of Figure \ref{fig:homotopy-cube};
note that we have precomposed
$\W\nabla$ 
so that $\W(\B A' \ox \B A')$ and $\W(\B A \ox \B A)$
rather than $\W\B(A' \ox A')$ and  $\W\B(A \ox A)$ appear.

%

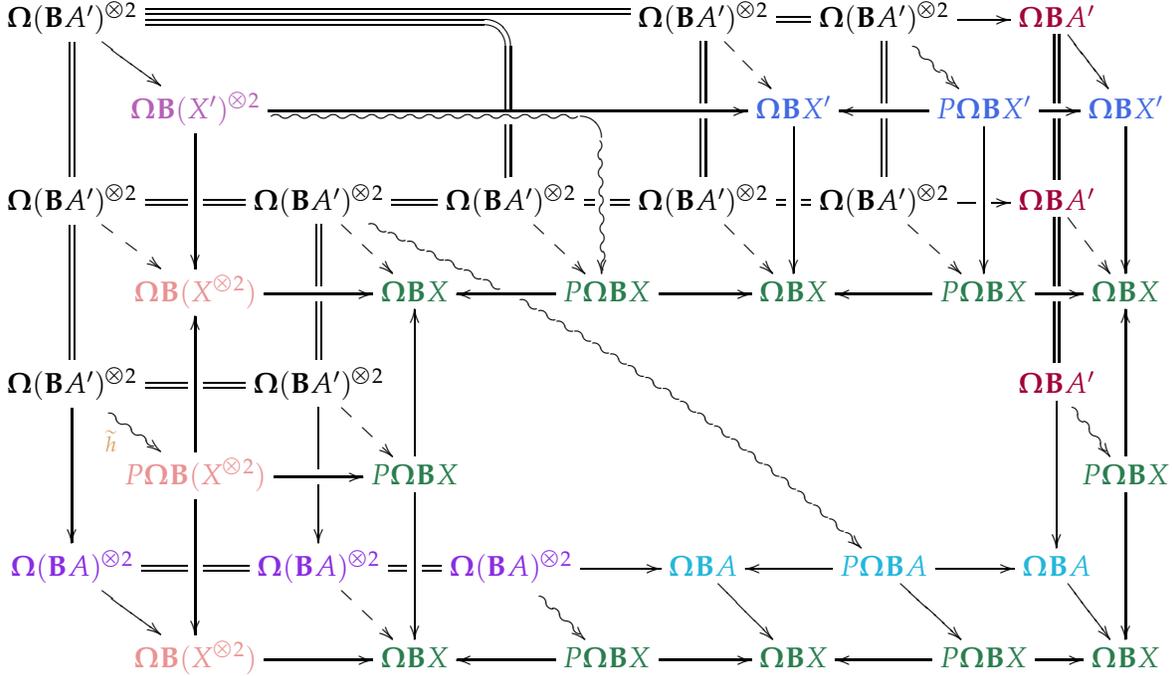
\begin{figure}
\centering
	$
	\begin{aligned}
	\xymatrix@C=-1em@R=1.5em{
		\W(\B A')\ot 	\ar@[jred][dr]
		\ar@<3pt>@{=}[rrrrrr]
		\ar@{=}[dd]
		\ar@{-}`r[rrrr][rrrrdd]						|(.45)\hole
		\ar@{-}@<-2pt>`r[rrrr][rrrrdd]		|(.45)\hole				
		&& &&
		&& \W(\B A')\ot \ar@{=}[rr]
		\ar@{-->}@[jcmp][dr]
		\ar@{=}[dd]|\hole
		&& \W(\B A')\ot \ar@[jred][rr]
		\ar@{=}[dd]|\hole
		&& \textcolor{jQ2}{\W\B A'} 		\ar@[jred][dr]
		\ar@{=}@[jQ2][dd]|\hole
		\\
		& \textcolor{compromise}{\W\B(X')\ot}	\ar@[jred][dd]
		\ar@[jred][rrrrrr]
		\ar@{~>}@[jhmt]@<-2pt>`r[rrrr][rrrrdd]						
		&& && 
		&&\textcolor{RoyalBlue}{ \W\B X' 	}	\ar@[jred][dd]
		&&\textcolor{RoyalBlue}{P\,\mn\W\B X' }	\ar@{<~}@[jhmt][ul]
		\ar@[RoyalBlue][ll]
		\ar@[RoyalBlue][rr]
		\ar@[jred][dd]
		&& \textcolor{RoyalBlue}{\W\B X'	}	\ar@[jred][dd]
		\\
		\W(\B A')\ot 		\ar@{=}[dd]
		\ar@{=}[rr]|\hole
		\ar@{-->}@[jcmp][dr]
		&&\W(\B A')\ot 	\ar@{=}[rr] 
		\ar@{=}[dd]|\hole
		\ar@{~>}@[jhmt]@/^1pc/[ddddrrrrrr]|(.31)\hole
		\ar@{-->}@[jcmp][dr]
		&&\W(\B A')\ot	\ar@{=}[rr]|(.475)\hole
		&&\W(\B A')\ot	\ar@{=}[rr]|\hole
		\ar@{-->}@[jcmp][dr]
		&&\W(\B A')\ot	\ar@{-->}@[jcmp][dr]
		\ar@[jred][rr]|(.575)\hole
		&&\textcolor{jQ2}{\W\B A'}		\ar@{=}@[jQ2][dd]|\hole
		\ar@{-->}@[jcmp][dr]
		\\
		&\textcolor{jX2p}{\W\B(X\ot)	}	\ar@[jred][rr]
		&&\textcolor{jgreen}{\W\B X}
		&&\textcolor{jgreen}{P\,\mn\W\B X} 		\ar@{<--}@[jcmp][ul]
		\ar@[jgreen][ll]
		\ar@[jgreen][rr]
		&&\textcolor{jgreen}{\W\B X} 
		&&\textcolor{jgreen}{P\,\mn\W\B X} 		
		\ar@[jgreen][ll]
		\ar@[jgreen][rr]
		&&\textcolor{jgreen}{\W\B X}
		\\
		\W(\B A')\ot 	\ar@[jred][dd]
		\ar@{=}[rr]|\hole
		&&\W(\B A')\ot 	\ar@[jred][dd]|\hole
		&& 
		&& 
		&& 
		&&\textcolor{jQ2}{\W\B A'}		\ar@[jred][dd]
		\\
		&\textcolor{jX2p}{P\,\mn\W\B(X\ot)}	\ar@{<~}@[jhmt][ul]^(.56){\textcolor{jhmt}{\wt H^P}}
		\ar@[jX2p][uu]
		\ar@[jX2p][dd]
		\ar@[jred][rr]
		&&\textcolor{jgreen}{P\,\mn\W\B X} 		\ar@{<--}@[jcmp][ul]
		\ar@[jgreen][uu]
		\ar@[jgreen][dd]
		&& && && 
		&& \textcolor{jgreen}{P\,\mn\W\B X}		\ar@{<~}@[jhmt][ul]
		\ar@[jgreen][uu]
		\ar@[jgreen][dd]
		\\
		\textcolor{BlueViolet}{\W(\B A)\ot}		\ar@[jred][dr]
		\ar@{=}@[BlueViolet][rr]|\hole
		&&\textcolor{BlueViolet}{\W(\B A)\ot}	
		\ar@{=}@[BlueViolet][rr]|\hole
		\ar@{-->}@[jcmp][dr]
		&& \textcolor{BlueViolet}{\W(\B A)\ot} 	\ar@[jred][rr]
		&& \textcolor{jcyan}{\W\B A }	\ar@[jred][dr]
		&& \textcolor{jcyan}{P\,\mn\W\B A }
		\ar@[jcyan][ll]
		\ar@[jcyan][rr]
		&& \textcolor{jcyan}{\W\B A	}	\ar@[jred][dr]
		\\
		&\textcolor{jX2p}{\W\B(X\ot)	}	\ar@[jred][rr]
		&&\textcolor{jgreen}{\W\B X} 
		&& \textcolor{jgreen}{P\,\mn\W\B X} 		\ar@{<~}@[jhmt][ul]
		\ar@[jgreen][ll]
		\ar@[jgreen][rr]
		&& \textcolor{jgreen}{\W\B X} 
		&& \textcolor{jgreen}{P\,\mn\W\B X}		\ar@{<-}@[jred][ul]
		\ar@[jgreen][ll]
		\ar@[jgreen][rr]
		&& \textcolor{jgreen}{\W\B X}
	}
	\end{aligned}$
	\caption{The assemblage of right homotopies implied by \Cref{fig:homotopy-cube}.}
	\label{fig:Phi-DGA-partial}
\end{figure}

The right homotopies witnessing these together fit into 
Figure \ref{fig:Phi-DGA-partial}.
We have color-coded the \DGAs by quasi-isomorphism type to match
 \eqref{fig:homotopy-cube}
 and colored the arrows coming from \eqref{fig:homotopy-cube} in red; 
we do not need to label them because they are uniquely determined
by their source and target.
Gold wavy 
arrows are right homotopies corresponding to the faces in \eqref{fig:homotopy-cube} 
and grey dashed arrows are the defined as the necessary composites
making the diagram commutative.
The projections from path objects are green, and arranged so that
$\pi_0$ always points up or left,
$\pi_1$ down or right.
The reader should convince themself
Figure \ref{fig:Phi-DGA-partial}
expresses only the existence of right homotopies
representing the homotopies we have just discussed.
We are not yet asserting anything about the front or back
of the large prism on the lower right.

Again by \Cref{thm:homotopy-compose}, 
the homotopies from $\W(\BA')\ot$
can be composed, and by \Cref{def:Sha}
the composite of two consecutive triples can 
be represented by a single right homotopy.
By \Cref{thm:homotopies-homotopic},
these composite right homotopies $\W(\BA')\ot \lt P\,\mn\W\B X$
are themselves homotopic,
and this is witnessed by a right homotopy $\W(\BA')\ot \lt P \mnn P\,\mn\W\B X$.
We can combine all the codomains into 
Figure \ref{fig:ring-map-path-plug}, 
to be thought of as a cone under $\W(\BA')\ot$.
\new{
As with 
\Cref{fig:assoc-path-plug},
\Cref{fig:comm-plug}, and 
\eqref{eq:diamonds}, 
this cone is not commutative, though by assumptions the faces apart from the base are.
Again, the issue can be repaired with the fix of
\Cref{thm:homotopy-independence-mk2},
at the cost of assuming
the right homotopy 
$H\: \W\B(A')\ot \lt PP\W\B X$
in the middle of the cone be {endpoint-fixing}.
In this case,
that means the maps
$P\pi_0 \o H$ and $P\pi_1 \o H\:\W\B(A\ott) \lt P\W\B X$
must factor respectively as
$\z\o\W\Phi_{X}\o\W(\l_X\T\l_X)\o\W(\xi'\T\xi')\o\W\nabla$
and
$\z\o\W\xi\o\W\l_A\o\W\Phi_{A'}\o\W\nabla$,
where $\z\: \W\B X \lt P \W\B X$ 
is the natural map 
of \Cref{def:P}.
Thus the relevant cone from $\W(\BA')\ot$ 
is again in fact over \Cref{fig:true-base}.
}

\new{
As before, the author does not know
when this is actually achievable.
}

\begin{figure}
\centering
$
\resizebox{146mm}{!}{
	\xymatrix@C=1.25em@R=4em{
		\smash{\us{\mathclap{\W(\B A')\ot}}\Tor}\ \,\,\, \big(\W\B (X')\ot\big)
				\ar@[nonisored][rrr]
				\ar@[nonisored][d]
				\ar@/^2pc/@<-.75ex>@[nonisored][rrd]
		&
		&
		&\smash{\us{\mathclap{\W(\B A')\ot}}\Tor}\ \, (\W\B X')\ar@[nonisored][d]
		&\smash{\us{\mathclap{\W(\B A')\ot}}\Tor}\ \, (P\,\mn\W\B X')\ar@[nonisored][d]\ar@[nonisored][r]\ar[l]
		&\smash{\us{\mathclap{\W\B A'}}\Tor}\ (\W\B X')\ar@[nonisored][d]
	\\
		\smash{\us{\mathclap{\W(\B A')\ot}}\Tor}\ \,\,\, \big(\W\B (X\ot)\big)
			\ar@[nonisored][r]
		& {\us{\mathclap{\W(\B A')\ot}}\Tor}\ \, (\W\B X){\vphantom{X_{X_X}}}
		& \smash{\us{\mathclap{\W(\B A')\ot}}\Tor}\ \, (P\,\mn\W\B X)
			{\vphantom{X^{X^X}}}
			\ar[l]\ar[r]
		& \smash{\us{\mathclap{\W(\B A')\ot}}\Tor}\ \, (\W\B X)
		& \smash{\us{\mathclap{\W(\B A')\ot}}\Tor}\ \, (P\,\mn\W\B X)\ar[l]\ar@[nonisored][r]
		& {\us{\mathclap{\W\B A'}}\Tor}\ (\W\B X)
	\\
		&
		&
		&
		& \smash{\us{\mathclap{\W(\B A')\ot}}\Tor}\ \, (\W\B T\W\B X)
			\ar@[nonisored][ru]
			\ar[d]
			\ar@/^1.5pc/@[nonisored][rddd]
			\ar@[nonisored][rd]
			\ar[u]
			\ar[ul]
			\ar[ull]
			\ar@/^/[ulll]
	\\
		\smash{\us{\mathclap{\W(\B A')\ot}}\Tor}\ \,\,\, \big(P\,\mn\W\B (X\ot)\big)
			\ar@[nonisored][dd]
			\ar[uu]
		\ar@[nonisored][r]
		&\smash{\us{\mathclap{\W(\B A')\ot}}\Tor}\ \,  (P\,\mn\W\B X){\vphantom{X_{X_X}}}
			\ar@[nonisored][dd]
			\ar[uu]
		& \smash{\us{\mathclap{\W(\B A')\ot}}\Tor}\ \, (P\,\mn\W\B X)
			\ar[luu]
			\ar@[nonisored][rrrdd]
		&\smash{\us{\mathclap{\W(\B A')\ot}}\Tor}\ \, (P \mnn P\,\mn\W\B X)
			\ar[r]
			\ar[l]
			\ar[lluu]
			\ar@[nonisored][rrdd]
		&\smash{\us{\mathclap{\W(\B A')\ot}}\Tor}\ \, (P\,\mn\W\B X)
			\ar[uulll]
			\ar@[nonisored][rdd]	
		&\smash{\us{\mathclap{\W\B A'}}\Tor}\ (P\,\mn\W\B X)
			\ar@[nonisored][dd]
			\ar[uu]	
	\\
		&
		& \smash{\us{\mathclap{\W(\B A')\ot}}\Tor}\ \, (\W\B T\W\B X)
			\ar@[nonisored][ld]
			\ar@[nonisored][d]
			\ar@[nonisored][rd]
			\ar@[nonisored][rrd]
			\ar@/^/@[nonisored][rrrd]
			\ar[u]
			\ar@/^1.75pc/[luuu]
			\ar[ul]
	\\
		\smash{\us{\mathclap{\W(\B A)\ot}}\Tor}\ \,\, \big(\W\B (X\ot)\big)\ar[r]
		& \smash{\us{\mathclap{\W(\B A)\ot}}\Tor}\ \, (\W\B X)
		& \smash{\us{\mathclap{\W(\B A)\ot}}\Tor}\ \, (P\,\mn\W\B X)\ar[l]\ar[r]
		& \smash{\us{\mathclap{\W\B A}}\Tor}\, (\W\B X)
		& \smash{\us{\mathclap{P\,\mn\W\B A}}\Tor}\ (P\,\mn\W\B X)\ar[l]\ar[r]
		& {\us{\mathclap{\W\B A}}\Tor}\, (\W\B X)
		}
	}
$
\caption{The completed $\Phi$ square.}
\label{fig:Phi-Tor-partial}
\end{figure}
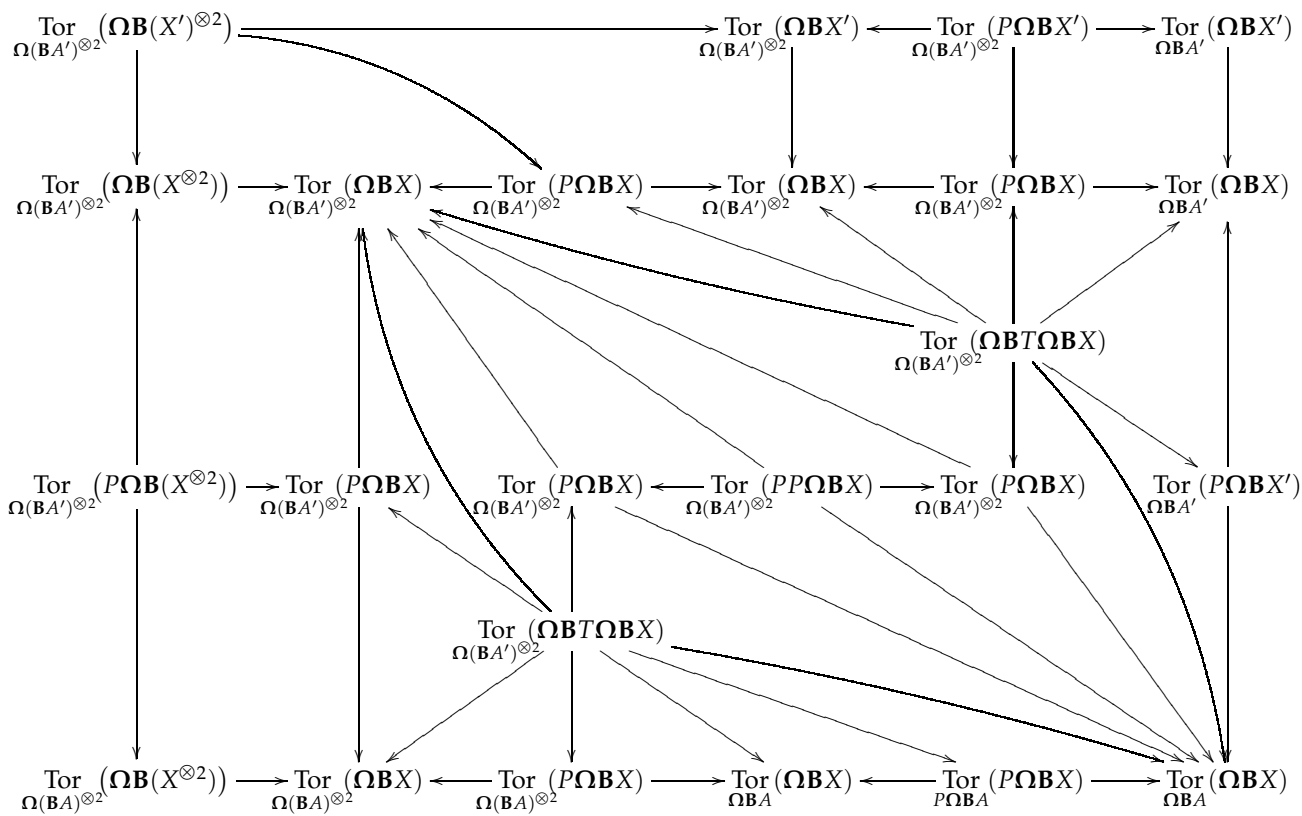

Using the factorizations of the maps along the right and bottom edges
through $\W\B(A\ot)$ and $\W\B A$,
we may insert this cone into Figure \ref{fig:Phi-DGA-partial},
and taking Tor, obtain Figure \ref{fig:Phi-Tor-partial},
in which the black arrows are isomorphisms and the red are not.
This is the $\Phi$ square of \eqref{eq:grand-scheme},
and using all subdividing commutative squares and triangles,
we see it commutes.

\bigskip

This completes the $\Phi$ diagram
and the proof $\Miso$ is multiplicative with respect to $\Pi'$ and $\Pi$.

\bs

\brmk\label{rmk:HGA}

It has long been known~\cite[\SS2]{baues1981double}\cite[Cor.~6]{gerstenhabervoronov1995} 
that the bar construction of the normalized cochain algebra 
$\B\C(X)$ of a connected simplicial set carries a differential graded 
Hopf algebra structure functorial in $X$
(also known as a homotopy Gerstenhaber algebra structure on $\C(X)$, 
and related~\cite[Rmk.~4.2]{franz2019shc} 
to the traditional \SHCA structure on $\C(X)$).
As a result, there is an idea for
an alternative construction of a product on Tor,
starting with the exterior product 
and $\gamma$ as we have done,
then following not with $\W\nabla$ and squares involving $\Phi$,
but $\W\mu$
for $\mu\mnn\: \B\C(B) \ox \B\C(B) \lt \B\C(B)$
the \DG Hopf algebra multiplication, 
which is a \DGC map.
As far as we can tell,
such a proof would require stronger hypotheses than the present one.
\ermk

\brmk[Sketches of other ways forward]\label{rmk:Cinfty} 
It has been pointed out to us that results in $E_n$-spectra
could likely also be used to prove an analogue to
\Cref{thm:Tor-CGA}, with some interpretation.
Basterra and Mandell show~\cite[Thm.~5.3]{basterramandell2011}
that the bar construction $B\wt A$ of a so-called 
augmented partial $\ms C_n$-algebra $\wt A$
is an augmented partial $\ms C_{n-1}$-algebra,
where $\ms C_n$ is the little $n$-cubes operad.
A variant of their proof likely establishes 
(but this should be checked)
that given
a span of partial $\ms C_n$-algebras $\wt X \from \wt A \to \wt Y$,
the two-sided bar construction $B(\wt X,\wt A,\wt Y)$
is also a partial $\ms C_{n-1}$-algebra.
A suitably enhanced version of the Dold--Kan correspondence
should take a span  $X \from A \to Y$ of $E_3$-algebras
to a span of augmented 
partial $\ms C_3$-algebras $\wt X \from \wt A \to \wt Y$,
so that $B_\bullet(\wt X,\wt A,\wt Y)$ becomes a partial
$\ms C_2$-algebra whose cohomology is $\Tor_A(X,Y)$,
and this should also give a \CGA structure on Tor functorial in
triples of $E_3$-algebra maps making the expected pair of squares commute.

There are minor attendant difficulties in formalizing this argument, 
which would involve generalizing the Basterra--Mandell proof rather than 
simply applying their result,
and the connection of true and partial $\ms C_n$-algebras is not direct,
but passes through a certain zigzag of equivalences.
An additional complication,
for our intended topological application in
\Cref{thm:main-topological},
is that
the existing \SHCA formality maps given by 
Munkholm inducing the additive 
isomorphism $\Tor_{\H(B)}\big(\H(X),\H(E)\big) 
\lt\Tor_{\C(B)}\big(\C(X),\C(E)\big)$ 
have not been shown to be $E_3$-algebra maps
(and as far as this author can see may not be),
and the squares only commute up to \Ai-homotopy anyway.
It may of course 
be that sufficient functoriality properties can be recovered
for this hypothetical other version of the product as well,
but this is not obvious.
A benefit of the existing argument 
leading to
\Cref{thm:main-topological}
is that the hypotheses seem to be minimal,
and the homotopy-commutative squares of \SHCA maps 
to be taken as input to the theorem 
are already known,
so that from our current \Cref{thm:Tor-CGA}, 
\Cref{thm:main-topological} is automatic, 
whereas with another approach,
as we have discussed, 
some additional massaging would be needed.

Another suggestion that has come to the author involves
a result of Fresse~\cite{Fressebar} 
that the bar construction of an $E_\infty$-algebra $A$
is another $E_\infty$-algebra. 
This likely applies as well to the two-sided
bar construction $\B(X,A,Y)$ of a span of $E_\infty$-algebras,
hence inducing a \CGA structure on its cohomology.
If so, this would 
induce the expected \CGA structures on 
$\Tor_{\C(B)}\big(\C(X),\C(E)\big)$
and 
$\Tor_{\H(B)}\big(\H(X),\H(E)\big)$
under mild flatness conditions.
Again, however, there is not a reason to expect this construction 
of a \CGA structure to be to be functorial under the sort of homotopy-commutative
squares of \Ai-algebra maps we already have,
and thus this approach as well
does not immediately yield our target application. 
\ermk

\counterwithin{figure}{section}
\counterwithin{theorem}{section}
\counterwithin{lemma}{section}
\counterwithin{corollary}{section}
\counterwithin{proposition}{section}
\counterwithin{equation}{section}
\counterwithin{notation}{section}







\bs

\newcommand{\etalchar}[1]{$^{#1}$}


{
\noindent {34 Evergreen St.}, \\
{Providence, RI 02906}, \\
{USA} } 

\smallskip

\nd\url{jdkcarlson@gmail.com}%

\end{document}